\documentclass[11pt]{article}
\usepackage[utf8]{inputenc}
\usepackage{geometry}
\usepackage{aligned-overset}
\usepackage{amsmath}
\usepackage{amsfonts}
\usepackage{amssymb}
\usepackage{amsthm}
\usepackage{graphicx}
\usepackage{enumitem,mdwlist}
\setlist[itemize]{topsep=0ex,itemsep=0ex,parsep=0.4ex}
\setlist[enumerate]{topsep=0ex,itemsep=0.1ex,parsep=0.4ex}
\usepackage{bm}
\usepackage{caption}
\usepackage{subcaption}
\usepackage{url}
\usepackage[linktoc = all, hidelinks, colorlinks, unicode=true, pagebackref=true]{hyperref} 
\hypersetup{linkcolor={blue}, citecolor={green!60!black}, urlcolor={black}}
\usepackage[table,xcdraw]{xcolor}
\usepackage{mathtools}
\usepackage{systeme}
\usepackage[capitalise]{cleveref}
\usepackage{thm-restate}
\usepackage{todonotes}
\usepackage{bbm}
\usepackage{comment}
\usetikzlibrary{shapes,arrows,calc}
\usepackage{setspace}
\usepackage[square,sort,comma,numbers]{natbib} 
\newtheorem{theorem}{Theorem}[section]
\newtheorem{proposition}[theorem]{Proposition}

\newtheorem{conjecture}[theorem]{Conjecture}
\newtheorem{lemma}[theorem]{Lemma}
\newtheorem{corollary}[theorem]{Corollary}
\newtheorem{fact}[theorem]{Fact}

\newtheorem*{remarks*}{Remarks}
\newtheorem*{remark*}{Remark}

\newtheorem{claim}{Claim}[theorem]
\crefname{equation}{}{}
\crefname{item}{}{}
\newcommand{\defn}[1]{\textcolor{red!65!black}{\emph{#1}}}
\numberwithin{equation}{section}
\newenvironment{proofclaim}[1][Proof of claim]{\begin{proof}[#1]}{\end{proof}}

\DeclarePairedDelimiter{\ceil}{\lceil}{\rceil}
\newcommand{\E}{\mathbb{E}}

\newcommand{\N}{\mathbb{N}}
\newcommand{\even}{\mathrm{Even}}
\newcommand{\odd}{\mathrm{Odd}}

\newcommand{\cP}{\mathcal{P}}
\newcommand*{\ext}{\textup{Ext}}
\newcommand*{\inter}{\textup{Int}}

\newcommand{\diam}{\mathrm{diam}}

\newcommand{\floor}[1]{\lfloor{#1}\rfloor}
\newcommand{\inv}{^{-1}}
\renewcommand{\geq}{\geqslant}
\renewcommand{\leq}{\leqslant}

\newcommand{\firstborn}[1]{\hat{#1}}
\Crefname{conjecture}{Conjecture}{Conjectures}
\Crefname{claim}{Claim}{Claims}
\Crefname{lemma}{Lemma}{Lemmas}
\Crefname{fact}{Fact}{Facts}
\Crefname{subsection}{Subsection}{Subsections}
\Crefname{figure}{Figure}{Figures}

\theoremstyle{definition}
\newtheorem{definition}[theorem]{Definition}

\def\namedlabel#1#2{\begingroup
    #2%
    \def\@currentlabel{#2}%
    \phantomsection\label{#1}\endgroup
}

\def\ew#1{}
	\renewcommand{\ew}[1]{\footnote{\textbf{EW:}#1}}
 
\def\eps{{\varepsilon}}
\title{Embedding trees using minimum and maximum degree conditions}

\author{
Alexey Pokrovskiy\footnotemark[1]
\and
Leo Versteegen\footnotemark[2]
\and 
Ella Williams\footnotemark[1]
}

\date{\today}

\begin{document}

\maketitle

\begin{abstract}
A variant of the Erd\H{o}s--S\'os conjecture, posed by Havet, Reed, Stein and Wood, states that every graph with minimum degree at least $\lfloor 2k/3 \rfloor$ and maximum degree at least $k$ contains a copy of every tree with $k$ edges. Both degree bounds are best possible. 

We confirm this conjecture for large trees with bounded maximum degree, by proving that for all~$\Delta\in \N$ and sufficiently large $k\in \N$, every graph $G$ with $\delta(G)\geq \lfloor2k/3 \rfloor$ and $\Delta(G)\geq k$ contains a copy of every tree $T$ with $k$ edges and $\Delta(T)\leq \Delta$. 

We also prove similar results where alternative degree conditions are considered. For the same class of trees, this verifies exactly a related conjecture of Besomi, Pavez-Sign\'e and Stein, and provides asymptotic confirmations of two others.
\end{abstract}

\renewcommand{\thefootnote}{\fnsymbol{footnote}}

\footnotetext[0]{\emph{MSC 2020}: 05C05, 05C07, 05C35}

\footnotetext[1]{Department of Mathematics, University College London, UK. Emails: \tt{\{\href{mailto:a.pokrovskiy@ucl.ac.uk}{a.pokrovskiy},\href{mailto:ella.williams.23@ucl.ac.uk}{ella.williams.23}\}@ucl.ac.uk}. \textrm{Research of EW supported by the Martingale Foundation.}}
\footnotetext[2]{Mathematics Institute, University of Warwick, UK. Email: \tt{\href{mailto:lversteegen.math@gmail.com}{lversteegen.math@gmail.com}}. \textrm{Research carried out while at Department of Mathematics, The London School of Economics, UK.}}

\section{Introduction}

Which degree conditions can be imposed on a host graph in order to guarantee that it contains a copy of every tree of a fixed size? A standard observation is that a graph of minimum degree at least $k$ contains a copy of every $k$-edge tree, as can be seen by greedily embedding the vertices of the tree one-by-one, using a degeneracy ordering. This bound cannot be lowered in general since there must exist a vertex of degree $k$ in order to embed the $k$-edge star, or alternatively, to avoid the components of the graph having fewer than $k+1$ vertices, in which case no $k$-edge tree can be embedded. Motivated by this observation, one can instead ask what happens when we consider alternative degree conditions within the host graph. A notable conjecture of Erd\H{o}s and S\'os from the 1960s (see \cite{Erdos63}) is a central point for much research in the area. It states
that every graph with average degree strictly larger than $k-1$ contains every $k$-edge tree. In the 1980s, Ajtai, Koml\'os, Simonovits and Szemer\'edi
announced a proof of this conjecture for large $k$. Although a publication of their result has not yet appeared, some of the main ideas have been communicated \cite{ES_slides}, and in the meantime several partial results have been obtained, the most general of which are an approximate proof for dense graphs and large trees with sublinear maximum degree due to Besomi, Pavez-Sign\'e and Stein \cite{BPS_ErdosSos} and independently Rozhoň \cite{rozhon19}, and for large trees with maximum degree bounded by a constant due to the first author \cite{AP_hyperstability}. A more detailed survey of progress towards this conjecture where specific cases have been considered can be found in \cite{AP_hyperstability}. 

As a natural adaptation of the Erd\H{o}s--S\'os conjecture, Havet, Reed, Stein and Wood \cite{HRSW02497401} considered a combination of minimum and maximum degree requirements for the host graph, and conjectured the following.

\begin{conjecture}[Havet, Reed, Stein and Wood {\cite{HRSW02497401}}]\label{conj:2k/3}
    Every graph with minimum degree at least $\lfloor 2k/3 \rfloor$ and maximum degree at least $k$ contains a copy of every $k$-edge tree.
\end{conjecture}

Both the minimum and maximum degree conditions in the conjecture are best possible, where the latter is necessary to again ensure there is a component in $G$ with at least $k+1$ vertices. To see that the minimum degree condition cannot be lowered, consider the following construction, assuming first that $k$ is divisible by $3$. Suppose $G$ is obtained by taking the union of two vertex-disjoint cliques, each with $2k/3-1$ vertices, and adding a universal vertex $x$. Observe that $\Delta(G) = 4k/3-2$ and $\delta(G)= 2k/3-1$. Let $T$ be a tree on $k$ edges containing a vertex $v$ with $d_T(v) = 3$ such that each of the three components of $T-v$ have exactly $k/3$ vertices. See \cref{fig:extremal2k/3}. No matter which vertex of $T$ we embed at $x$, all vertices in two of the components of $T-v$ must be embedded into the same clique in $G$, but neither clique is big enough. It is easy to generalise this example to demonstrate the sharpness of \cref{conj:2k/3} when $k$ is not divisible by $3$.  

\begin{figure}
\centering
            \begin{minipage}[b]{0.47\linewidth}
            \centering
\begin{tikzpicture}[scale=0.8]
    \draw [line width=0.5mm, draw=black, fill=black!5]
       (0,3) -- (2.97,0.74) -- (0.78,0) -- cycle;
    \draw [line width=0.5mm, draw=black, fill=black!5]
       (0,3) -- (-2.97,0.74) -- (-0.78,0) -- cycle;
\filldraw[black] (0,3) circle (4pt) node[anchor=south, yshift=0.2em]{$x$};
\draw[black, fill =black!5,line width=0.5mm] (-2,0) circle (35pt) node[anchor=south, yshift=-1em]{\large{$\frac{2k}{3} -1$}};
\filldraw[black, fill = black!5,line width=0.5mm] (2,0) circle (35pt) node[anchor=south, yshift=-1em]{\large{$\frac{2k}{3} -1$}};
            \end{tikzpicture}
        \end{minipage}
        \hspace{0.4cm}
        \begin{minipage}[b]{0.47\linewidth}
            \centering
             \begin{tikzpicture}[scale=1,rotate=-30]
\foreach \y in {1,2,3}{
\filldraw ({120*\y}:1) circle (2.5pt);
\draw[line width =  0.5mm] ({120*\y}:1) -- (0,0);
\draw[loosely dotted, line width =0.5mm] ({120*\y}:1) -- ({120*\y-22}:2);
\draw[loosely dotted, line width =0.5mm] ({120*\y}:1) -- ({120*\y+22}:2);}
\node[] at (1.8,0) {$k/3$};
\node[] at (-1.2,-1.9) {$k/3$};
\node[] at (-0.95,1.7) {$k/3$};
\filldraw[black] (0,0) circle (3pt) node[anchor=east, xshift=-0.1em, yshift=0.3em]{$v$};
            \end{tikzpicture}
 \end{minipage}
 \vspace{0.5cm}

 \begin{minipage}[b]{0.47\linewidth}
     \centering \small $G$ is the union of two cliques \\and a universal vertex
 \end{minipage}
 \hspace{0.3cm}
 \begin{minipage}[b]{0.47\linewidth}
    \centering \small Three components of $T-v$ \\ each have $k/3$ vertices
 \end{minipage}
 \caption{Extremal example for \cref{conj:2k/3}}\label{fig:extremal2k/3}
\end{figure}

Reed and Stein \cite{ReedStein_spanning1,ReedStein_spanning2} proved an exact version of \cref{conj:2k/3} for large spanning trees, that is, when the host graph has exactly $k+1$ vertices. Havet, Reed, Stein and Wood \cite{HRSW02497401} proved that there exists a function $f:\N \rightarrow \N$ such that $\delta(G)\geq \floor{2k/3}$ and $\Delta(G)\geq f(k)$ is sufficient for $G$ to contain a copy of every $k$-edge tree.
The bound on the maximum degree of $G$ was later improved for dense graphs by Besomi, Pavez-Sign\'{e} and Stein \cite{BPSdegreeconditions2019}, who proved an approximate version of the conjecture for large trees with maximum degree bounded by a root of $k$.

\begin{theorem}[Besomi, Pavez-Sign\'{e} and Stein {\cite{BPSdegreeconditions2019}}]\label{thm:dense2k/3}
For all $\varepsilon>0$ there exists $k_0$ such that for all $n$ and $k>k_0$ with $ k \geq \eps n$ the following holds.
Every graph on $n$ vertices with minimum degree at least $(2/3+\varepsilon)k$ and maximum degree at least $(1+\eps)k$
contains a copy of every $k$-edge tree $T$ with $\Delta(T)\leq k^{1/49}$.
\end{theorem}

It is natural to consider more combinations of minimum and maximum degree conditions. The same three authors posed the following conjecture.

\begin{conjecture}[Besomi, Pavez-Sign\'{e} and Stein {\cite{BPSmaxmin2020}}]\label{conj:alpha}
Let $k \in \mathbb{N}$ and let $\alpha \in (0,1/3)$.
Every graph with minimum degree at least $(1+\alpha)k/2$ and maximum degree at least $2(1-\alpha)k$ contains a copy of every $k$-edge tree.
\end{conjecture}

For every odd integer $\ell\geq 5$, Besomi, Pavez-Sign\'{e} and Stein \cite{BPSmaxmin2020} showed that the bounds are asymptotically best possible when $\alpha =1/\ell$ (see also \cref{sec:conclusion}). In favour of \cref{conj:alpha}, an approximate version has been proven, again in the dense setting and for trees with maximum degree bounded by a root of $k$.

\begin{theorem}[Besomi, Pavez-Sign\'{e} and Stein {\cite{BPSmaxmin2020}}]\label{thm:denseBPSalpha}
For all $\varepsilon>0$ there exists $k_0$ such that for all $n$ and $k>k_0$ with $k \geq \eps n$ and for each $\alpha \in (0,1/3)$ the following holds.
Every graph on $n$ vertices with minimum degree at least $(1+\eps)(1+\alpha)k/2$ and maximum degree at least $2(1+\eps)(1-\alpha)k$ contains a copy of every $k$-edge tree $T$ with $\Delta(T)\leq k^{1/67}$.
\end{theorem}

\cref{conj:alpha} was originally also posed for the case where $\alpha =0$, that is, when the host graph has minimum degree at least $k/2$ and maximum degree at least $2k$. However, a counterexample to this was given by Hyde and Reed \cite{HYDE2023}, who found a $k$-edge tree with maximum degree $(k-3)/2$ that cannot be embedded into a host graph with the given degree conditions. The host graph is constructed using a random graph. Conversely, in the same paper they also proved that minimum degree $k/2$ is sufficient for embedding all $k$-edge trees provided that the host graph has maximum degree at least $g(k)$ where $g$ is some large function of $k$. It was commented there that taking \cref{conj:alpha} with $\alpha = 0$ may be true if the inequalities are made strict in the statement.  Furthermore, since the tree in the counterexample has high maximum degree, one could ask what happens if we restrict the setting only to bounded degree trees. In this direction, Besomi, Pavez-Sign\'{e} and Stein suggested a stronger statement could hold for trees with maximum degree bounded by a constant. 

\begin{conjecture}[Besomi, Pavez-Sign\'{e} and Stein {\cite{BPSdegreeconditions2019}}]\label{conj:boundedk/2}
Let $k, \Delta \in \mathbb{N}$.
Every graph with minimum degree at least $k/2$ and maximum degree at least $2(1-1/\Delta)k$ contains a copy of every $k$-edge tree $T$ with $\Delta(T)\leq \Delta$.
\end{conjecture}

In \cite{BPSdegreeconditions2019}, an example is provided showing that the bounds in \cref{conj:boundedk/2} are close to best possible. In the same paper, as evidence towards the conjecture holding, an approximate version was given again in the dense setting.
\begin{theorem}[Besomi, Pavez-Sign\'{e} and Stein {\cite{BPSdegreeconditions2019}}]\label{thm:2k-1/Delta_dense}
For all $\varepsilon>0$ and $\Delta \geq 2$  there exists $k_0$ such that for all $n$ and $k>k_0$ with $k \geq \eps n$ the following holds. 
Every graph on $n$ vertices with minimum degree at least $(1+\eps)k/2$
and maximum
degree at least $2(1-1/\Delta+\eps)k$ contains a copy of every $k$-edge tree $T$ with $\Delta(T)\leq \Delta$.
\end{theorem}
One final conjecture that we consider, also due to Besomi, Pavez-Sign\'{e} and Stein, now replaces the maximum degree condition with a requirement for the existence of a vertex with many first and second neighbours. For a graph $G$ and a vertex $x\in V(G)$, the \defn{second neighbourhood} of $x$ in $G$, denoted by \defn{$N^2_G(x)$}, is the set of vertices $y\in V(G)\setminus \{x\}$ sharing a common neighbour with $x$.

\begin{conjecture}[Besomi, Pavez-Sign\'{e} and Stein {\cite{BPSdegreeconditions2019}}]\label{conj:secondneighbourhood}
Let $k\in \mathbb{N}$. If $G$ is a graph with $\delta(G)\geq k/2$ containing a vertex $x\in V(G)$ such that $\min\{|N_G(x)|,|N^2_G(x)|\}\geq 4k/3$, then $G$ contains a copy of every $k$-edge tree.
\end{conjecture}
An approximate version in the dense setting has been proven by the same authors, for trees with maximum degree bounded by a root of $k$.

\begin{theorem}[Besomi, Pavez-Sign\'{e} and Stein {\cite{BPSdegreeconditions2019}}]\label{thm:secondneighbourhood_dense}
 For all $\varepsilon \in(0,1)$ there exists $k_0$ such that for all $n$ and $k>k_0$ with $ k \geq \eps n$ the following holds. 
If $G$ is a graph on $n$ vertices with $\delta(G)\geq (1+\eps)k/2$ and containing a vertex $x\in V(G)$ such that $\min\{|N_G(x)|,|N_G^2(x)|\}\geq (1+\eps)4k/3$, then $G$ contains a copy of every $k$-edge tree $T$ with $\Delta(T)\leq k^{1/67}$.
\end{theorem}

\subsection{Our results: extensions to the sparse setting}
The main contribution of this paper to prove \cref{conj:2k/3} exactly for large bounded degree trees.
\begin{theorem}\label{thm:exact2/3}
    For all $\Delta \in \N$ and sufficiently large $k\in \N$ the following holds. If $G$ is a graph satisfying $\delta(G)\geq \lfloor2k/3\rfloor$ and $\Delta(G)\geq k$, then $G$ contains a copy of every $k$-edge tree $T$ with $\Delta(T)\leq \Delta$.
\end{theorem}
Regarding alternative combinations of minimum and maximum degree, our methods allow us to prove exactly \Cref{conj:alpha} for the same class of trees.
\begin{theorem}\label{thm:mainalpha-exact} 
       For all $ \Delta \in \N$, $\alpha \in (0,1/3)$ and sufficiently large $k\in \N$ the following holds. 
     If~$G$ is a graph satisfying  $\delta(G)\geq(1+\alpha)k/2$ and $\Delta(G)\geq2(1-\alpha)k$, 
   then $G$ contains a copy of every $k$-edge tree $T$ with $\Delta(T)\leq\Delta$.     
\end{theorem}

Finally, we asymptotically prove both \cref{conj:boundedk/2,conj:secondneighbourhood}, again for large bounded degree trees.
\begin{theorem}\label{thm:approx_k/2}
       For all $ \Delta \in \N$, $\varepsilon >0$ and sufficiently large $k\in \N$ the following holds. 
     If $G$ is a graph satisfying  $\delta(G)\geq(1+\varepsilon)k/2$ and $\Delta(G)\geq2(1-1/\Delta +\eps)k$, 
   then $G$ contains a copy of every $k$-edge tree~$T$ with $\Delta(T)\leq\Delta$. 
\end{theorem}

\begin{theorem}\label{thm:secondnbhd_result}
For all $ \Delta \in \N$, $\eps>0$ and sufficiently large $k\in \N$ the following holds. 
If $G$ is a graph satisfying $\delta(G)\geq (1+\eps)k/2$ and containing a vertex $x\in V(G)$ such that $\min\{|N_G(x)|,|N_G^2(x)|\}\geq (1+\eps)4k/3$, then $G$ contains a copy of every $k$-edge tree~$T$ with $\Delta(T)\leq \Delta$.
\end{theorem}

\subsection{Organisation of the paper}
In \cref{sec:preliminaries}, we outline the proofs of each of our results, all of which begin with the same unified strategy, and then require additional arguments that are specific to the scenario. We also introduce notation and some basic tools that will be used in the remainder of the paper. \cref{sec:tree_splitting} focuses on finding different ways to partition trees into smaller subtrees to satisfy various constraints. In \cref{sec:regularity}, we recall Szemer\'edi's regularity lemma and some related results, and use these to prove an embedding lemma for trees that will be required later. At the end of this section, we introduce some helpful tools due to Besomi, Pavez-Sign\'e and Stein about the structure of reduced graphs of tree-avoiding graphs. The general part of our strategy is detailed in  \cref{sec:mainlemma}, culminating with our main ingredient, \cref{lemma:no-proper-bad}, that allows us to convert the problems into the dense setting. \cref{thm:approx_k/2} will be proved at the end of this section without further difficulty, exemplifying the use of this main lemma. In \cref{sec:alpha}, we prove \cref{thm:mainalpha-exact}, by combining our methods with machinery of Besomi, Pavez-Sign\'e and Stein.
The stability arguments needed to achieve the exact bound in \cref{thm:exact2/3} will be provided in \cref{sec:2/3}, and we prove this theorem here. In \cref{sec:second_nbhd}, we prove some additional lemmas that allow us to deduce \cref{thm:secondnbhd_result}. Finally, we conclude in \cref{sec:conclusion} with a discussion of the methods, tightness of degree combinations, and related extensions.

\section{Preliminaries}\label{sec:preliminaries}
Throughout the asymptotic arguments, we omit integer parts and harmless additive $O(1)$ terms whenever they are not crucial to the argument. 
We first introduce standard notation. For $n,k\in \N$ with $k\leq n$ we write $[n] = \{1,\dots,n\}$,  $[n]_0 = [n]\cup \{0\}$ and $[n]^{(k)} = \{S\subseteq [n]: |S| = k\}$. Given integers $\ell$ and $m$ with $\ell<m$, we sometimes write $j\in [\ell,m]$ to mean $j\in [\ell,m]\cap \mathbb{Z}$ where it is clear from context that $j$ must be integer valued. If we say that a statement holds whenever $0 < a \ll b \leq1$, then there exists a non-decreasing function $f:(0,1] \rightarrow (0,1] $ such that the statement holds for all $0<a,b\leq1$ with $a\leq f(b)$. We similarly consider a hierarchy of constants $0 < b_1 \ll b_2 \ll \dots \ll b_k <1 $, and these constants must be chosen from right to left.

Given a graph $G$, we denote by $|G|$ and $e(G)$ the size of its vertex set and edge set respectively. The degree and neighbourhood of a vertex $x \in V(G)$ are denoted by $d_G(x)$ and $N_G(x)$ respectively, and the neighbourhood of a vertex set $S\subseteq V(G)$ is $N_G(S) = \bigcup_{x \in S}N_G(x)$.
The graph $G[S]$ is the subgraph of $G$ induced by $S$, and $G \setminus S$ is $G[V(G)\setminus S]$.
Given a subgraph $H\subset G$, we will sometimes write $G[H]$ to mean $G[V(H)]$ and $G- H$ to mean $G[V(G)\setminus V(H)]$. For a vertex $x\in V(G)$  we write $G - x$ instead of $G \setminus \{x\}$, and sometimes use $H\cup \{x\}$ instead of $G[V(H) \cup \{x\}]$ where it is clear from context what the underlying graph is. For convenience we often use $N_G(x)\cap H$ in place of $N_G(x)\cap V(H)$ and similarly $N_{G}(x)\setminus H$ for $N_{G}(x)\setminus V(H)$. For two sets $A,B \subset V(G)$, we write $e_G(A,B)$ for the number of edges in $G$ that go between $A$ and $B$. We may omit the subscripts where context is clear.

We now provide a few definitions. A \defn{cover} of a graph $G$ is a vertex set $U\subseteq V(G)$ such that every edge in $G$ has at least one endpoint in $U$. 
For two partitions $\mathcal{P}$ and $\mathcal{Q}$ of the same set, we say that $\mathcal{P}$ \defn{refines} $\mathcal{Q}$ if for every $P\in \mathcal{P}$ there exists $Q\in\mathcal{Q}$ such that $P\subseteq Q$.

For a tree $T$, note that there is a unique path between every pair of vertices in $T$, and the length of this path defines the distance between such a pair. For a vertex $v\in V(T)$, let \defn{$\textnormal{Even}_T(v)$} be the set of vertices in $T$ with even distance from $v$ (excluding $v$ itself), and let \defn{$\textnormal{Odd}_T(v)$} be the set of vertices in $T$ with odd distance from $v$. Additionally, if $T$ is rooted, with root $r$, then we think of the vertices in $T$ as having a partial order defined by their distance from $r$ in $T$. For vertices $u,v\in V(T)$, we say that $v$ is a \defn{descendant} of $u$ if $u$ lies on the unique path from $r$ to $v$ in $T$. If also $uv\in E(T)$, we say that $u$ is the \defn{parent} of $v$, and $v$ is a \defn{child} of $u$. For a vertex $u\in V(T)$, we define \defn{$T(u)$} to be the subtree of $T$ induced by $u$ and all of its descendants. 

At various points, we would like to embed sections of trees into different parts of our host graph, and it would be helpful if these parts possess some nice connectivity property.
The notion of connectedness that we will consider throughout is that of cut-density, as defined by Conlon, Fox, and Sudakov \cite{conlon2014cycle}. For $\rho\geq 0$, a graph $G$ is said to be \defn{$\rho$-cut-dense} if there is no partition of $V(G)$ into two sets $A$ and $B$ such that $e(A,B)<\rho\vert A \vert \vert B\vert$. Note that every graph is trivially $0$-cut-dense.

\subsection{Proof outline}\label{subsec:outline}

We remark that almost everything previously known about \cref{conj:2k/3,conj:alpha,conj:boundedk/2,conj:secondneighbourhood} holds only for dense graphs, that is, when $|G|= O(k)$. In order to tackle each of Theorems~\ref{thm:exact2/3},~\ref{thm:mainalpha-exact},~\ref{thm:approx_k/2}~and~\ref{thm:secondnbhd_result}, we proceed by one general strategy, that in some sense allows us to convert the problems into the dense setting. A key ingredient is the following structural decomposition theorem for tree-avoiding graphs due to the first author \cite{AP_hyperstability,APembedding_covers}, telling us that it is possible to delete few edges from $G$ to find a subgraph whose connected components each have small covers.
\begin{theorem}[Hyperstability theorem, Pokrovskiy {\cite{AP_hyperstability}}]\label{thm:CoverV2}
   For all $\Delta\in \N$, $\eps>0$ and sufficiently large $k\in \N$ the following holds. Let $T$ be a tree with $k$ edges and $\Delta(T)\leq \Delta$. For any graph $G $ not containing a copy of $T$, it is possible to delete $\eps k|G|$ edges to obtain a graph $G'$ each of whose connected components has a cover of size at most $ (2 + \eps)k$.
\end{theorem}
Motivated by this, we seek to find a collection of vertex-disjoint subgraphs in $G$, obtained by making refinements to the components obtained from an application of \cref{thm:CoverV2}, each satisfying various helpful properties that will later allow us to embed sections of trees into these subgraphs. For this purpose, we define the notion of a `rich' subgraph in $G$ (see \cref{defn:good_subgraph}), to be one that has a small cover, high minimum degree, bounded size, and some nice connectivity property, formalised as positive cut-density. This latter property will be useful for the additional arguments required for \cref{thm:exact2/3,thm:mainalpha-exact,thm:secondnbhd_result}. 

Building on this strategy, our proofs centre around a structural lemma (\cref{lemma:no-proper-bad}), which roughly says that if a graph $G$ with minimum degree $\delta k$ does not contain some bounded degree tree, then there exists a collection of vertex-disjoint rich subgraphs in $G$, each with minimum degree close to $\delta k$, and importantly, such that for every vertex $v\in V(G)$ there are two of these subgraphs that together contain almost all neighbours of $v$. Furthermore, we show that for every rich subgraph $C$ in the collection, there are only few vertices in $G$ that have many neighbours inside $C$ but that do not already lie within it. Together, these properties give us a strong picture about the behaviour of vertices with respect to the collection of rich subgraphs.

Before sketching the proof of \cref{lemma:no-proper-bad}, let us first provide a short explanation of how the lemma can be used in an approximate sense, with the aim to provide a heuristic overview. Throughout, we assume $k\in \N$ is sufficiently large with respect to parameters $\eps\in (0,1)$ and $\Delta \in \N$. 

Suppose we are given a (potentially sparse) graph $G$ satisfying $\delta(G)\geq (1+\eps)ak$ for some $a\geq 1/2$ such that $G$ does not embed a tree $T$ with $\Delta(T)\leq \Delta$. \cref{lemma:no-proper-bad} will allow us to find a collection $(C_i)_{i\in [m]}$ of vertex-disjoint subgraphs in $G$ that satisfy $\delta(C_i)\geq (1+\eps/4)ak$ and $|C_i|<100k$ for each $i$. If additionally $G$ satisfies $\Delta(G)\geq (1+\eps)bk$ for some $b\in [a,2]$ (think of the pair $a,b$ as a degree combination from one of the conjectures, e.g.\ $a=1/2$ and $b= 2(1-1/\Delta)$ from \cref{conj:boundedk/2}), then taking a maximum degree vertex $x\in V(G)$, the lemma tells us that there are rich subgraphs $C_i$ and $C_j$ such that $x$ has at most $\eps k/4$ neighbours outside of these two subgraphs. Taking $H$ to be the subgraph of $G$ induced on $V(C_i)\cup V(C_j) \cup \{x\}$, we have $\Delta(H)\geq d_G(x) - \eps k/4 \geq (1+\eps/4)bk$ and $\delta(H)\geq (1+\eps/4)ak$. Since $|H|\leq 1+|C_i|+|C_j|\leq 200k$, then we have found a dense graph $H$ that still has minimum and maximum degree asymptotically above the starting thresholds. Therefore, we can apply previous dense results such as \cref{thm:2k-1/Delta_dense} for appropriate $a,b$ pairs to embed the tree $T$. Thus, \cref{lemma:no-proper-bad} by itself is sufficient to prove \cref{thm:approx_k/2}. 

This strategy would also provide asymptotic resolutions of \cref{conj:2k/3,conj:alpha} for trees with maximum degree bounded by a constant, noting that we do not require the additional conclusion of the lemma as stated above. However, we must work harder to achieve the exact bounds in \cref{thm:exact2/3,thm:mainalpha-exact}.
For \cref{thm:mainalpha-exact}, we combine our use of \cref{lemma:no-proper-bad} with tools of Besomi, Pavez-Sign\'e and Stein (see \cref{subsec:BPS}) that allow us to analyse the structure of dense graphs that do not embed all bounded degree trees. 

To prove \cref{thm:exact2/3}, we essentially proceed in the same way, now applying \cref{lemma:no-proper-bad} to find a dense subgraph $H$ in $G$ with $\delta(H)\geq (2/3-\eps)k$ and $\Delta(H)\geq (1-\eps)k$. We require a sequence of stability arguments that allow us to study the structure of the subgraph $H$, and how other vertices in $G$ behave with respect to $H$, more closely. The subgraph $H$ is in some sense quite similar to the extremal graph depicted in \cref{fig:extremal2k/3}. Each subsection in \cref{sec:2/3} focusses on a different scenario for the structure of $H$ (and its surrounding vertices), and we consider these cases in turn to prove the theorem at the end of the section. More details are provided there.

Finally, for \cref{thm:secondnbhd_result}, we also start by applying \cref{lemma:no-proper-bad} to find the corresponding dense subgraph $H$ with minimum degree asymptotically above $k/2$, and containing the high degree vertex $x$ from the theorem. We know that $x$ also has many second neighbours, and if almost all of them lie within $H$, then we can apply \cref{thm:secondneighbourhood_dense} to $H$ to reach our desired conclusion. In the alternative case where $x$ has many second neighbours outside of $H$, we require some additional arguments to gain control over where these vertices may lie with respect to the remaining subgraphs in $(C_i)_{i\in [m]}$.  

Let us now sketch the proof strategy for \cref{lemma:no-proper-bad}, starting with a graph $G$ of minimum degree $(1+\eps)ak$ for some $a\geq 1/2$ and not containing a bounded degree $k$-edge tree $T$. Our first aim is to show that a collection of rich vertex-disjoint subgraphs can be found in $G$, so that additionally their union covers $(1-o(1))|G|$ vertices.
Having applied the machinery of the first author given by \cref{thm:CoverV2} with suitable parameters, and viewing the obtained components as a collection of vertex-disjoint subgraphs, we make some further refinements to ensure that each subgraph has minimum degree close to $(1+\eps)ak$ and is positively cut-dense, deleting few vertices from their union in the process. 
The cover property is also used to show that each subgraph is small, namely by showing that otherwise we can find high degree vertices and reduce the problem to the dense setting to embed $T$. Let us observe that only the minimum degree condition of $G$ is required for this first step.

Our second aim is to show that from a fixed collection of rich subgraphs, every vertex in $G$ is associated with two of these subgraphs that together contain almost all of its neighbours. We use a maximality argument, whereby the collection $(C_i)_{i\in[m]}$ of vertex-disjoint rich subgraphs is chosen such that the number of vertices covered by their union is as large as possible, subject to some pre-specified parametrisation. Recall that we also wish to understand the set of vertices that have many neighbours within each of these subgraphs. 
For each $i\in [m]$ and $d\in \mathbb{R}$, we define the \defn{$d$-periphery} set $L_d(C_i) := \{v\in V(G): \vert N_G(v)\cap V(C_i)\vert \geq d\}$, noting that $V(C_i)\subseteq L_d(C_i)$ for all $d\leq \delta(C_i)$ (see further \cref{defn:periphery}). We fix a small constant $\gamma$ and consider the subgraphs where the size of $L_{\gamma k}(C_i) \setminus V(C_i)$ is large separately from those subgraphs which do not expand much. Essentially, we prove that $\bigcup_{i\in [m]}L_{\gamma k}(C_i)$ contains all vertices of $G$, and that $|L_{\gamma k}(C_i)\setminus V(C_i)|\leq \gamma k$ for every $i\in [m]$. 

The main idea here is to show that if there are vertices outside of $\bigcup_{i\in [m]}L_{\gamma k}(C_i)$ or there are rich subgraphs $C_i$ having many vertices in the outer-periphery $L_{\gamma k}(C_i)\setminus V(C_i)$, then we can either control where their neighbours lie and embed $T$, or we can find a subgraph amongst these vertices that has minimum degree very close to $\delta k$. In the latter case, we can apply the same logic in the first step, to find a new collection of rich subgraphs, which, when combined appropriately with $(C_i)_{i\in [m]}$, yield a collection covering a larger number of vertices, contradicting the maximality assumption.
Once we have that $V(G) = \bigcup_{i\in [m]}L_{\gamma k}(C_i)$, we choose for each vertex $v\in V(G)$ the two subgraphs $C_i$ and $C_j$ from the collection that contain the first and second highest amount of neighbours of $v$. We then bound the number of neighbours of $v$ lying outside of these subgraphs, by showing that if there are many, there is a way to embed $T$ in $G$. This uses some fairly standard arguments about splitting up the tree $T$ into smaller subforests, and a more involved embedding lemma found in \cref{sec:regularity}.

\subsection{Basic results}

At various points throughout the proof, we will make use of the following fact, formalising a way to greedily embed a rooted tree into a graph with high minimum degree.

\begin{fact}[Greedy tree embedding, see e.g. \cite{AP_hyperstability}]\label{fact:greedy_embedding}
     Let $k,\Delta\in \N$. Let $T$ be a rooted tree on at most $k$ edges satisfying $\Delta(T) \leq \Delta$, and let $r$ denote its root. If $G$ is a graph containing a vertex $x$ such that $\delta(G -x)\geq k-1$ and $d(x) \geq d(r)$,  then~$G$ contains a copy of $T$ rooted at $x$.
\end{fact}

We will also need to bound the diameter of dense graphs.
\begin{theorem}[Erd\H{o}s, Pach, Pollack and Tuza {\cite{radiusErdos}}]\label{thm:mindeg_diameter}
Let $G$ be a connected graph on $n$ vertices with $\delta(G)\geq 2$. Then $$\textnormal{diam}(G)\leq \left\lfloor\frac{3n}{\delta(G)+1}\right\rfloor-1.$$ 
\end{theorem}

\section{Tree splitting}\label{sec:tree_splitting}

In this section, we introduce and prove several straightforward results about splitting up trees into pieces to satisfy various constraints, which will be combined later with tools for embedding these pieces separately into a host graph $G$. 
We start by noting that the bipartition classes of a bounded degree tree cannot be too unbalanced, as observed (in slightly weaker form) in \cite{APembedding_covers}.

\begin{fact}\label{fact:tree_part_sizes}
    Let $k,\Delta\in \N$. In every $k$-edge tree $T$ with maximum degree $\Delta$, the bipartition classes of $T$ both have at least $k/\Delta$ vertices.
\end{fact}
\begin{proof}
    Let $A$ be a bipartition class in $T$. We have $k = e(T) = \sum_{a\in A}d(a) \leq |A|\Delta$, and rearranging gives $|A|\geq k/\Delta$, as desired.
\end{proof}

Next, we introduce another standard fact, see for example \cite[Observation 2.3]{HRSW02497401}.
\begin{fact}\label{lem:deletedvertex_compsk/2}
    Every tree $T$ contains a vertex $v$ such that all components of $T -v$ have at most $|T|/2$ vertices. 
\end{fact}
We further wish to partition the components of $T-v$ into subforests of appropriate sizes. We apply the following lemma (see \cite[Lemma 4.4 and Remark 4.5]{BPSdegreeconditions2019}. The version below also has an $\ell\ge 2$ assumption which the lemma needs and which the proof in \cite{BPSdegreeconditions2019} implicitly used), in order to gain control of these sizes.
\begin{lemma}[Besomi, Pavez-Sign\'e and Stein {\cite{BPSdegreeconditions2019}}]\label{lem:sum_partition}
    Let $\ell,m\in \N$ with $\ell\ge 2$, let $(a_i)_{i\in [m]}$ be a sequence of nonnegative integers such that $a_i \leq \lceil \ell/2\rceil$ for each $i \in [m]$, and $\sum_{i\in[m]}a_i \leq \ell$. Then 

\begin{enumerate}[label = \upshape{(\roman*)}]
    \item there exists a partition $\{J_1,J_2\}$ of $[m]$ such that $\sum_{j\in J_2}a_j \leq \sum_{j\in J_1}a_j \leq \lfloor 2\ell/3 \rfloor$, and \label{item:2parts}
    \item there exists a partition $\{I_1,I_2,I_3\}$ of $[m]$ such that $\sum_{i\in I_3}a_i \leq \sum_{i\in I_2}a_i \leq \sum_{i\in I_1}a_i \leq \lceil \ell/2\rceil$, and moreover $|I_3|\leq 1$. \label{item:3parts}
\end{enumerate}
\end{lemma}

\begin{corollary}\label{cor:partition_subforests_T-v}
    Let $k\in \N$ with $k\ge 2$. Every $k$-edge tree $T$ contains a vertex $v \in V(T)$ such that
    \begin{enumerate}[label = \upshape{(\arabic*)}]
        \item there exists a partition of the forest $T - v$ into two vertex-disjoint subforests $F_1,F_2$ satisfying $ k/2 \leq |F_1|\leq \lfloor 2k/3 \rfloor$, and \label{item:Fsize1}
        \item there exists a partition of the forest $T - v$ into three vertex-disjoint subforests $F_1',F_2',F_3'$ each with at most $\lceil k/2 \rceil$ vertices and such that $F_3'$ is either empty or a tree. \label{item:Fsize2}
    \end{enumerate}
\end{corollary}

\begin{proof}
    Apply \cref{lem:deletedvertex_compsk/2} to obtain a vertex $v \in V(T)$ such that all components of $T - v$ each have at most $\lceil k/2 \rceil$ vertices. Denote the components by $A_1,\dots, A_m$ and let $a_i = |A_i|$ for each $i \in [m]$. Then $1\leq a_i\leq \lceil k/2 \rceil$  for each $i\in [m]$ and $\sum_{i\in[m]}a_i= |T -v| = k$, so we can apply \cref{lem:sum_partition} and denote by $\{J_1,J_2\}$ and $\{I_1,I_2,I_3\}$ the partitions of $[m]$ obtained from the first and second part respectively. For (1), let $F_j = \bigcup_{i\in J_j}A_i$ for $j\in \{1,2\}$. It follows from part (i) of the lemma that $|F_2|\leq |F_1|\leq \floor{2k/3}$. Since $|F_1 \cup F_2| = k$, then $|F_1|\geq k/2  $. For (2), consider the three subforests given by $F_j' = \bigcup_{i\in I_j}A_i$ for $j\in [3]$. Since $|I_3|\leq 1$, then $F_3'$ is either empty or a single component. The upper bounds on the sizes of the subforests follow directly from part (ii) of \cref{lem:sum_partition}.
\end{proof}

We require the following result which allows us to divide a tree at a specified vertex into two subtrees of controlled size.
\begin{proposition}[Montgomery {\cite{MONTGOMERY2019106793}}]\label{thm:subtree_m3m}
Let $m,t \in \mathbb{N}$ satisfy $m \leq t/3$. For every tree $T$ on $t$ vertices and any vertex $v \in V(T)$, there exist two subtrees $S_1$ and $S_2$ whose edge sets partition $E(T)$, that intersect in exactly one vertex, and that satisfy $v \in V(S_1)$ and $\vert S_2\vert \in [m,3m]$.
\end{proposition}

It is not too hard to deduce the following corollary, which divides the tree further and imposes an exact size for one of the subtrees. 

\begin{corollary}\label{cor:logk_subtrees}
Let $m,t\in \mathbb{N}$ satisfy $m\leq t$. For every tree $T$ on $t$ vertices, there exist subtrees $S_0,S_1,\dots,S_\ell$ whose edge sets partition $E(T)$, such that $\vert S_0\vert = m$ and $\ell\leq \lceil\log_{3/2}t\rceil+5$. Moreover we have $|V(S_0)\cap V(S_j)|=1$ for every $j\in [\ell]$.
\end{corollary}

\begin{proof}
Let $m$, $t$, and $T$ be as given.
Let $T_0$ be the empty graph and $R_0 = T$. We successively construct a sequence of pairs of subtrees $(T_i,R_i)_{i\in \N}$ such that for each $i\in \N$, all of the following hold:
\begin{enumerate}[label = (\roman*),leftmargin =\widthof{R1000}]
    \item $T_i \cup R_i = R_{i-1}$,
    \item $T_i$ and $R_i$ are edge-disjoint,
    \item $|R_i|\geq m$,
    \item $|R_i|-m-4\leq (2/3)^i(t-m-4)$
\end{enumerate} 
We first run this process until we have found an integer $r_0\geq0$ such that $|R_{r_0}|-m\leq4$.
Suppose $i\in \N$ and we have found all previous pairs satisfying the properties. Whilst $|R_{i-1}|-m>4$,  apply \cref{thm:subtree_m3m} with
$T_{\ref{thm:subtree_m3m}}=R_{i-1}$, $m_{\ref{thm:subtree_m3m}} = \left\lfloor\frac{|R_{i-1}|-m+1}{3}\right\rfloor$ and $v_{\ref{thm:subtree_m3m}}$ arbitrary in $V(R_{i-1})$ to obtain subtrees $T_i$ and $R_i$ of $R_{i-1}$ whose edge sets partition $E(R_{i-1})$ and such that $|T_i|\in \left[\left\lfloor\frac{|R_{i-1}|-m+1}{3}\right\rfloor,3\left\lfloor\frac{|R_{i-1}|-m+1}{3}\right\rfloor\right]$. In particular, we have
$|R_i|\geq |R_{i-1}|-3\left\lfloor\frac{|R_{i-1}|-m+1}{3}\right\rfloor+1\geq m$ and 
\begin{equation*}
    |R_i|-m -4\leq  |R_{i-1}|-\left\lfloor\frac{|R_{i-1}|-m+1}{3}\right\rfloor+1-m-4\le 
    \frac{2}{3}(|R_{i-1}|-m-4)\leq \left(\frac{2}{3}\right)^{i}(t-m-4).
\end{equation*}

It follows from (iv) that this first part of the process stops after at most $\lceil\log_{3/2}t\rceil+1$ steps. We then delete leaves of the current $R_{r_0}$ one at a time, taking the deleted edge as the next $T_i$, until we obtain $R_r$ with $|R_r|=m$. This requires at most four further steps, and hence $r\leq \lceil\log_{3/2}t\rceil+5$. 

Observe that, by construction, the collection of trees $\{T_i:i\leq r\}$ are pairwise edge-disjoint. Let $F$ be the subforest of $T$ with edge set $\bigcup_{i\leq r}E(T_i)$. Let $S_1,S_2,\ldots,S_\ell$ denote the connected components of $F$, so that each $S_j$ is formed by a union of trees from $\{T_i:i\leq r\}$. Note by (i) and (ii) that for each $i\leq r$, we have $E(T_{i})= E(R_{i-1})\setminus E(R_i)$. It follows that $E(F)\cup E(R_r) = (\bigcup_{i\leq r}E(T_i))\cup E(R_r)$, and therefore $ E(S_1)\cup \dots\cup E(S_\ell)\cup E(R_r)$, forms a partition of $E(T)$.
Finally, note that each $S_j$ is a tree and must intersect $R_r$ exactly once, as else a cycle is obtained in $T$, or $T$ is disconnected. 
Clearly $\ell \leq r\leq \lceil\log_{3/2}t\rceil+5$ and so taking $S_0 = R_r$ proves the corollary.
\end{proof}

\section{Regularity}\label{sec:regularity}

Let $G$ be a graph and $A,B $ be disjoint sets of vertices of $G$. The \defn{density} of the pair $(A,B)$ is given by $d(A,B) = \frac{e(A,B)}{|A||B|}$. For $\gamma>0$, the pair $(A,B)$ is \defn{$\gamma$-regular} if for all $A'\subseteq A$ and $B'\subseteq B$ satisfying $|A'|\geq \gamma|A|$ and $|B'|\geq \gamma|B|$ we have $d(A',B') = d(A,B) \pm \gamma,$
and such a pair $(A,B)$ is $(\gamma,\eta)$-regular if we additionally have $d(A,B)>\eta$.

We will use the degree form of Szemer\'edi's regularity lemma, see for instance \cite{KomlosSimonovits_reg_lemma}.
A vertex partition $V(G) = V_1 \cup  \dots  \cup  V_m$ is said to be an \defn{$(\gamma , \eta )$-regular partition} of $G$ if the following hold:
\begin{enumerate}
    \item $|V_1|  = | V_2|  =\dots  = | V_m |$,
    \item $V_i$ is an independent set for all $i \in  [m]$, and
    \item for all $1 \leq  i < j \leq  m$, the pair $(V_i, V_j )$ is $\gamma$-regular with density either $d(V_i, V_j ) > \eta$  or $d(V_i, V_j ) = 0$.
\end{enumerate}

\begin{lemma}[Degree form of Szemer\'edi's regularity lemma]\label{lemma:regularity-degree}
For all $\gamma>0$ and $m_0\in \N$, there are $N_0, M_0\in \N$ such that the following holds for all $\eta\in [0,1]$ and $n>N_0$. Any $n$-vertex graph $G$ has a subgraph $H$ with $\vert G\vert - \vert H\vert <\gamma n$ and $d_H(v)> d_G(v)- (\gamma+\eta)n$ for all $v\in V(H)$ such that $H$ has an $(\gamma, \eta)$-regular partition $V(H)=V_1\cup \dots \cup V_m$ with $m_0\leq m \leq M_0$.
\end{lemma}

The \defn{$(\gamma , \eta )$-reduced graph} $R$ of $G$, with respect to the $(\gamma , \eta )$-regular partition given by \cref{lemma:regularity-degree}, is the graph with vertex set $\{ V_i : i \in  [m]\}$ in which $V_iV_j$ is an edge if and only if $d(V_i, V_j ) > \eta$. We will sometimes refer to a nondescript $(\gamma , \eta )$-reduced graph $R$ without explicitly referring to the associated $(\gamma , \eta )$-regular partition of the graph. It is well known that the reduced graph inherits a similar minimum degree property as $G$.

\begin{fact}\label{fact:reduced-min-degree}
    Let $0<2\gamma\leq \eta \leq \alpha/2$. If $G$ is an $n$-vertex graph with $\delta(G)\geq \alpha n$, and $R$ is a $(\gamma, \eta)$-reduced graph of $G$, then $\delta(R)\geq (\alpha-2\eta)\vert R\vert$.
\end{fact}
The next two observations will allow us to preserve some kind of connectivity property for both the subgraph $H$ and the reduced graph $R$ obtained when applying the regularity lemma. Recall that $H$ is $\rho$-cut-dense if there is no partition of $V(H)$ into two sets $A$ and $B$ such that $e(A,B)<\rho\vert A \vert \vert B\vert$.

\begin{fact}[Pokrovskiy {\cite{APembedding_covers}}]\label{fact:cutdense_subgraph}
    Let $\rho>2\alpha>0$ and $n\in \N$.   Let $G$ be an $n$-vertex, $\rho$-cut-dense graph and $H$ a subgraph of $G$ such that every vertex $v\in V(H)$ has $d_G(v)-d_{H}(v)\leq \alpha n$. Then $H$ is $(\rho-2\alpha)$-cut-dense.
\end{fact}

\begin{fact}[Pokrovskiy  {\cite{APembedding_covers}}]\label{fact:cutdense_reduced_connected}
    Let $\eta>0$ and $\rho > \gamma>0$. Let $G$ be a $\rho$-cut-dense graph, and $R$ be a $(\gamma,\eta)$-reduced graph of $G$. Then $R$ is connected.
\end{fact}

Our main application of the regularity lemma here will be to prove \cref{lem:embedding_injection}, which allows us to embed trees such that certain vertices are mapped to a strategically chosen set. Before we introduce this, let us first observe the following simple proposition.

\begin{proposition}\label{lemma:short-even-path}
For all $\alpha\in (0,1)$ and $n\in \N$ the following holds. Let $G$ be a graph on $n$ vertices with $\delta(G)\geq\alpha n$. If $uv\in E(G)$ and there is an even walk from $u$ to a vertex $w$, then, for every even integer $L\geq 4/\alpha+2$, there is a $uw$-walk of length exactly $L$ whose first edge is $uv$.
\end{proposition}
\begin{proof}
   Let $P=w_0\dots w_d$ be a shortest walk of even length between $u$ and $w$, where $w_0=u$ and $w_d=w$. If $i,j\in [d]_0$ and $|i-j|$ is even and larger than $2$, then $N_G(w_i)\cap N_G(w_j)=\emptyset$, since otherwise we could replace the part of $P$ between $w_i$ and $w_j$ by a walk of length $2$ and obtain a shorter even $uw$-walk. Hence the sets $N_G(w_0),N_G(w_4),\dots,N_G\bigl(w_{4\lfloor d/4\rfloor}\bigr)$ are pairwise disjoint. Therefore
\[
n\geq \sum_{i=0}^{\lfloor d/4\rfloor}|N_G(w_{4i})|\geq \bigl(\lfloor d/4\rfloor+1\bigr)\alpha n,
\]
and consequently $d<4/\alpha$.

Since $L$ and $d$ are even and $L\geq 4/\alpha+2>d$, the number $q=(L-d)/2$ is a positive integer. Starting at $u$, traverse the two-edge walk $uvu$ exactly $q$ times, and then follow $P$. The resulting walk has length $2q+d=L$, starts with the edge $uv$, and ends at $w$.
\end{proof}

Suppose we have some tree $T^*$ and a subtree $T \subset T^*$. As mentioned in the introduction, having $\delta(G)\geq |T| - 1$ is sufficient for a graph $G$ to contain a copy of $T$, by a simple greedy argument. Sometimes, we will need to gain more control over the embedding of $T$ in $G$, for example ensuring some `special' vertices of $T$ are embedded into a predetermined part in $G$. This will help us to extend the embedding of $T$ to an embedding of $T^*$. This provides the motivation for the following lemma, where now we boost the minimum degree of $G$ to $\delta(G)\geq (1+\eps)|T|$, but we are more selective about where we embed our set of `special' vertices of $T$ (the set $L$ in the statement).

\begin{lemma}\label{lem:embedding_injection}
    For all $\eps\in (0,1)$, there exists $\delta\in (0,1)$ such that the following holds for all sufficiently large $t$. Let $T$ be a tree on $t$ vertices, let $S$ be a subtree of $T$ of size at least $(1-\delta)t$ and let $L$ be a set of leaves of $T$ that are not vertices of $S$. Let $G$ be a graph of size $n\leq 200t$ with $\delta(G)>(1+\eps)t$ and let $Y\subseteq V(G)$ have size at least $\eps t$.

    If all leaves in $L$ are in the same bipartition class of $T$ and have distance at least $1000$ from all vertices in $S$, then there exists an embedding $\phi\colon T\hookrightarrow G$ such that $\phi(v)\in Y$ for all $v\in L$.
\end{lemma}

\begin{proof}

    Let $\gamma \ll \eta \ll \eps$ and let $N_0,M_0\in \N$ be the outputs of \Cref{lemma:regularity-degree} with $\gamma_{\ref{lemma:regularity-degree}}=\gamma$ and $(m_0)_{\ref{lemma:regularity-degree}} = \gamma^{-1}$. Let $\delta = \eta^2/(8M_0)$ and let $t>N_0$. 
    Suppose that $S$, $T$, $G$, $Y$ and $L$ are as in the statement of the lemma and let $A$ be the set of vertices of $S$ that are adjacent to vertices in $T-S$. We may assume without loss of generality that all vertices in $A$ are in the same bipartition class of $T$ as the leaves in $L$. If this is not the case, there needs to be a vertex  $v\in A$ that has distance at least 1001 to all leaves in $L$, and we can add all neighbours of $v$ to $S$ (thus removing $v$ from $A$) without violating the condition that $S$ and $L$ have distance at least 1000.

    By our choice of $m_0$, $M_0$, and $N_0$, \cref{lemma:regularity-degree} yields a subgraph $H$ of $G$ with $\vert H\vert>\vert G\vert-\gamma n$, $\delta(H)\geq \delta(G)-(\eta+\gamma)n\geq (1+\eps/2)t$, and admitting a $(\gamma,\eta)$-regular partition $V(H) = V_1\cup\dots \cup V_m$. Observe that by \Cref{fact:reduced-min-degree}, the reduced graph $R$ on $\{V_i:i\in [m]\}$ corresponding to this partition has minimum degree at least $(\frac{\delta(H)}{n} - 2\eta)\vert R\vert>\vert R \vert/200$. 
  
    By a simple counting argument, there must exist $\alpha_0 \in [m]$ such that $\vert V_{\alpha_0}\cap Y\vert\geq \eps\vert V_{\alpha_0}\vert/400$. For each $\alpha\in [m]$, fix an arbitrary set $W_\alpha\subseteq  V_\alpha$ such that $\vert W_\alpha\vert =\eta \vert V_\alpha\vert$ with the additional restriction that $W_{\alpha_0}\subset Y$. Let $H'$ be the induced subgraph of $H$ on vertex set $V(H)\setminus \bigcup_{\alpha\in [m]} W_\alpha$, and note that $\delta(H')>(1+\eps/4)t$.

    Next, we fix an arbitrary vertex $u\in S$ that is in the same bipartition class of $T$ as the vertices in $A$ and $L$. Starting by mapping $u$ into $V_{\alpha_0}\setminus W_{\alpha_0}$ arbitrarily, we can construct an embedding $\phi_S\colon S\hookrightarrow H'$, by greedily embedding vertices in $S$ one-by-one, in increasing order of distance from $u$. This is possible since $\delta(H')>t$ and $\vert S\vert \leq t$. 

    For each $a\in A$, let $T_a$ be the component of $T- (S - a)$ that contains $a$ and note that the trees $T_a$ are disjoint for different $a$, as else $S$ would not be connected. We will successively construct embeddings $\phi_a: T_a \hookrightarrow H$ in such a way that for all $a\in A$, we have
    \begin{itemize}
        \item $\phi_a(a)=\phi_S(a)$,
        \item $\phi_a(V(T_a))\cap \left(\phi_S(V(S))\cup \bigcup_{b\in A\setminus \{a\}} \phi_b(V(T_b))\right)=\{\phi_a(a)\}$, and 
        \item $\phi_a(v)\in Y$ for each $v\in L\cap V(T_a)$.
    \end{itemize}
    
      Combining $\phi_S$ with all the embeddings $(\phi_a)_{a\in A}$, we obtain an embedding $\phi\colon T\hookrightarrow H$ with the desired properties.

      Suppose we have constructed such embeddings for a subset $B\subset A$, and let $a\in A\setminus B$. Put $x=\phi_S(a)$, fix $\beta_0\in[m]$ with $x\in V_{\beta_0}$, and choose $\beta_1\in[m]$ such that $|(N_H(x)\cap V_{\beta_1})\setminus\phi_S(S)|\geq \eps|V_{\beta_1}|/800$.

    Note that since the path from $u$ to $a$ in $T$ has even length, then the path from $\phi_S(u)$ to $\phi_S(a)=x$ in $H$ has even length as well. Since all edges in $H$ are between dense pairs of the regular partition, we may conclude that $R$ contains a walk of even length between $V_{\alpha_0}$ and $V_{\beta_0}$. Moreover, $V_{\beta_0}V_{\beta_1}\in E(R)$, since $x$ has a neighbour in $V_{\beta_1}$ and every pair of partition classes has density either zero or greater than $\eta$. Since $\delta(R)\geq |R|/200$, \cref{lemma:short-even-path} gives a walk $V_{\beta_0}V_{\beta_1}\dots V_{\beta_k}$ with $\beta_k=\alpha_0$ and $k=1000$.
    
    For a vertex $v\in V(T_a)$, let $h(v)$ be the distance of $v$ to $a$ in $T_a$ (in particular, $h(a)=0$). We now build $\phi_a$ by embedding vertices of $T_a$ in increasing order of $h$ such that for all $i\in \N$ and $v\in V(T_a)$ with $h(v)=i$ the following hold.
     \begin{itemize}
        \item $\phi_a(v)\notin \phi_S(S)\cup \bigcup_{b\in B} \phi_b(T_b)$ unless $v=a$ in which case $\phi_a(v)=\phi_S(a)$.
        \item If $i < k$, then $\phi_a(v)\in V_{\beta_i}$ and $\vert N_H(\phi_a(v))\cap V_{\beta_{i+1}}\setminus \phi_S(S)\vert \geq \eta^2 \vert V_{\beta_{i+1}}\vert/2$. When $i=k-1$, we additionally require $\vert N_H(\phi_a(v))\cap W_{\beta_k}\vert\geq \eta\vert W_{\beta_k}\vert/2$.
        \item If $i\geq k$ and $i$ is even, then $\phi_a(v)\in W_{\beta_k}$ and $\vert N_H(\phi_a(v))\cap V_{\beta_{k-1}}\setminus \phi_S(S)\vert \geq \eta^2 \vert V_{\beta_{k-1}}\vert/2$.
        \item If $i\geq k$ and $i$ is odd, then $\phi_a(v)\in V_{\beta_{k-1}}$ and $\vert N_H(\phi_a(v))\cap W_{\beta_{k}}\setminus \phi_S(S)\vert \geq \eta \vert W_{\beta_k}\vert/2$.
    \end{itemize}
    Since for all $v\in L\cap V(T_a)$ we have $h(v)\geq 1000$ and $h(v)$ is even, the resulting embedding satisfies $\phi_a(L\cap V(T_a))\subset W_{\beta_k} = W_{\alpha_0}\subset Y$.

    We start by embedding $\phi_a(a)=x\in V_{\beta_0}$ and note that by choice of $\beta_1$, $x$ has more than $\eta^2\vert V_{\beta_1}\vert/2$ neighbours in $V_{\beta_{1}}\setminus \phi_S(S)$. Suppose now that we have a partial embedding $\phi_a^*$ of $T_a$ that satisfies the above, and let $v\neq a$ be a minimal vertex of $T_a$ that has not yet been embedded. We discuss only the case that $i=h(v)<k$ as the other cases follow by analogous arguments.
    
    Let $w$ be the neighbour of $v$ on the path to $a$ and let $y=\phi_a^*(w)\in V_{\beta_{i-1}}$. Consider the set 
    \begin{equation*}
        P=N_H(y)\cap V_{\beta_{i}}\setminus \left(\phi_S(S) \cup \bigcup_{b\in B}\phi_b(T_b) \cup \text{im}(\phi_a^*)\right).
    \end{equation*}
     Since $\vert T - S\vert < \delta t \leq \eta^2t/8m \leq \eta^2 \vert V_{\beta_{i}}\vert/4$, and $\vert N_H(y)\cap V_{\beta_{i}}\setminus \phi_S(S)\vert \geq \eta^2 \vert V_{\beta_{i}}\vert/2$, we have $\vert P\vert \geq \eta^2 \vert V_{\beta_{i}}\vert/4$. Furthermore, $W_{\beta_{i+1}}\subset V_{\beta_{i+1}}$ is disjoint from $\phi_S(S)$ and has size $\eta \vert V_{\beta_{i+1}}\vert$. By regularity, $P$ must contain a vertex $z$ with at least $\eta \vert W_{\beta_{i+1}}\vert/2\geq \eta^2 \vert V_{\beta_{i+1}}\vert/2$ neighbours in $W_{\beta_{i+1}}$. In particular it has this many neighbours in $V_{\beta_{i+1}}$, and when $i=k-1$ this also gives the additional requirement involving $W_{\beta_k}$. Embedding $v$ onto $z$ therefore yields the desired properties for $\phi_a(v)$.
\end{proof}

\subsection[Tools of Besomi, Pavez-Signé and Stein]{Tools of Besomi, Pavez-Sign\'e and Stein}\label{subsec:BPS}
We end this section by stating a few important embedding results about reduced graphs.
In order for the additional analysis required to prove \cref{thm:mainalpha-exact,thm:secondnbhd_result}, we will use two further related results.
Besomi, Pavez-Sign\'e and Stein \cite{BPSmaxmin2020} showed that although having $\delta(G)\geq (1+\eps)k/2$ and $\Delta(G)\geq (1+\eps)4k/3$ is not sufficient to guarantee that every $k$-edge tree embeds in $G$, we can still gain some insight into the structure of $G$ when $G$ avoids some bounded degree tree. Roughly speaking, they show that, for a vertex $x$ of highest degree in $G$, in a reduced graph $R$ of $G-x$, there are two components in which $x$ has many neighbours, and $x$ does not have any neighbours in a third component. Moreover, the component in which $x$ has most neighbours is bipartite, $x$ has neighbours only in the larger bipartition class, and this class has size at least $2k/3$. To formalise this, we first need a simple definition from \cite{BPSdegreeconditions2019}. Given $\theta\in (0,1)$, a graph $G$, a vertex $x\in V(G)$, and a reduced graph $R$ of $G-x$, we say that $x$ \defn{$\theta$-sees} a component $S$ of $R$, if in $G$, $x$ has at least $\theta|\bigcup V(S)|$ neighbours in $\bigcup V(S)$.

The following is a compilation of results.  
More precisely, the first two parts \ref{item:BPS1} and \ref{item:BPS2} are given in \cite[Lemma 7.3]{BPSdegreeconditions2019}. We will make use of these two properties during the proof of \cref{thm:exact2/3}. The next two parts \ref{item:BPS3} and \ref{item:BPS4} are given in \cite[Theorem 1.4 and Definition 4.1]{BPSmaxmin2020}. (We implicitly use the definition of $(\eps,x)$-extremal graphs given in \cite{BPSmaxmin2020} in the statement below.) These will be required in the proof of \cref{thm:mainalpha-exact,thm:secondnbhd_result}. We remark that the original results cover several more cases that we choose not to include here, and that their work also extends to trees with maximum degree bounded by a root of $k$. 

\begin{theorem}[Besomi, Pavez-Sign\'e and Stein {\cite{BPSmaxmin2020}}]\label{thm:BPS_combined} 
For all $\eps \in(0,1)$, $\gamma \in (0,\eps^4/10^{10}]$ and $M_0\in \N$ there exists $k_0\in \N$ such that for all $n,k\geq k_0$ with $n\geq k\geq \eps n$, the following holds for every $n$-vertex graph $G$. 
Suppose $\delta(G)\geq (1+\eps)k/2$ and there exists a $k$-edge tree $T$ with $\Delta(T)\le k^{\frac{1}{67}}$ that does not embed in $G$.
Let $x \in V(G)$ and $R$ be a $(\gamma,5\sqrt{\gamma})$-reduced graph of $G - x$ with $|R|\leq M_0$. Then
\begin{enumerate}[label=\upshape{(R\arabic*)},leftmargin =\widthof{R10000}]
		\item every non-bipartite component of $R$ has at most $ (1+\sqrt[4]{\gamma})k \cdot \tfrac{|R|}{|G|}$ vertices; \label{item:BPS1}
        \item every bipartite component of $R$ has at most $ (1+\sqrt[4]{\gamma})k_1 \cdot \tfrac{|R|}{|G|}$ vertices in each of its bipartition classes, where $k_1$ is the size of the larger bipartition class of $T$; \label{item:BPS2}
        \end{enumerate}
		and furthermore if $d_G(x)\geq (1+\eps)4k/3$, then
		\begin{enumerate}[label=\upshape{(R\arabic*)},leftmargin =\widthof{R10000},resume]
		\item $x$ $\sqrt{\gamma}$-sees two components $S_1$ and $S_2$ of $R$ and $x$ does not have neighbours in any other component of~$R$; \label{item:BPS3}
		\item if additionally $|N_G(x)\cap \bigcup V(S_1)|\geq |N_G(x)\cap\bigcup V(S_2)|$, then $S_1$ is bipartite and its larger part has at least $(1+\sqrt[4]{\gamma})\tfrac{2k}{3} \cdot \tfrac{|R|}{|G|}$ vertices, and $x$ has only neighbours in the largest part. \label{item:BPS4}
	\end{enumerate}
\end{theorem}

We'll use the following routine corollary.
\begin{lemma}\label{lem:single-rich-BPS}
For all $\Delta\in\N$ and $\xi,\rho>0$, and all sufficiently large $k$, the following holds. Let $T$ be a $k$-edge tree with $\Delta(T)\leq\Delta$, let $C$ be a $\rho$-cut-dense graph with $|C|<100k$ and $\delta(C)\geq (1+\xi)k/2$, and let $x$ be a vertex of a graph $G$ containing $C$ such that
$ |N_G(x)\cap V(C)|\geq (1+\xi)4k/3.$
Then $G[V(C)\cup\{x\}]$ contains a copy of $T$.
\end{lemma}
\begin{proof}
    Suppose not. Let $\mu=\min\{\xi/10,1/200\}$ and choose $
 0<\gamma\ll \mu,\rho$ with $
 \gamma\leq \mu^4/10^{10}.$ 
Apply \cref{lemma:regularity-degree} to $C-x$ with parameters $\gamma$ and $5\sqrt\gamma$. We obtain a subgraph $H\subseteq C-x$ with $|C-x|-|H|<\gamma|C|$, a $(\gamma,5\sqrt\gamma)$-regular partition with reduced graph $R$ of bounded order, and $d_H(v)>d_{C-x}(v)-(\gamma+5\sqrt\gamma)|C|$
 for every $v\in V(H)$.
For $\gamma$ sufficiently small and $k$ sufficiently large this gives $\delta(H)\geq (1+2\mu)k/2$. Moreover every surviving vertex loses, on passing from $C$ to $H$, at most $1+(\gamma+5\sqrt\gamma)|C|$ incident edges. Hence \cref{fact:cutdense_subgraph} gives that $H$ is $\rho/2$-cut-dense, and therefore $R$ is connected by \cref{fact:cutdense_reduced_connected}.

Let $H_x=G[V(H)\cup\{x\}]$. Since fewer than $\gamma|C|$ vertices of $C-x$ were discarded,
we have $d_{H_x}(x)\geq (1+\xi)4k/3-\gamma|C|\geq (1+2\mu)4k/3$
while $\delta(H_x)\geq(1+\mu)k/2$. Also $d_{H_x}(x)>4k/3$ implies $|H_x|>k$, and $|H_x|<100k+1$ implies $k\geq\mu|H_x|$. For sufficiently large $k$ we have $\Delta(T)\leq k^{1/67}$. Thus \cref{thm:BPS_combined}, with its parameter $\eps$ equal to $\mu$, applies to $H_x$, $T$, and $R$. Property \ref{item:BPS3} requires $x$ to $\sqrt\gamma$-see two distinct components of $R$, contradicting that $R$ is connected.
\end{proof}

Finally, once we have the structure given by \cref{thm:BPS_combined}, it will be helpful to have an extra argument allowing us to embed trees using a bipartite reduced graph. As defined in \cite{BPSdegreeconditions2019}, a \defn{$(k_1, k_2, d)$-forest} is a forest $F$ with bipartition classes $A_1$ and $A_2$ such that $|A_1|\leq k_1$, $|A_2|\leq k_2$ and $\Delta(F)\leq (k_1+k_2)^d$. 
\begin{lemma}[Besomi, Pavez-Sign\'e and Stein {\cite[Corollary 5.4]{BPSdegreeconditions2019}}]\label{lem:BPSprescribed_forest_embedding_bipartite}
   For all $\gamma  \in  (0, 10^{-8})$ and for all $d, M_0 \in  \N$  there is $k_0$ such that for all $n, k_1, k_2 \geq  k_0$ the following holds. Let $G$ be an $n$-vertex graph having a $(\gamma , 5\sqrt{\gamma} )$-reduced graph $R$ that satisfies $| R|  \leq  M_0$, such that 
   \begin{itemize}
       \item $R = (X, Y ) $ is connected and bipartite,
       \item $\diam(R)\leq d$,
       \item $d_R(x) \geq  (1 + 100\sqrt{ \gamma})k_2 \cdot \frac{|R|} {n}$ for all $x \in  X$,
       \item $| X|  \geq  (1 + 100\sqrt{\gamma})k_1 \cdot  \frac{|R|} {n}$.
   \end{itemize}
   Then any $(k_1, k_2, 1/d)$-forest $F$, with bipartition classes $A_1$ and $A_2$ can be embedded into $G$, with $A_1$ going to $\bigcup V(X)$ and $A_2$ going to $\bigcup V(Y)$.
   Moreover, if $F$ has at most $ \frac{\gamma n}{| R|}$  roots, then the images of the roots going to $\bigcup V(X)$ can be mapped to any prescribed set of size at least $2\gamma|\bigcup V(X)|$  in $\bigcup V(X)$, and the images of the roots going to $\bigcup V(Y)$ can be mapped to any prescribed set of size at least $2\gamma |\bigcup V(Y)|$ in $\bigcup V(Y)$.
\end{lemma}

\section{Converting to the dense setting}\label{sec:mainlemma}
As mentioned in \cref{subsec:outline}, this section is dedicated to transforming our problem to the dense setting, with the end goal of proving \cref{lemma:no-proper-bad}. Playing an important role is that of the hyperstability theorem given by \cref{thm:CoverV2}.
Motivated by this theorem, and our earlier remarks about preserving minimum degree and cut-density, let us now formally introduce the notion of a \textit{rich} subgraph. 
\begin{definition}\label{defn:good_subgraph}
    For $k\in \N$ and $c,\rho \geq 0$, a subgraph $H\subseteq G$ is \defn{$(c,\rho,k)$-rich} if the following hold.
    \begin{itemize}
        \item The minimum degree of $H$ is at least $c k$,
        \item $H$ has a cover of size at most $3k$, 
        \item $H$ is $\rho$-cut-dense, and
        \item $|H|<100k$.
    \end{itemize}
\end{definition}

Our next three lemmas are developed to help us find rich subgraphs. First, we show that every large graph with few low degree vertices and a small cover can embed bounded degree trees.
\begin{lemma}\label{lem:size_subgraphcover}
    For all $\Delta\in \N$, $\eps\in (0,1)$ and sufficiently large $k \in \N$ the following holds. If $H$ is a graph such that $\vert H\vert \geq 100k$, $H$ has a cover of size at most $3k$, and all but at most $\eps^2|H|/50$ vertices in $H$ have degree at least $(1/2+\eps)k$, then $H$ contains a copy of every $k$-edge tree $T$ with $\Delta(T)\leq \Delta$.
\end{lemma}

\begin{proof}
Let $k$ be sufficiently large with respect to $\Delta$ and $\eps$, and suppose to the contrary that $H$ is a graph satisfying the assumptions which does not embed some $k$-edge tree $T$ with $\Delta(T)\leq \Delta$. Let $A$ denote a cover in $H$ of size at most $3k$. Let $U$ denote the set of vertices in $H$ of degree less than $(1/2+\eps)k$, noting that $|U|\leq \eps^2|H|/50$. First, we prove a weak upper bound on the number of vertices in $H$, which depends on $\eps$.

\begin{claim}\label{claim:weak_bound_H_size}
    $|H|\leq 25\eps^{-1}k$.
\end{claim}
\begin{proofclaim}
Suppose $|H|>25\eps^{-1}k$. Choose $B_0 \subseteq V(H) \setminus (A\cup U)$ of size $20\eps^{-1}k$ arbitrarily. Let $A_0 \coloneqq A$ and define $H_0 = H[A_0\cup B_0]$. We will successively construct subsets $A_0\supset A_1 \supset \dots \supset A_t $ and  $B_0\supset B_1 \supset \dots \supset B_t$, and define subgraphs $H_0 \supset H_1 \supset \dots \supset H_t$ so that $H_i = H[A_i \cup B_i]$ for each $i$. 
 We stop this process at step $t$ when we have arrived at $H_t$ that additionally satisfies $d_{H_t}(a)\geq 3k$ for every $a\in A_t$. 

For $i\geq 1$, suppose we have already constructed the sets $A_{i-1}$, $B_{i-1}$, and the graph $H_{i-1}$ as above. If all vertices in $A_{i-1}$ have degree at least $3k$ and all vertices in $B_{i-1}$ have degree at least  $(1/2+\eps/2)k$ in $H_{i-1}$, 
    then stop the process with $t \coloneqq i-1$. Otherwise, apply the following steps to $H_{i-1}$:
    \begin{enumerate}
        \item Delete all vertices from $A_{i-1}$ which have degree less than $3k$ in $H_{i-1}$. Call the resulting subset $A_i$.
        \item Delete all vertices from $B_{i-1}$ which have degree less than $(1/2+\eps/2)k$ in $H_{i-1}$. Call the resulting subset $B_i$.
        \item Define $H_i = H[A_i\cup B_i]$.
    \end{enumerate}
    Consider the graph $H_t$ at the time this process stops. 
    First we show that $H_t$ is non-empty. In order to obtain $H_t$ from $H_0$, by deleting vertices from $A$ we have removed at most $3k|A| \leq 9k^2$ edges. Since each vertex deleted from $B_0$ must have intersected at least $\eps k /2$ of these removed edges, then in total we have deleted at most $9k^2/(\eps k /2) = 18\eps^{-1}k$ vertices from $B_0$. Since $|B_0| >18\eps^{-1}k$ there are still some vertices remaining, implying $H_t \neq \emptyset$.
    
    By definition of the process, it follows that $\delta(H_t)\geq (1/2+\eps/2)k>(1+\eps/25)k/2$. Moreover $B_t$ is independent, so $B_t\neq\emptyset$ and the positive minimum degree imply $A_t\neq\emptyset$. Hence the stopping condition gives $\Delta(H_t)\geq 3k$. Thus $\Delta(H_t)\geq 3k> 2(1-1/\Delta+\eps/25)k$.  Since $|H_t|\leq |A_0|+|B_0|\leq 25\eps^{-1}k$ then assuming $k$ is sufficiently large, $H_t$ satisfies the assumptions of
    \cref{thm:2k-1/Delta_dense} with $\eps_{\ref{thm:2k-1/Delta_dense}}=\eps/25$, 
    implying $H_t$ contains a copy of $T$, a contradiction. This proves the claim.
\end{proofclaim}

Now with \cref{claim:weak_bound_H_size} in hand, we will check that the minimum degree, maximum degree, and size of the graph 
$H^* = H \setminus U$ also satisfy the assumptions of \cref{thm:2k-1/Delta_dense} in order to find a copy of $T$ in this subgraph. For the minimum degree condition, note by \cref{claim:weak_bound_H_size} that $|U|\leq \eps^2|H|/50 \leq \eps k/2$, and so $\delta(H^*)\geq (1/2+\eps)k-|U| > (1+\eps/25)k/2$. 

Next, suppose every vertex in $A$ has degree less than $3k$ in $H^*$, and let $X= V(H^*)\setminus A$. Then
\begin{equation*}
    e_{H^*}(A,X)\leq \sum_{a\in A}|N_{H^*}(a)\cap X| < |A|\cdot 3k \leq 9k^2.
\end{equation*}
On the other hand, we have $|X| \geq |H|-|U|-|A| \geq 96k$. Since $A$ is a cover, each vertex in $X$ has all of its neighbours in $A$. Therefore
\begin{equation*}
    e_{H^*}(A,X)\geq \sum_{x\in X}d_{H^*}(x) \geq |X|\cdot k/2  > 9k^2,
\end{equation*}
a contradiction. Thus there exists $a\in A$ of degree at least $3k$, and in particular $\Delta(H^*)\geq 3k>2(1 - 1/\Delta+\eps/25)k$.  Finally we have $|H^*|\leq |H|\leq 25\eps^{-1}k$. Thus $H^*$ satisfies the assumptions of \cref{thm:2k-1/Delta_dense} again with $\eps_{\ref{thm:2k-1/Delta_dense}}=\eps/25$, and we find a copy of $T$ in $H^* \subseteq H$, again a contradiction. 
\end{proof}

We remark that, when we look to find a rich subgraph in a graph $G$ that does not embed every bounded degree $k$-edge tree, it will suffice to find a subgraph $H$ satisfying only the first three properties of \cref{defn:good_subgraph}. This is because, given \Cref{lem:size_subgraphcover}, the final property is implied immediately by the first two properties. This will be seen in action later, in the proof of \cref{lemma:good-components}. 
Before introducing this, we first prove that every graph with high minimum degree compared to its order has a large subgraph, also with high minimum degree, whose connected components each have positive cut-density.

\begin{lemma}\label{lem:cutdense/mindegree}
    For all $\eps\in (0,1)$, $\delta\in (0,\eps/400)$ and sufficiently large $k\in \N$ the following holds for all $a\geq1/2$. Let $G$ be a graph with $\delta(G)\geq (a+\eps)k$ and $|G|\leq100k$. Then there exists a subgraph $H$ of $G$ obtained by deleting at most $200\delta \vert G\vert$ vertices, such that $\delta(H)\geq (a+\eps-400\delta)k$ and each component of $H$ is $(\delta^2/20000)$-cut-dense.
\end{lemma}

\begin{proof}
Let $\eps$ and $\delta$ be as in the statement, let $\rho =\delta^2/20000$, and let $k$ be sufficiently large. Let $G_0 \coloneqq G$. We will successively construct subgraphs $G_0 \supset G_1 \supset G_2 \supset \dots$ such that $\vert G_i\vert\geq \vert G\vert - i\delta k$ and $\delta(G_i)\geq (a+\eps-2i\delta)k$ for each $i\geq 0$ until we arrive at a subgraph $G_t$ with the property that each of its components is $\rho$-cut-dense. Each subgraph $G_i$ will have at least $i$ connected components, and since each component of $G_i$ has at least $\delta(G_i)\geq k/2$ vertices, the process must terminate with a subgraph $G_t$ where $t\leq \vert G\vert/(k/2)\leq 200$. Choosing $H \coloneqq G_t$ will give the desired properties. 

For $i\geq 1$, suppose we have already constructed $G_{i-1}$ to satisfy the listed conditions. If every component of $G_{i-1}$ is $\rho$-cut-dense, then stop the process with $t \coloneqq i-1$. So suppose otherwise that $G_{i-1}$ has a component $C$ which is not $\rho$-cut-dense. Consider a partition $A\cup B$ of $V(C)$ such that $e_{G_{i-1}}(A,B) = e_{C}(A,B)\leq \rho|A||B|$. In order to construct $G_{i}$ from $G_{i-1}$, we first delete all edges between $A$ and $B$ and call the resulting graph $\hat{G}$, noting that $\hat{G}$ has at least one more component than $G_{i-1}$. Let $S = \{v\in C: d_{\hat{G}}(v)<(a+\eps-(2i-1)\delta)k\}$. Recall that $\delta(G_{i-1})\geq (a+\eps-2(i-1)\delta)k$, so that every $v\in S$ must have been incident to at least $\delta k$ edges of $E_{C}(A,B)$ in order for its degree to have dropped below the threshold. Therefore, we have
\begin{equation*}\label{eqn:sizeV_jbad}
    |S|\leq \frac{2e_{C}(A,B)}{\delta k}\leq \frac{2\rho|C|^2}{\delta k}\leq \delta k,
\end{equation*}
where the last inequality holds since $\vert C\vert\leq 100k$ and $\rho=\delta^2/20000$.

Let $G_{i}$ be the graph obtained from $\hat{G}$ by deleting all vertices in $S$. Thus, $\vert G_{i}\vert = \vert G_{i-1}\vert - \vert S\vert\geq \vert G\vert - i\delta k$. Furthermore, the minimum degree of $G_{i}$ is at least $(a+\eps-(2i-1)\delta)k-\vert S\vert\geq (a+\eps-2i\delta)k$. 

Finally it remains to check that $G_{i}$ has at least $i$ components. Recall that $G_{i-1}$ had at least $i-1$ components, and so as already observed, $\hat{G}$ must have at least $i$ components. By choice of the partition, we have
\begin{equation*}
    \rho|A||B|\geq e_C(A,B) \geq \delta(G_{i-1})|A| - 2e(G_{i-1}[A])\geq (a+\eps-2(i-1)\delta)k|A| - |A|^2,
\end{equation*}
which we can rearrange to obtain $|A|\geq (a+\eps-2(i-1)\delta)k-\rho|B|$. Since $|B|\leq 100k$ it follows that $|A|\geq (a+\eps-2(i-1)\delta- 100\rho)k\geq ak>|S|$. An identical calculation also shows that $|B|>|S|$. Hence deleting $S$ from $\hat{G}$ cannot remove either component $\hat{G}[A]$ or $\hat{G}[B]$ completely. Also, no component of $\hat{G}$ other than $C$ is altered in obtaining $G_{i}$. Therefore $G_{i}$ is as desired and this completes the proof.
\end{proof}

We are now ready to apply \cref{thm:CoverV2} to `almost decompose' a graph $G$ into vertex-disjoint rich subgraphs. 

\begin{lemma}\label{lemma:good-components}
For all $\Delta\in \N$, $\eps\in(0,1)$ and $\delta\in (0,\eps/6)$ the following holds for all sufficiently large $k\in \N$ \textup{(}in terms of $\Delta$, $\eps$ and $\delta$\textup{)} and all $a\geq 1/2$. Let $T$ be a $k$-edge tree with $\Delta(T)\leq \Delta$. If 
$G$ is a graph with $\delta(G)\geq(a+\eps)k$ such that $T$ cannot be embedded into $G$, then there exist vertex-disjoint $(a+\eps/2,\delta^2/10^{10}, k)$-rich subgraphs $C_1,\dots, C_m$ whose union contains all but at most $\delta\vert G\vert$ vertices in $G$.
\end{lemma}
\begin{proof}
    Let $\Delta$, $\eps$ and $\delta$ be as in the statement and let $\gamma = \delta \eps^3/19200$. Let $k$ be sufficiently large in terms of all parameters, let $a\geq 1/2$, and suppose $G$ and $T$ are as in the statement of the lemma, with $n = |G|$. Since $G$ contains no copy of $T$, we may apply \cref{thm:CoverV2} with $\eps_{\ref{thm:CoverV2}} =\gamma/2$ to find a spanning subgraph $G' \subseteq G$ obtained by deleting less than $\gamma kn$ edges from $G$, whose connected components $D_1,\dots, D_t$ each have a cover of order at most $3k$.
    Let $U = \{v \in V(G): d_{G'}(v) < (a+3\eps/4)k\}$, and note that every vertex in $U$ must belong to at least $\eps k /4$ edges of $E(G)\setminus E(G')$. Therefore we have
\begin{equation*}
    2\gamma k n \geq \sum_{u \in U}d_G(u)-d_{G'}(u) \geq |U| \left((a+\eps)k - (a+3\eps/4)k\right) = |U|\frac{\eps k}{4},
\end{equation*}
    and so $|U|\leq 8\gamma \eps^{-1}n = \delta \eps^2 n/2400$.

For each $i\in [t]$, let $U_i = U \cap D_i$. Let $J=\{j\in[t]: |U_j|>\eps^2|D_j|/1200\}$ and let $I = [t]\setminus J$.
Suppose there exists $i\in I$ such that $|D_i|\geq100k$. Since $D_i$ has a cover of size at most $3k$, and all but $|U_i|\leq \eps^2|D_i|/1200 \leq (\eps/4)^2|D_i|/50$ vertices of $D_i$ have degree at least $(a+3\eps/4)k > (1/2+\eps/4)k$, then an application of \cref{lem:size_subgraphcover} with $H_{\ref{lem:size_subgraphcover}} = D_i$, $U_{\ref{lem:size_subgraphcover}} = U_i$ and $\eps_{\ref{lem:size_subgraphcover}}=\eps/4$ shows that every $k$-edge tree $T$ with $\Delta(T)\leq \Delta$ embeds in $D_i$, and hence in $G$, a contradiction. Thus we may assume $|D_i|<100k$ for all $ i\in I$.

Now, for each $i\in I$, let $H_i = G'[V(D_i)\setminus U_i]$. Every $v \in V(H_i)$ satisfies 
\begin{equation*}
d_{H_i}(v) \geq d_{G'}(v)-|U_i| \geq (a+3\eps/4)k - \eps^2|D_i|/1200 \geq (a+3\eps/4)k - \eps^2 k/12\geq (a+2\eps/3)k.    
\end{equation*}
Therefore $\delta(H_i)\geq (a+2\eps/3)k$. Recall that $D_i$ has a cover of size at most $3k$ and since this cover property is inherited by all subgraphs, $H_i$ also has a cover of size at most $3k$.

Now, for each $i\in I$, let us apply \cref{lem:cutdense/mindegree} with $\eps_{\ref{lem:cutdense/mindegree}} = 2\eps/3$, $\delta_{\ref{lem:cutdense/mindegree}} = \delta/400$ and $G_{\ref{lem:cutdense/mindegree}} = H_i$ to obtain a subgraph $H_i'$ of $H_i$ such that each component of $H_i'$ is $(\delta^2/10^{10})$-cut-dense, $\vert H_i'\vert \geq (1-200\delta/400)\vert H_i\vert>(1-\delta/2)\vert H_i\vert$, and $\delta(H_i')\geq (a+2\eps/3-\delta)k>(a+\eps/2)k$. In particular, each component of $H_i'$ has at most $|H_i'|\leq |H_i|\leq |D_i|<100k$ vertices, and a cover of size at most $3k$. Thus every component of $H_i'$ is an $(a+\eps/2,\delta^2/10^{10}, k)$-rich subgraph.
Let us label all components of $\bigcup_{i\in I}H_i'$ by $C_1,\dots,C_m$.

All vertices in $\bigcup_{i\in I}V(H_i')$ are contained in one of these rich subgraphs. The number of vertices outside of this set is
\begin{align*}
    \sum_{j\in J}|D_j| + \sum_{i\in I} \left(\vert V(D_i)\setminus V(H_i)\vert + \vert V(H_i)\setminus V(H_i')\vert \right) < \sum_{j\in J}|D_j| + \sum_{i\in I} \vert U_i\vert + \delta n/2.
\end{align*}

Since $|D_j|\leq 1200\eps^{-2}|U_j|$ for every $j\in J$ and $\vert U \vert \leq \delta \eps^2 n/2400$, we have
\begin{equation*}
    \sum_{j\in J}|D_j|+\sum_{i\in I}|U_i| \leq \sum_{j\in J}1200\eps^{-2}|U_j| +\sum_{i\in I}|U_i| \leq \sum_{\ell \in [t]}1200\eps^{-2}|U_{\ell}|= 1200\eps^{-2}|U| <\delta n/2,
\end{equation*}
which completes the proof. 
\end{proof}

We now build up results that will help us to prove the `second aim' as described in \cref{subsec:outline}, with the end goal of proving \cref{lemma:no-proper-bad}. To understand this lemma, we need the following definition.
\begin{definition}\label{defn:periphery}
    Let $G$ be a graph and let $S\subseteq V(G)$. For $d\in \mathbb{R}$, the \defn{$d$-periphery} of $S$ is the set
\begin{align*}
    L_d(S)=\{v\in V(G): \vert N_G(v)\cap S\vert \geq d\}.
\end{align*}
For a subgraph $H$ of $G$ we write $L_d(H)$ to mean $L_d(V(H))$. 
\end{definition}

Note that $L_d(S)\subseteq L_{b}(S)$ for all $d\geq b$, and if $H$ has minimum degree at least $d$, then $V(H)\subseteq L_d(H)$. Given some vertex-disjoint rich subgraphs within a graph $G$, we are interested in understanding how periphery sets of these subgraphs can interact with each other.
\begin{lemma}\label{lem:intersectionproperties}
    For all $\Delta\in \N$, $\eps \in(0,1)$, and sufficiently large $k \in \N$ the following holds. Let $T$ be a $k$-edge tree with $\Delta(T)\leq \Delta$. If 
$G$ is a graph with vertex-disjoint $(1/2+\eps,0,k)$-rich subgraphs $C_1, \dots, C_m$ such that $T$ cannot be embedded into $G$, then all of the following hold.
    \begin{enumerate}[label = \upshape{(L\arabic*)}, leftmargin = \widthof{L10000}]
            \item For all distinct $i,j,\ell \in [m]$, we have $L_\Delta(C_i)\cap L_\Delta(C_j)\cap L_1(C_\ell) = \emptyset$. \label{item:periperhery_intersection1}
        \item For all $I\in [m]^{(2\Delta)}$, we have $\bigcap_{i\in I} L_1(L_{2\Delta}(C_i))=\emptyset$. In particular, $\bigcap_{i\in I} L_1(C_i)=\emptyset$. \label{item:periperhery_intersection2}
        \item For all $I\in [m]^{(2)}$, we have $\vert\bigcap_{i\in I}  L_{\eps k}(C_i)\vert < \eps k$. \label{item:periperhery_intersection3}
    \end{enumerate}
    
\end{lemma}
\begin{figure}[bp]
    \centering
\begin{tikzpicture}[u/.style={fill=black,minimum size =6pt,ellipse,inner sep=1pt},node distance=1.5cm,scale=0.8]
\node[u] (x) at (0.8,1){};
\node[] (v2) at (-0.9,-0.8){};
\node[] (v3) at (-1.58,-0.8){};
\node[] (v4) at (-1.9,-0.6){};
\node[] (v5) at (0.2,-1){};
\node[] (v6) at (0.8,-1){};
\node[] (v7) at (1.4,-1){};
\node[] (v8) at (3,-0.8){};

\draw[line width=1pt] (x) -- (v2);
\draw[line width=1pt] (x) -- (v3);
\draw[line width=1pt] (x) -- (v4);
\draw[line width=1pt] (x) -- (v5);
\draw[line width=1pt] (x) -- (v6);
\draw[line width=1pt] (x) -- (v7);
\draw[line width=1pt] (x) -- (v8);
\node[above=1mm,font=\small] at (x) {$x$};
\node[] at (0.8,-3) {\ref{item:periperhery_intersection1}};

  \def\r{0.8}
  \filldraw[line width=1pt, fill=white] (-1.3,-1.1) circle (\r);
  \filldraw[line width=1pt, fill=white] (0.8,-1.1) circle (\r);
  \filldraw[line width=1pt, fill=white] (2.9,-1.1) circle (\r);
  \node[] at (-1.3,-1.1) {$C_1$};
  \node[] at (0.8,-1.1) {$C_2$};
  \node[] at (2.9,-1.1) {$C_3$};
    \node[font=\footnotesize] at (-0.8,0.5) {$\Delta$};
\end{tikzpicture}\hspace{1.2cm}
\begin{tikzpicture}[u/.style={fill=black,minimum size =5pt,ellipse,inner sep=1pt},node distance=1.5cm,scale=0.8]
\node[u] (x) at (2.5,1.5){};
\node[u] (y1) at (0,0.6){};
\node[u] (y2) at (2,0.6){};
\node[u] (y3) at (5,0.6){};
\node[] (v2) at (0.3,-1){};
\node[] (v3) at (0.8,-0.9){};
\node[] (v4) at (-0.3,-1){};
\node[] (v4b) at (-0.8,-0.9){};
\node[] (v5) at (2.8,-0.9){};
\node[] (v6) at (2.3,-1){};
\node[] (v7) at (1.7,-1){};
\node[] (v7b) at (1.2,-0.9){};
\node[] (v8) at (4.2,-0.9){};
\node[] (v9) at (4.7,-1){};
\node[] (v10) at (5.3,-1){};
\node[] (v11) at (5.8,-0.9){};

\draw[line width=1pt] (x) -- (y1);
\draw[line width=1pt] (x) -- (y2);
\draw[line width=1pt] (x) -- (y3);
\draw[line width=1pt] (y1) -- (v2);
\draw[line width=1pt] (y1) -- (v3);
\draw[line width=1pt] (y1) -- (v4);
\draw[line width=1pt] (y1) -- (v4b);
\draw[line width=1pt] (y2) -- (v5);
\draw[line width=1pt] (y2) -- (v6);
\draw[line width=1pt] (y2) -- (v7);
\draw[line width=1pt] (y2) -- (v7b);
\draw[line width=1pt] (y3) -- (v8);
\draw[line width=1pt] (y3) -- (v9);
\draw[line width=1pt] (y3) -- (v10);
\draw[line width=1pt] (y3) -- (v11);
\node[above=1mm] at (x) {$x$};
\node[] at (2.5,-2.8) {\ref{item:periperhery_intersection2}};
\node[] at (3.5,0.6) {$\dots$};

  \def\r{0.8}
  \filldraw[line width=1pt, fill=white] (0,-1) circle (\r);
  \filldraw[line width=1pt, fill=white] (2,-1) circle (\r);
  \filldraw[line width=1pt, fill=white] (5,-1) circle (\r);
  \node[] at (0,-1) {$C_1$};
  \node[] at (2,-1) {$C_2$};
  \node[] at (5,-1) {$C_{2\Delta}$};
  \node[font=\footnotesize] at (-0.72,0.12) {$2\Delta$};
  \node[font=\footnotesize, left] at (y1) {$y_1$};
  \node[font=\footnotesize, left] at (y2) {$y_2$};
  \node[font=\footnotesize, right] at (y3) {$y_{2\Delta}$};

\end{tikzpicture}\hspace{1.4cm}
\begin{tikzpicture}[u/.style={fill=black,minimum size =4pt,ellipse,inner sep=1pt},node distance=1.5cm,scale=0.8]

  \def\r{0.8}
 \filldraw[line width=1pt, fill=white] (1,1) ellipse (4mm and 7mm);
\node[u] (v0) at (1,0.8){};
\node[] (v1) at (0.25,-0.8){};
\node[] (v2) at (0.5,-0.8){};
\node[] (v3) at (-0.25,-0.8){};
\node[] (v4) at (0,-0.8){};
\node[] (v5) at (-0.5,-0.8){};
\node[] (v6) at (1.75,-0.8){};
\node[] (v7) at (1.5,-0.8){};
\node[] (v8) at (2.5,-0.8){};
\node[] (v9) at (2,-0.8){};
\node[] (v10) at (2.25,-0.8){};

\draw[line width=1pt] (v0) -- (v1);
\draw[line width=1pt] (v0) -- (v2);
\draw[line width=1pt] (v0) -- (v3);
\draw[line width=1pt] (v0) -- (v4);
\draw[line width=1pt] (v0) -- (v5);
\draw[line width=1pt] (v0) -- (v6);
\draw[line width=1pt] (v0) -- (v7);
\draw[line width=1pt] (v0) -- (v8);
\draw[line width=1pt] (v0) -- (v9);
\draw[line width=1pt] (v0) -- (v10);

  \filldraw[line width=1pt, fill=white] (0,-1) circle (\r);
  \filldraw[line width=1pt, fill=white] (2,-1) circle (\r);

  \node[] at (0,-1) {$C_1$};
  \node[] at (2,-1) {$C_2$};
  \node[font=\footnotesize] at (1.7,1.5) {$X$};
  \node[font=\footnotesize] at (-0.05,0.21) {$\eps k$};

\node[u] (v0) at (1,0.8){};
\node[right,font=\small,rotate=90] at (v0) {$\dots$};
\node[below=2cm] at (1,0) {\ref{item:periperhery_intersection3}};
\end{tikzpicture}
\caption{\centering Forbidden structures in $G$ as detailed in the proof of \cref{lem:intersectionproperties}. Note that $x$, $y_1,\dots,y_{2\Delta}$ and $X$ need not actually be disjoint from the rich subgraphs.}
\end{figure}

\begin{proof}
\ref{item:periperhery_intersection1}: Suppose not, by relabelling, without loss of generality we may assume that there exists $x \in L_\Delta(C_1)\cap L_\Delta(C_2)\cap L_1(C_3)$. By \cref{cor:partition_subforests_T-v}\ref{item:Fsize2} there exists a vertex $r \in V(T)$ such that $T-r$ can be partitioned into three vertex-disjoint subforests $F_1$, $F_2$, $F_3$, each with at most $\ceil{k/2}$ vertices, and such that $F_3$ is either empty or a tree. For $i \in \{1,2\}$, we have $|N(x)\cap C_i|\geq \Delta$ and $\delta(C_i - x) \geq \ceil{k/2}\geq |F_i|$. By \Cref{fact:greedy_embedding}, and considering $r$ to be the root of the tree $F_i \cup \{r\}$, there is a copy of $F_i\cup \{r\}$ in $C_i\cup \{x\}$ where $r$ is mapped at $x$. If $F_3$ is empty, the first two embeddings already give $T$. Otherwise, $r$ has exactly one neighbour $u$ in $F_3$ since it is a tree, and we consider $u$ to be the root of the subtree $F_3$. Let $y$ be a neighbour of $x$ in $V(C_3)$. Since $\delta(C_3 - \{x,y\})\geq \ceil{k/2}\geq |F_3|$ and $y$ has at least $\Delta$ neighbours in $C_3-x$, again applying \cref{fact:greedy_embedding} we find a copy of $F_3$ in $C_3 -x$ such that $u$ is embedded at $y$. Since the three embeddings of the subforests $F_1$, $F_2$ and $F_3$ are compatible with $r$ being mapped at $x$, and $C_1$, $C_2$ and $C_3$ are vertex-disjoint, this yields a copy of $T$ in $G$, a contradiction.

\vspace{0.2cm}

\ref{item:periperhery_intersection2}: Suppose not, without loss of generality $I = \{1,\dots,2\Delta\}$ is such that there exists $x\in \bigcap_{i\in I}L_1(L_{2\Delta}(C_i))$. For each $i\in [2\Delta]$, denote by $y_i$ a neighbour of $x$ in $L_{2\Delta}(C_i)$. We may assume that there is a subset $I' \subseteq [2\Delta]$ of size $|I'|=\Delta$ such that the set of $y_i$ for $i\in I'$ are pairwise distinct. Indeed, if not, there must be some $\ell \in I$, and distinct $r,s,t\in I \subseteq [m]$ such that $y_{\ell}\in L_{2\Delta}(C_r)\cap L_{2\Delta}(C_s)\cap L_{2\Delta}(C_t)$. In particular we have $y_\ell\in L_{\Delta}(C_r)\cap L_{\Delta}(C_s)\cap L_1(C_t)$, a contradiction to \ref{item:periperhery_intersection1}. By reordering if necessary, we can further assume that $I' = \{1,2,\dots,\Delta\}$.

By \Cref{lem:deletedvertex_compsk/2}, there exists $r\in V(T)$ such that the connected components $S_1,\dots,S_q$ of $T-r$ each have order at most $\ceil{k/2}$, where $q\leq \Delta$. For $j\in [q]$, let $u_j$ denote the unique vertex in $N_T(r)\cap V(S_j)$. 

Now let us construct our embedding of $T$. For this purpose, we fix an embedding of $S_i$ into $H_i = G[(V(C_i)\setminus \{x,y_1,\dots,y_{q}\})\cup \{y_i\}]$ in such a way that $u_i$ is mapped to $y_i$ for each $i\in [q]$. Such embeddings exist by \Cref{fact:greedy_embedding} as $d_{H_i}(y_i)\geq \Delta$ and $\delta(H_i - y_i) \geq \delta(C_i) - q-1 \geq \ceil{k/2}$. Since the graphs $H_i$ have disjoint vertex sets and none of them contain $x$, we can combine these embeddings to an embedding of $T$ into $G$ by mapping $r$ to $x$, which is adjacent to all $y_i$.

The `in particular' part of the statement follows by the observation that $x\in L_1(C_i)$ implies $x\in L_1(L_{2\Delta}(C_i))$ for each $i\in [m]$.

\vspace{0.2cm}

\ref{item:periperhery_intersection3}:
Let $\delta$ be the output of \cref{lem:embedding_injection} when applied with $\eps/8$ and $\Delta$ and assume $k$ is sufficiently large so that $\Delta^{1004}(\lceil\log_{3/2}(k+1)\rceil+5)\leq \min \{\delta \floor{k/2},\eps^2k/2000\}$. Suppose to the contrary, without loss of generality $I = \{1,2\}$ is such that the set $\bigcap_{i\in I}  L_{\eps k}(C_i)$ has size at least $\eps k$. 
Let $X$ be a subset of $\bigcap_{i\in I}  L_{\eps k}(C_i)$
such that $|X| = \lfloor \eps k/3 \rfloor$.
Since for every $x\in X$ and $i\in I$ we have $|N_G(x)\cap V(C_i)| - |X| \geq \eps k - \eps k/3 \geq \eps k/2$, then every vertex in $X$ still has at least $\eps k/2$ neighbours in both $V(C_1)\setminus X$ and $V(C_2)\setminus X$.

Let $Y$ be the set of vertices in $V(C_1)\setminus X$ that have at least $\eps^2k/1600$ neighbours in $X$. We have
\begin{align*}
 \frac{\eps^2 k^2}{8} \leq \frac{\eps k }{2}|X| \leq e(X,V(C_1)\setminus X) \leq \frac{\eps^2 k}{1600}\vert V(C_1)\vert+ \frac{\eps k\vert Y\vert}{2}\leq \frac{\eps^2k^2}{16}+ \frac{\eps k\vert Y\vert}{2},
\end{align*}
from which it follows that $\vert Y\vert\geq \eps k/8$.

Fix a bipartition class in $T$, denoted by $A$. By \cref{cor:logk_subtrees} there exist edge-disjoint subtrees $S_0,S_1,\dots, S_\ell$ that cover $T$, where $\ell\leq \lceil\log_{3/2}(k+1)\rceil+5$ and $|S_0|=\floor{k/2}$, and for each $i \in [\ell]$, there is a unique vertex $w_i$ such that $V(S_0) \cap V(S_i) = \{w_i\}$. For each $i\in [\ell]$ define $W_i \subseteq V(S_i)$ to be the set of all vertices at distance at most $1000$ from $S_0$ if $w_i\in A$, or at distance at most $1001$ from $S_0$ if $w_i\notin A$. 

Consider the tree $S_0' = T[S_0 \cup \bigcup_{i
\in [\ell]}W_i]$ and note that $|\bigcup_{i\in [\ell]}W_i|\leq \ell \Delta^{1002}\leq \delta \floor{k/2} \leq \delta |S_0'|$. Let $L$ denote the set of leaves in $S_0'$ that are not leaves in $T$, and note that they cannot be leaves of $S_0$ either. Indeed, a leaf of $S_0$ is either a leaf of $T$ or it is adjacent to a vertex in $W_i$ for some $i\in [\ell]$. Furthermore, due to the way we defined the set $W_i$, we have $L\subset A$ and every vertex in $L$ is at distance at least 1000 from $S_0$.
We have $|S_0|\geq (1- \delta)|S_0'|$ and $|L|\leq \ell \Delta^{1002}$. 

Since $k$ is sufficiently large, we have $\delta(C_1\setminus X)\geq (1/2+\eps)k-|X|>(1+\eps/8)(k/2 +\ell\Delta^{1002})\geq(1+\eps/8)|S_0'|$. Applying \Cref{lem:embedding_injection} with $S_0$, $S_0'$, $C_1\setminus X$, $\eps/8$ and $|S_0'|$ playing the roles of $S$, $T$, $G$, $\eps$ and  $t$ respectively, there exists an embedding $\phi_0:S_0' \hookrightarrow C_1\setminus X$ such that $\phi_0(u) \in Y$ for each $u \in L$.

We show next that all remaining vertices of $T$ can be embedded into $X\cup V(C_2)$, maintaining compatibility with $\phi_0$. First, consider $U_1 = \left(\bigcup_{v\in L}N_T(v)\right)\setminus V(S_0')$, and note that $|U_1|\leq \Delta|L| \leq \Delta^{1003}\ell$. For each $u\in L$, since $\phi_0(u)\in Y$ and by choice of $Y$, it follows that
\begin{equation*}
    |N_G(\phi_0(u))\cap X|- \left|U_1\right| \geq \eps^{2}k/1600 - \Delta^{1003}\ell\geq \Delta.
\end{equation*}
So by a greedy argument there exists an embedding $\phi_1:T[U_1] \hookrightarrow X$ such that, for each $v\in U_1$ we have $\phi_0(u)\phi_1(v)\in E(G)$ for the unique $u\in L\cap N_T(v)$.

Now let $U_2$ denote the set of vertices in $T$ that have a neighbour in $U_1$, and that have not already been embedded under $\phi_0$ or $\phi_1$. Note that $|U_2|\leq \Delta|U_1|\leq \Delta^{1004}\ell$.
For each $v\in U_1$, since $\phi_1(v)\in X$ and by choice of $X$, we have
\begin{equation*}
    |N_G(\phi_1(v))\cap (V(C_2)\setminus X)|- |U_2|\geq \eps k/2 - \Delta^{1004}\ell \geq \Delta.
\end{equation*}
Again a greedy argument tells us that there exists an embedding $\phi_2:T[U_2]\hookrightarrow C_2\setminus X$ such that, for each $v\in U_2$ we have  $\phi_1(w)\phi_2(v)\in E(G)$ for every $w \in U_1 \cap N_T(v)$.

Finally, let $U_3 =V(T)\setminus (V(S_0')\cup U_1\cup U_2)$ be the set of vertices in $T$ yet to be embedded, noting that $|U_3|<k/2$. Observe that 
\begin{equation*}
    \delta(C_2) - |\phi_1(U_1) \cup \phi_2(U_2)| \geq (1/2+\eps)k - 2\Delta^{1004}\ell >k/2.
\end{equation*}
It follows that all vertices in $C_2$ have at least $k/2>|U_3|$ neighbours in $C_2$ that do not already belong to the image of $\phi_1\cup \phi_2$. Therefore we can greedily find an embedding $\phi_3:T[U_3]\hookrightarrow C_2$ such that for each $w\in U_3$ we have $\phi_2(v)\phi_3(w)\in E(G)$ for every $v\in U_2\cap N_T(w)$.

In particular, since $\phi_0(S_0')\subseteq V(C_1)\setminus X$, we know $\text{im}(\phi_3) \cap \text{im}(\phi_0) = \emptyset$. Also $\operatorname{im}(\phi_0)\cap\operatorname{im}(\phi_1)=\emptyset$, since $\operatorname{im}(\phi_1)\subseteq X$. It follows that the images of $\phi_0, \phi_1, \phi_2$ and $\phi_3$ are pairwise disjoint. 
So, taking the union $\phi_0\cup \phi_1\cup\phi_2\cup \phi_3$ defines an embedding of $T$ in $G$, a contradiction.
\end{proof}

For $s\in \N$ and a collection of disjoint subgraphs $C_1, \dots, C_m$ of $G$, we say that a vertex $v\in V(G)$ is \defn{$s$-split} (with respect to $(C_i)_{i\in [m]}$), if there exist distinct $i,j\in [m]$ such that $v\in L_{s}(C_i)\cap L_{s}(C_j)$. 

\begin{lemma}\label{lemma:split-neighbours}
    For all $\Delta\in \N$, $\eps \in (0,1)$ and sufficiently large $k \in \N$ the following holds. Let $T$ be a $k$-edge tree with $\Delta(T)\leq \Delta$. If $G$ is a graph with vertex-disjoint $(1/2+\eps,0,k)$-rich subgraphs $C_1, \dots, C_m$ such that $T$ cannot be embedded into $G$, then the set $S$ of $\eps k$-split vertices with respect to $(C_i)_{i\in [m]}$ satisfies $L_{\eps k}(S)=\emptyset$.
\end{lemma}

\begin{proof}
Set $\eps'=\eps/(2\binom{2\Delta}{2})$.  Since an $(1/2+\eps,0,k)$-rich graph is also $(1/2+\eps',0,k)$-rich, \cref{lem:intersectionproperties} may be applied with parameter $\eps'$. 
Let $k\in \N$ be sufficiently large in terms of $\Delta$ and $\eps$. Suppose towards a contradiction that $L_{\eps k}(S)$ is non-empty, meaning that there exists $x\in V(G)$ such that the size of $N=N_G(x)\cap S$ is at least $\eps k$. Let $I$ be the set of those $i\in [m]$ for which there exists $y\in N$ such that $y\in L_{\eps k}(C_i)$. Since $x\in \bigcap_{i\in I} L_1(L_{2\Delta}(C_i))$, by part \ref{item:periperhery_intersection2} of \Cref{lem:intersectionproperties} we must have that $\vert I \vert < 2\Delta$. Furthermore, because for each $y\in N$ there are distinct $i,j\in I$ such that $y\in L_{\eps k}(C_i)\cap L_{\eps k}(C_j)$, then there exist distinct $i,j\in I$ such that $\vert L_{\eps k}(C_i)\cap L_{\eps k}(C_j)\vert \geq |N|/\binom{2\Delta}{2}>\eps k/2\Delta^2$. 
 This is at least $2\eps'k$, while $L_{\eps k}(C_i)\subseteq L_{\eps'k}(C_i)$ and $L_{\eps k}(C_j)\subseteq L_{\eps'k}(C_j)$, contradicting \cref{lem:intersectionproperties}\ref{item:periperhery_intersection3} with parameter $\eps'$.
\end{proof}

We remark that \Cref{lem:intersectionproperties,lemma:split-neighbours} are stronger if $\eps$ is smaller. This means that we may apply \Cref{lemma:split-neighbours}, for example, to see that if $C_1,\dots,C_m$ are $(1/2+\eps_1,0,k)$-rich graphs then the set of $\eps_2 k$-split vertices $S$ satisfies $L_{\eps_3 k}(S)=\emptyset$ for any $\eps_1,\eps_2,\eps_3>0$ as long as $k$ is sufficiently large in terms of $\min(\eps_1,\eps_2,\eps_3)$. 

We say that $H\subseteq G$ is \defn{$t$-closed} if $\vert L_t(H)\setminus V(H)\vert< t$ and \defn{$t$-open} if $\vert L_t(H)\setminus V(H) \vert \geq t$. We are now ready to state our key lemma in this section.

\begin{lemma}\label{lemma:no-proper-bad}
For all $\Delta \in \N$ and $ \eps, \eta \in (0,1)$, there exists $\rho\in (0,1)$ such that the following holds for all sufficiently large $k\in \N$ \textup{(}in terms of $\Delta$, $\eps$, $\eta$ and $\rho$\textup{)} and all $a\geq1/2$. 
Let $T$ be a $k$-edge tree with $\Delta(T)\leq \Delta$. 
If $G$ is a graph with $\delta(G)\geq(a+\eps)k$ such that $T$ cannot be embedded into $G$, then $G$ contains vertex-disjoint $(a+\eps/4,\rho, k)$-rich subgraphs $C_1,\dots, C_m$ with the following properties: 
\begin{enumerate}[label=\upshape{(\Roman*)}, leftmargin =\widthof{110000}]
    \item no other collection of $(a+\eps/4,\rho, k)$-rich subgraphs covers more vertices than $(C_i)_{i\in [m]}$, \label{item:main1}
    \item all $C_i$ are $\eta k$-closed, i.e., $\vert L_{\eta k}(C_i)\setminus C_i\vert<\eta k$ for all $i\in [m]$, and \label{item:main2}
    \item for all $v\in V(G)$ there exist $i,j\in [m]$ such that $\vert N_G(v) \setminus (C_i \cup C_j)\vert < \eta k$. \label{item:main3}
\end{enumerate}
\end{lemma}

\begin{proof}
Let $\Delta\in \N$ and $\eps, \eta >0$, and note that properties \ref{item:main2} and \ref{item:main3} are monotonic in $\eta$, so without loss of generality we may assume that $\eta \ll\eps,\Delta^{-1}$. Additionally, let $\gamma=\eta/5\Delta$, $\delta=\gamma/800$, and $\rho=\delta^2/10^{10}$. Finally, let $k$ be sufficiently large to satisfy \Cref{lemma:good-components} with $\eps_{\ref{lemma:good-components}}=\eps/2$ and $\delta_{\ref{lemma:good-components}} = \delta$, \Cref{lem:intersectionproperties} with $\eps_{\ref{lem:intersectionproperties}}=\eps/4 $, and \Cref{lemma:split-neighbours} with $\eps_{\ref{lemma:split-neighbours}}=\min(\eps/4, \gamma, \delta\gamma/200)$. Now, suppose that $G$ is a graph with $\delta(G)\geq(a+\eps)k$ and that $T$ is a $k$-edge tree with $\Delta(T)\leq \Delta$ such that $T$ does not embed into $G$. We fix a family $(C_i)_{i\in [m]}$ of vertex-disjoint $(a+\eps/4,\rho,k)$-rich subgraphs whose union contains the maximum amount of vertices of $G$.

    Let $S$ be the set of $\gamma k$-split vertices with respect to $(C_i)_{i\in [m]}$, and define $I=\{i\in [m]: C_i \text{ is $\gamma k$-open}\}$ and $J=[m]\setminus I$. We say that a vertex is \defn{free}, if it is not contained in the set
    \begin{equation*}
        S\cup \bigcup_{i=1}^m V(C_i) \cup \bigcup_{j\in J} L_{\gamma k}(C_j).
    \end{equation*}
    \begin{claim}\label{claim:few-closed-comp-neighbours}
    If $v$ is free or $v\in V(C_i)\setminus S$ for some $i\in I$, then
    \begin{equation*}
        \left\vert N_G(v)\cap \left(\bigcup_{j\in J} L_{\gamma k}(C_j)\right)\right\vert <\frac{\eps k}{4}.
    \end{equation*}
    \end{claim}
    \begin{proofclaim}
    Let $K$ be the set of indices $j\in J$ such that $N_G(v)\cap L_{\gamma k}(C_j)\neq \emptyset$, i.e., such that $v\in L_1( L_{\gamma k}(C_j))$. As $\gamma k>2\Delta$ for sufficiently large $k$, we may infer from property \ref{item:periperhery_intersection2} of \Cref{lem:intersectionproperties} that the size of $K$ is less than $2\Delta$. Furthermore, if $v$ is free, then $v\notin L_{\gamma k}(C_j)$ for any $j\in J$ by definition, and if $v\in V(C_i)\setminus S$ for some $i\in I$, then $v\notin L_{\gamma k}(C_j)$ for any $j\in J$ as otherwise $v$ would be $\gamma k$-split --- so in both cases we establish $|N_G(v)\cap V(C_j)|\le \gamma k$. Since $\vert L_{\gamma k}(C_j)\setminus V(C_j)\vert <\gamma k$ for all $j\in J$ by definition of $J$, we have
    \begin{align*}
        \left\vert N_G(v)\cap \left(\bigcup_{j\in J} L_{\gamma k}(C_j)\right)\right\vert\leq         \sum_{j\in K}\left(\left\vert N_G(v)\cap V(C_j)\right\vert+\left\vert L_{\gamma k}(C_j)\setminus V(C_j)\right\vert\right)<4\Delta \gamma k < \frac{\eps k}{4},
    \end{align*}
    as desired.
    \end{proofclaim}

    \begin{claim}\label{claim:split-bound}
        For all $i\in I$, we have $\vert L_{\gamma k}(C_i)\cap S\vert \leq \delta k/2 \leq \delta|C_i|$.
    \end{claim}

    \begin{proofclaim}
    Since $|C_i|\geq \delta(C_i)\geq k/2$ for all $i\in [m]$, the second inequality is immediate, so it suffices to show that the first holds. Suppose to the contrary that $\vert L_{\gamma k}(C_i)\cap S\vert > \delta k/2$. By double counting of edges between $V(C_i)$ and $L_{\gamma k}(C_i)\cap S$, there exists $u\in V(C_i)$ with at least $\delta\gamma k/200$ neighbours in $L_{\gamma k}(C_i)\cap S$, and therefore in $S$, which by \Cref{lemma:split-neighbours} and choice of $k$ is not possible.
    \end{proofclaim}

    Let $F$ be the set of free vertices.
    \begin{claim}\label{claim:many-free-in-layer}
        For all $i\in I$, we have $\vert L_{\gamma k}(C_i) \cap F\vert \geq \vert L_{\gamma k}(C_i)\setminus V(C_i)\vert/2 \geq \gamma k/2$.
    \end{claim}

    \begin{proofclaim}
    By choice of $I$, we have $\vert L_{\gamma k}(C_i)\setminus V(C_i)\vert \geq \gamma k$, so that the second inequality holds for all $i \in I$. Next, we show that $L_{\gamma k}(C_i)\setminus V(C_i)\subset S\cup F$. Indeed, since $V(C_j)\subset L_{\gamma k}(C_j)$ for all $j\in [m]$, if $v\in L_{\gamma k}(C_i)\setminus V(C_i)$ is not $\gamma k$-split, then it cannot be in $\bigcup_{j\in [m]}V(C_j)$ or in $\bigcup_{j\in J} L_{\gamma k}(C_j)$. Hence, by definition, $v$ is free.

    By \Cref{claim:split-bound}, we obtain $\vert L_{\gamma k}(C_i) \cap F\vert \geq\vert L_{\gamma k}(C_i)\setminus V(C_i)\vert - \vert L_{\gamma k}(C_i)\cap S\vert \geq \vert L_{\gamma k}(C_i)\setminus V(C_i)\vert - \delta k/2 \geq \vert L_{\gamma k}(C_i)\setminus V(C_i)\vert/2$ as required. 
    \end{proofclaim}

    \begin{claim}\label{claim:no-free-vertices}
        $F$ is empty.
    \end{claim}

    \begin{proofclaim}
    Suppose towards a contradiction that $F$ is not empty. Since no vertex in $F$ is $\gamma k$-split by definition, $(L_{\gamma k}(C_i) \cap F)_{i\in I}$ is a disjoint family of subsets. Recalling that $\vert C_i\vert < 100k$ for all $i\in [m]$ by definition of a rich subgraph, we observe that by \Cref{claim:many-free-in-layer},
    \begin{equation}\label{eq:open-comp-bound}
    \left\vert \bigcup_{i\in I} V(C_i)\right\vert \leq \sum_{i\in I} 100k \leq \frac{200}{\gamma}\sum_{i\in I} \vert L_{\gamma k}(C_i) \cap F\vert \leq \frac{200\vert F\vert}{\gamma}.
    \end{equation}
    Consider now the subgraph $H$ that $G$ induces on 
    \begin{equation*}
        V(H)=F\cup \bigcup_{i\in I} V(C_i)\setminus S.
    \end{equation*}
    By \Cref{claim:split-bound}, we have
    \begin{equation}\label{eq:key-lemma-lower-bound-H}
        \vert V(H)\vert \geq \vert F\vert + \left\vert \bigcup_{i\in I} V(C_i)\setminus S\right\vert \geq \vert F\vert + (1-\delta)\left\vert \bigcup_{i\in I} V(C_i)\right\vert.
    \end{equation}
    In addition, we will show now that the minimum degree in $H$ is at least $(a+\eps/2)k$. The minimum degree of $G$ is at least $(a+\eps)k$ by assumption, so it suffices to show that for each $v\in V(H)$, we have $\vert N_G(v)\setminus V(H)\vert \leq \eps k/2$. This is the case because by \Cref{claim:few-closed-comp-neighbours}, no vertex in $V(H)$ has more than $\eps k/4$ neighbours in $\bigcup_{j\in J} L_{\gamma k}(C_j)$, and by \Cref{lemma:split-neighbours}, no vertex in $V(G)$ has more than $\eps k/4$ neighbours in $S$.

    Applying \Cref{lemma:good-components} to $H$, we obtain $(a+\eps/4,\rho,k)$-rich subgraphs $D_1,\dots,D_{m'}$ that cover all but $\delta\vert H\vert$ vertices in $H$. Together with the $(a+\eps/4,\rho,k)$-rich subgraphs $(C_j)_{j\in J}$ this gives a collection $\mathcal{C}$ of disjoint $(a+\eps/4,\rho,k)$-rich subgraphs. The number of vertices covered by $\mathcal{C}$ is
    \begin{align*}
        \left\vert \bigcup_{j\in J} V(C_j)\cup \bigcup_{i=1}^{m'} V(D_i)\right\vert &\geq \left\vert \bigcup_{j\in J} V(C_j)\right\vert +(1-\delta)\left\vert V(H)\right\vert \\
        \overset{\eqref{eq:key-lemma-lower-bound-H}}&{\geq}\left\vert \bigcup_{j=1}^m V(C_j)\right\vert-2\delta\left\vert \bigcup_{i\in I} V(C_i)\right\vert +(1-\delta) \vert F\vert.
    \end{align*}
    By \Cref{eq:open-comp-bound}, the latter is at least
    \begin{equation*}
        \left\vert \bigcup_{j=1}^m V(C_j)\right\vert-\frac{400\delta}{\gamma}\vert F\vert +(1-\delta) \vert F\vert,
    \end{equation*}
    which, as $\delta= \gamma/800$ and $\vert F\vert >0$, is more than the number of vertices that were covered by the original collection of rich subgraphs $(C_i)_{i\in [m]}$, contradicting its maximality.
    \end{proofclaim}
    Combining \cref{claim:many-free-in-layer,claim:no-free-vertices}, we deduce that $I=\emptyset$ as otherwise $\gamma k \leq 0$. Since $\gamma <\eta$ and $L_{\eta k}(C_i) \subseteq L_{\gamma k}(C_i)$, this shows that $(C_i)_{i\in [m]}$ has property \ref{item:main2}. 
    
    Since there are no free vertices and $J = [m]$, it follows that 
    \begin{equation}\label{eq:all_closed_layers}
        V(G)=\bigcup_{j\in J} L_{\gamma k}(C_j).
    \end{equation}
We now complete the proof by showing that for every $v\in V(G)$ there exist $i,j\in [m]$ such that $|N_G(v)\setminus (C_i\cup C_j)|<\eta k$. Fix $v\in V(G)$ and consider the set $M=\{\ell\in [m]: v\in L_1(L_{\gamma k}(C_\ell))\}$. By part \ref{item:periperhery_intersection2} of \Cref{lem:intersectionproperties}, $M$ has size less than $2\Delta$. If $M=\{i\}$ is a singleton, then by \eqref{eq:all_closed_layers} every neighbour of $v$ lies in $L_{\gamma k}(C_i)$, and hence $|N_G(v)\setminus C_i|\leq |L_{\gamma k}(C_i)\setminus C_i|<\gamma k<\eta k.$ Thus property \ref{item:main3} holds in this case by taking $j=i$. We may therefore assume that $|M|\geq2$, and fix the two indices $i,j\in M$ for which $\vert N_G(v)\cap V(C_i)\vert$ and $\vert N_G(v) \cap V(C_j)\vert$ are largest and second largest, respectively, settling ties arbitrarily. If $\vert N_G(v) \cap V(C_j)\vert\geq \Delta$, it follows from \ref{item:periperhery_intersection1} of \Cref{lem:intersectionproperties} that $N_G(v)\cap V(C_\ell)=\emptyset$ for all $\ell \in M\setminus \{i,j\}$. If $\vert N_G(v) \cap V(C_j)\vert< \Delta$, we have $\vert N_G(v) \cap V(C_\ell)\vert< \Delta$ for all $\ell \in M\setminus\{i\}$ by choice of $j$. Either way, we may conclude that 
    \begin{equation}\label{eq:few-component-neighboursv2}
        \left\vert N_G(v)\cap \bigcup_{\ell\in [m]\setminus\{i,j\}} V(C_\ell)\right\vert=\left\vert N_G(v)\cap \bigcup_{\ell\in M\setminus\{i,j\}} V(C_\ell)\right\vert < 2\Delta^2.
    \end{equation}

    Furthermore, by choice of $J$, all components are $\gamma k$-closed, so we have
    \begin{equation}\label{eq:few-peripheral-neighboursv2}
        \left\vert  N_G(v)\cap \left( \bigcup_{\ell\in M} (L_{\gamma k}(C_\ell)\setminus V(C_\ell))\right)\right\vert \leq 2\Delta\gamma k.
    \end{equation}

Combining \eqref{eq:all_closed_layers}, \eqref{eq:few-component-neighboursv2}, and \eqref{eq:few-peripheral-neighboursv2}, we obtain
    \begin{align*}
        \left\vert N_G(v)\setminus(C_i\cup C_j)\right\vert &\leq   \left\vert N_G(v)\cap \left(\bigcup_{\ell\in [m]\setminus \{i,j\}} V(C_\ell) \right)\right\vert+\left\vert N_G(v)\cap \left(\bigcup_{\ell\in [m]} L_{\gamma k}(C_\ell)\setminus V(C_\ell)\right)\right\vert\\
        &\leq 2\Delta^2+2\Delta\gamma k,
    \end{align*}
    but as we chose $\gamma$ to be $\eta/5\Delta$, it follows that the above is at most $\eta k$, as desired for \ref{item:main3}. This proves the lemma. 
\end{proof}

\subsection[Applying the sparse to dense lemma]{Proof of \cref{thm:approx_k/2}}

We will now demonstrate the usefulness of  \cref{lemma:no-proper-bad} by providing a short proof of \cref{thm:approx_k/2}. The idea is to find a dense subgraph in $G$ that still has minimum and maximum degree asymptotically exceeding the respective thresholds from the corresponding conjecture (\cref{conj:boundedk/2}). Applying the previously known approximate result for dense graphs (\cref{thm:2k-1/Delta_dense}) will give our desired conclusion. For this proof, we only require part \ref{item:main3} of \cref{lemma:no-proper-bad}, and do not use cut-density within the rich subgraphs. We remark that an identical argument would also provide approximate versions of \cref{conj:2k/3,conj:alpha} for large bounded degree trees, but we omit these here since Sections~\ref{sec:alpha}~and~\ref{sec:2/3} are dedicated to proving the stronger versions, \cref{thm:exact2/3,thm:mainalpha-exact} respectively.

\begin{proof}[Proof of \Cref{thm:approx_k/2} using \cref{lemma:no-proper-bad}]
    Let $\eps>0$, $\Delta \in \mathbb{N}$, and $k$ be sufficiently large.
    Let $T$ be a $k$-edge tree with $\Delta(T)\leq \Delta$. Suppose that $G$ is a graph with $\delta(G)\geq(1+\eps)k/2$ and $\Delta(G)\geq 2(1 - 1/\Delta+\eps)k$ such that $T$ does not embed in $G$. Let $x\in V(G)$ be chosen such that $d_G(x) = \Delta(G)$. Without loss of generality we may assume that $\eps<1/50$. 
    
    Apply \cref{lemma:no-proper-bad} to $G$ with $\eps_{\ref{lemma:no-proper-bad}}=\eps/2$, $\eta_{\ref{lemma:no-proper-bad}}=\eps$ and $a_{\ref{lemma:no-proper-bad}} = 1/2$ to obtain a collection of vertex-disjoint $(1/2+\eps/8,0,k)$-rich subgraphs $C_1,\dots, C_m$. By \ref{item:main3}, there exist $i,j\in[m]$ such that $\vert N_G(x)\setminus(C_i\cup C_j)\vert<\eps k$. If $i=j$, let $G'=G[V(C_i)\cup\{x\}]$. Then $\Delta(G')\geq \Delta(G)-\eps k>2(1-1/\Delta+\eps/4)k,$ while $\delta(G')\geq(1+\eps/4)k/2$ and $|G'|<100k+1<4\eps^{-1}k$. Hence \cref{thm:2k-1/Delta_dense}, with parameter $\eps/4$, gives a copy of $T$, a contradiction. Thus $i\neq j$, and by reordering we may assume that $i=1$ and $j=2$.  Consider the subgraph $G' = G[V(C_1)\cup V(C_2) \cup \{x\}]$. We have $\Delta(G') \geq |N_G(x)\cap (C_1 \cup C_2)|\geq \Delta(G)-\eps k >  2(1-1/\Delta +\eps/4)k$. Furthermore, by definition of a rich subgraph, $\delta(C_1)\geq (1+\eps/4)k/2$ and $\delta(C_2)\geq (1+\eps/4)k/2$, so that $\delta(G')\geq (1+\eps/4)k/2$. Finally, we have $|G'|\leq |C_1|+|C_2|+1< 200k <4\eps^{-1}k$. Therefore $G'$ satisfies the conditions of \cref{thm:2k-1/Delta_dense} with $\eps_{\ref{thm:2k-1/Delta_dense}} = \eps/4$, and applying this result tells us that $T$ embeds in $G'$, and thus in $G$, a contradiction.
\end{proof}

\section[Mixed degree combinations]{Mixed degree combinations: Proof of \cref{thm:mainalpha-exact}}\label{sec:alpha}

In order to prove \cref{conj:alpha} exactly, we first prove the following lemma which has a very similar flavour to that of \cref{thm:mainalpha-exact}, but that allows a slightly lower maximum degree bound and has an asymptotic relaxation.

\begin{lemma}\label{thm:alpha_weakened}
     For all $\Delta\in \N$, $\alpha \in (0,1/3)$ and $\beta<\min\{2/(3\Delta),1/3-\alpha\}$ the following holds. For all $\eps>0$ and sufficiently large $k\in \N$ \textup{(}in terms of $\Delta$, $\eps$ and $\beta$\textup{)}, if $G$ is a graph with $\delta(G)\geq (1+\alpha+\eps)k/2$ and $\Delta(G)\geq 2(1-\alpha-\beta+\eps)k$, then $G$ contains a copy of every $k$-edge tree $T$ with $\Delta(T)\leq \Delta$.
\end{lemma}

Now, let us show that \cref{thm:mainalpha-exact} is easily deduced from \cref{thm:alpha_weakened}.

\begin{proof}[Proof of \Cref{thm:mainalpha-exact} using \Cref{thm:alpha_weakened}]
     Let $\alpha\in (0,1/3)$ and $\Delta\in \N$. Let $\beta>0$ be such that $\beta<\min\{2/(3\Delta),1/3-\alpha,2\alpha\}$. 
    Define $\alpha' = \alpha -\beta/2$, noting that $\alpha'\in (0,1/3)$. Let $k$ be sufficiently large with respect to all other parameters. 
    Suppose $G$ is a graph satisfying $\delta(G)\geq(1+\alpha)k/2$ and $\Delta(G)\geq 2(1-\alpha)k$. We have $\delta(G)\geq (1+\alpha)k/2 = (1+\alpha'+\beta/2)k/2$ and $\Delta(G)\geq 2(1-\alpha)k 
    = 2(1-\alpha'-\beta+\beta/2)k$. Since $\beta<\min\{2/(3\Delta),1/3-\alpha'\}$, we may apply \cref{thm:alpha_weakened} with $\alpha_{\ref{thm:alpha_weakened}} = \alpha'$ and $\eps_{\ref{thm:alpha_weakened}} = \beta/2$ to deduce that $G$ contains a copy of every $k$-edge tree $T$ with $\Delta(T)\leq \Delta$, as desired.
\end{proof}

Our method for proving \cref{thm:alpha_weakened} again begins with an application of \cref{lemma:no-proper-bad} to a graph $G$ satisfying the minimum and maximum degree conditions. For a highest degree vertex $x$ in $G$ this yields two rich subgraphs in $G$ (as in \cref{defn:good_subgraph}), each with minimum degree asymptotically above $(1+\alpha)k/2$, in which $x$ has almost all of its neighbours. Let us denote by $C_1$ and $C_2$ these rich subgraphs. We consider a reduced graph of $C_1 \cup C_2$, and use tools of Besomi, Pavez-Sign\'e and Stein (see \cref{subsec:BPS}) to gain insight into the structure of this reduced graph. We again only require part \ref{item:main3} from the conclusions of \cref{lemma:no-proper-bad}, but unlike in the proof of \cref{thm:approx_k/2}, we also make use of the positive cut-density of $C_1$ and $C_2$ in order to use the connectivity within the reduced graph.

\begin{proof}[Proof of \cref{thm:alpha_weakened}]
     Let $\Delta\in \N$, $\alpha\in (0,1/3)$, $\beta <\min\{2/(3\Delta),1/3 -\alpha\}$ and $\eps >0$. Note that we may assume that $\eps\leq 1/100$ since the lemma only gets stronger as $\eps$ decreases. 
    Let $\rho>0$ be the output of \cref{lemma:no-proper-bad} when applied with $\eps_{\ref{lemma:no-proper-bad}} = \eta_{\ref{lemma:no-proper-bad}}= \eps/2$, and let $\gamma>0$ be sufficiently small in terms of $\eps$ and $\rho$.  
    Let $M_0$ and $N_0$ be the outputs of \cref{lemma:regularity-degree} when applied with $\gamma_{\ref{lemma:regularity-degree}}=\gamma$ and $(m_0)_{\ref{lemma:regularity-degree}}=\gamma\inv$. Let $k$ be sufficiently large in terms of all other parameters.

     Let $T$ be a $k$-edge tree with $\Delta(T)\leq \Delta$. Let $G$ be a graph satisfying the assumptions of \cref{thm:alpha_weakened}, but such that $G$ does not contain a copy of $T$. Since $\delta(G)\geq (1+\alpha+\eps)k/2 = (1/2+\alpha/2+\eps/2)k$, we may apply \cref{lemma:no-proper-bad} to $G$ with $\eps_{\ref{lemma:no-proper-bad}} = \eps/2$, $\eta_{\ref{lemma:no-proper-bad}}= \eps/2$ and $a_{\ref{lemma:no-proper-bad}} = 1/2+\alpha/2$ to obtain a collection of vertex-disjoint $(1/2+\alpha/2+\eps/8, \rho,k)$-rich subgraphs $C_1,\dots,C_m$ such that by \ref{item:main3}, for every $v\in V(G)$ there exist $i,j\in [m]$ for which $|N_G(v)\setminus (C_i\cup C_j)|< \eps k/2$. Let $x\in V(G)$ satisfy $d_G(x) = \Delta(G)\geq 2(1-\alpha-\beta+\eps)k$. By \ref{item:main3}, choose $i,j\in[m]$ such that $|N_G(x)\setminus(C_i\cup C_j)|<\eps k/2$. If $i=j$, then
$|N_G(x)\cap C_i|\geq 2(1-\alpha-\beta+3\eps/4)k>(1+\eps/16)4k/3$, using $\beta<1/3-\alpha$. Since $C_i$ is $\rho$-cut-dense and $\delta(C_i)\geq (1+\eps/8)k/2$, \cref{lem:single-rich-BPS} gives a copy of $T$, a contradiction. Hence $i\neq j$, and after relabelling we may assume $i=1,j=2$.

    Let $H$ be the subgraph of $G$ obtained by taking the disjoint union of $C_1 -x$ and $C_2-x$, not including any edges that lie between these two subgraphs. Hence, $H$ has two connected components, each of which has minimum degree at least $(1/2+\alpha/2+\eps/8)k - 1$, size less than $100k$, and by \cref{fact:cutdense_subgraph}, is $\rho/2$-cut-dense. Applying \cref{lemma:regularity-degree} to $H$ with $\gamma_{\ref{lemma:regularity-degree}} = \gamma$ and $\eta_{\ref{lemma:regularity-degree}} = 5\sqrt{\gamma}$, we obtain a subgraph $H'$ of $H$ such that $|H'|>(1-\gamma)|H|$ and $d_{H'}(v)>d_{H}(v) - (\gamma +5\sqrt{\gamma})|H|\geq (1/2+\alpha/2+\eps/16)k$ for all $v\in V(H')$, and a $(\gamma,5\sqrt{\gamma})$-partition $\cP$ of $V(H')$ with $|\cP|\leq M_0$. 
    Letting $C_i' = H' \cap C_i$ for both $i\in \{1,2\}$, observe that $e_{H'}(C_1',C_2') = 0$, $\delta(C_i')\geq (1/2+\alpha/2+\eps/16)k$ and $x\notin V(C_i')$. We make the following claim.

    \begin{claim}
        $\mathcal{P}$ refines the partition $V(H') = V(C_1')\cup V(C_2')$
    \end{claim}
    
   \begin{proofclaim}
       Suppose otherwise, that there exist $y_1\in V(C_1')$ and $y_2 \in V(C_2')$ that belong to the same part $P_0\in \mathcal{P}$. Now suppose for each $P\in \mathcal{P}$, $y_1$ has at most $|P|/400$ neighbours in $P$. 
    Noting that $|H'|<200k$, this yields a contradiction to $\delta(H') >k/2$ since
    \begin{equation*}
        d_{H'}(y_1) = \sum_{P\in \mathcal{P}}|N_{H'}(y_1)\cap V(P)| \leq \frac{1}{400}\sum_{ P\in \mathcal{P}}|P| = \frac{|H'|}{400} < \frac{k}{2}.
    \end{equation*}
So, there must exist $P_1\in \mathcal{P}$ such that $|N(y_1)\cap P_1|\geq |P_1|/400$. Since each vertex class in a regular partition is an independent set, $P_1\neq P_0$, and since $e(P_0,P_1)>0$, we have that $d(P_0,P_1)>5\sqrt{\gamma}$.  Similarly, there exists $P_2\in \mathcal{P} \setminus \{P_0\}$ such that $|N(y_2)\cap P_2|\geq |P_2|/400$ and $d(P_0,P_2)>5\sqrt{\gamma}$. 

Note that $|N(y_1)\cap P_1|\geq \gamma |P_1|$. Since $(P_0,P_1)$ is a $\gamma$-regular pair, we have $d(P_0, N(y_1)\cap P_1)> 5\sqrt{\gamma} - \gamma> 4\sqrt{\gamma}$, from which it follows that the set $Q_0 \subseteq P_0$ of vertices $v\in P_0$ with at least one neighbour in $N(y_1)\cap P_1$, has size at least $4\sqrt{\gamma}|P_0|$. Now, every vertex in $Q_0$ is connected to $y_1$ in $H'$ via a path of length two, and so $Q_0\subseteq V(C_1')$. Since $(P_0,P_2)$ is a $\gamma$-regular pair, we have $d(Q_0, N(y_2)\cap P_2)>4\sqrt{\gamma} >0$. This contradicts the fact that $y_2\in V(C_2')$ and $e_{H'}(C_1',C_2') = 0$, therefore the claim holds.
\end{proofclaim} 

For $i\in \{1,2\}$, let $\mathcal{P}_i$ consist of the set of vertex classes $P\in \mathcal{P}$ for which $P\cap V(C_i') \neq \emptyset$. It follows from the claim that $\mathcal{P}_1$ and $\mathcal{P}_2$ are disjoint, that $\mathcal{P} = \mathcal{P}_1\cup \mathcal{P}_2$, and that $\mathcal{P}_i$ is a $(\gamma,5\sqrt{\gamma})$-regular partition of $C_i'$.

Let $R$ be the reduced graph of $H'$ with respect to the partition $\cP$, and for $i\in \{1,2\}$ let $R_i$ be the reduced graph of $C_i'$ with respect to the partition $\cP_i$. Since $|H|<200k$ while $|C_i-x|>k/2$, every surviving vertex of $C_i-x$ loses fewer than $400(\gamma+5\sqrt{\gamma})|C_i-x|$ incident edges on passing to $C_i'$. Thus, for $\gamma$ sufficiently small in terms of $\rho$, \cref{fact:cutdense_subgraph} implies that $C_i'$ is $\rho/4$-cut-dense. 
 By \cref{fact:cutdense_reduced_connected}, each $R_i$ is connected. Therefore, by definition of the reduced graph, $R$ is composed of precisely the connected components $R_1$ and $R_2$. Note that $|R|\leq M_0$. Moreover by \cref{fact:reduced-min-degree}, each of $R$, $R_1$ and $R_2$ has minimum degree at least $(1/2+\alpha/2+\eps/32)k \cdot \frac{|R|}{|H'|}$. Note also, that $\frac{|H'|}{|R|} = \frac{|C_1'|}{|R_1|}= \frac{|C_2'|}{|R_2|}$ since this is precisely the number of vertices in every element of the partition $\cP$. Since $\delta(R_1)\geq k/2 \cdot |R_1|/|C_1'|$ and $|C_1|\leq 100k$, by \cref{thm:mindeg_diameter} we have $\diam(R_1)\leq 600$.

Let $H'_x$ be the subgraph obtained by adding $x$ to $H'$ as well as all edges that lie between $x$ and $V(H')$ in the original graph $G$. The number of neighbours of $x$ in $V(H')$ is at least $|N(x)\cap (C_1\cup C_2)| - |H\setminus H'| 
\geq 2(1-\alpha-\beta+\eps/2)k$, so we have $\Delta(H_x')\geq 2(1-\alpha-\beta+\eps/2)k \geq (1+\eps/2)4k/3$, since $\beta<1/3-\alpha$. Furthermore, $\delta(H_x')\geq (1+\alpha+\eps/16)k/2$ and $|H'_x|<200k$.

Now, the conditions of \cref{thm:BPS_combined} are satisfied with $G_{\ref{thm:BPS_combined}} = H_x'$, $\eps_{\ref{thm:BPS_combined}} = \eps/16$, $\gamma_{\ref{thm:BPS_combined}} = \gamma$ and $R_{\ref{thm:BPS_combined}} = R$, and therefore $R$ satisfies properties \ref{item:BPS1} through \ref{item:BPS4}. The two components $S_1$ and $S_2$ from \ref{item:BPS3} must in fact be $R_1$ and $R_2$, since we observed earlier that these are the only two components of $R$. Assume without loss of generality that $x$ has at least as many neighbours in $C_1'=\bigcup V(R_1)$ as it has in $C_2'=\bigcup V(R_2)$. Then we have the following.  
\begin{itemize}
    \item By \ref{item:BPS3}, $x$ $\sqrt{\gamma}$-sees both $R_1$ and $R_2$. In other words, $x$ has at least $\sqrt{\gamma}|C_1'|$ neighbours in $C_1'$ and at least $\sqrt{\gamma}|C_2'|$ neighbours in $C_2'$.
    \item By \ref{item:BPS4}, the component $R_1$ is bipartite with parts $(X,Y)$ such that$|X|\geq (1+\sqrt[4]{\gamma})\frac{2k}{3}\cdot \frac{|R|}{|H_x'|}$, and  $x$ has  no neighbours in $\bigcup V(Y)$.  Since $|H_x'|=|H'|+1$ and $|H'|=\Omega(k)$, for sufficiently large $k$ this implies $|X|\geq (1+100\sqrt{\gamma})\frac{2k}{3}\cdot\frac{|R|}{|H'|}$. 
\end{itemize}

Now let us consider the structure of the tree $T$. Let $r\in V(T)$ be obtained by \cref{cor:partition_subforests_T-v} so that both conclusions \ref{item:Fsize1} and \ref{item:Fsize2} are satisfied. Consider the sets $\odd_T(r)$ and $\even_T(r)$.

    \textbf{Case 1: $|\odd_T(r)|\geq |\even_T(r)|$.}
    
    By \cref{cor:partition_subforests_T-v}\ref{item:Fsize1}, there exists a partition of $T-r$ into vertex-disjoint subforests $F_1$ and $F_2$ such that $k/2\leq |F_1|\leq \floor{2k/3}$. Let us construct an embedding $\phi$ of $T$ in $G$. First, let $r$ be mapped to $x$. 

        We now embed $F_1$ appropriately into $C_1'$. 
    Note that $ V(F_1)\cap \odd_T(r)$ and $ V(F_1)\cap \even_T(r)$ form a bipartition of the forest $F_1$, and let $k_1$ and $k_2$ denote the respective sizes of these sets, so that $k_1+k_2 = |F_1|$.  Clearly we have $|X|\geq (1+100\sqrt{\gamma})k_1 \cdot |R_1|/|C_1'|$. By assumption of Case 1, $k_2\leq k/2$ and so using the minimum degree condition of $R_1$, we have $d_{R_1}(w) \geq (1+100\sqrt{\gamma})k_2\cdot \frac{|R_1|}{|C_1'|}$ for every $w\in X$. Therefore, the conditions of \cref{lem:BPSprescribed_forest_embedding_bipartite} are satisfied with
    $G_{\ref{lem:BPSprescribed_forest_embedding_bipartite}} = C_1'$ and $R_{\ref{lem:BPSprescribed_forest_embedding_bipartite}} = R_1$. Applying the conclusion of this lemma with $F_{\ref{lem:BPSprescribed_forest_embedding_bipartite}} = F_1$, $(A_1)_{\ref{lem:BPSprescribed_forest_embedding_bipartite}} = V(F_1)\cap \odd_T(r)$,
    and $(A_2)_{\ref{lem:BPSprescribed_forest_embedding_bipartite}} = V(F_1)\cap \even_T(r)$, we have that $F_1$ can be embedded into $C_1'$ such that $V(F_1)\cap \odd_T(r)$ goes into $\bigcup V(X)$ and $V(F_1)\cap \even_T(r)$ goes into $\bigcup V(Y)$. Moreover, since in $G$, $x$ has at least $(1-\alpha-\beta+\eps/4)k$ neighbours in $C_1'$, and at least $(1-\alpha-\beta+\eps/8)k$ of these lie in $\bigcup V(X)$, considering the set $V(F_1)\cap N_T(r)$ to be the `roots' of the forest $F_1$, \cref{lem:BPSprescribed_forest_embedding_bipartite} also tells us that we can further ensure our embedding of $F_1$ maps the vertices in $V(F_1)\cap N_T(r)$ into the set $\bigcup V(X)\cap N_G(x)$. This ensures the embedding is compatible with $r$ being mapped to $x$.

    It remains to embed $F_2$. The graph $T[V(F_2)\cup\{r\}]$ is a tree rooted at $r$. Since $\delta(C_2')\geq k/2\geq |F_2|$ and $x$ has at least $\sqrt{\gamma}|C_2'|\geq \Delta$ neighbours in $C_2'$, \cref{fact:greedy_embedding} embeds this tree into $G[V(C_2')\cup\{x\}]$ with $r$ mapped to $x$. Thus, altogether we have found a valid embedding of $T$ in $G$, a contradiction. 
    
    \textbf{Case 2: $|\odd_T(r)|< |\even_T(r)|$.}

    By choice of $r$ and \cref{cor:partition_subforests_T-v}\ref{item:Fsize2}, all components of $T-r$ have at most $\lceil k/2 \rceil$ vertices. 
    Let $B$ be a component of $T-r$ of largest size.
    First suppose that $|B|>\alpha k$. In this case, $T-B$ is a tree on less than $k-\alpha k+1$ vertices, so that by \cref{fact:tree_part_sizes}, the larger bipartition class of $T-B$ has at most $(1-1/\Delta)(1-\alpha +1/k)k$ vertices.
    Thus, denoting the sizes of bipartition classes $ (V(T-B)\cap \even_T(r))\cup \{r\}$ and $ V(T-B)\cap \odd_T(r)$ of $T-B$, by $\ell_1$ and $\ell_2$ respectively, we have $\ell_1\leq (1-1/\Delta)(1-\alpha +1/k)k$, as well as $\ell_2\leq k/2$ by assumption of Case 2.
    
       Similarly as in Case 1, we have $d_{R_1}(w) \geq (1+100\sqrt{\gamma})\ell_2\cdot \frac{|R_1|}{|C_1'|}$ for every $w\in X$. Since more than half of all neighbours of $x$ in $H_x'$ lie in $C_1'$, and since no neighbours of $x$ in $C_1'$ lie outside $\bigcup V(X)$, and since $\beta<2/(3\Delta) <(1-\alpha)/\Delta$, it follows that
    \begin{align*}
        |X|&\geq \bigl(|N_G(x)\cap V(C_1')|\bigr)\cdot \tfrac{|R_1|}{|C_1'|}\geq (1-\alpha-\beta+\eps/4)k \cdot \tfrac{|R_1|}{|C_1'|} \\
        &\geq
        (1+100\sqrt{\gamma})(1-1/\Delta)(1-\alpha +1/k)k \cdot \tfrac{|R_1|}{|C_1'|} 
        \geq (1+100\sqrt{\gamma})\ell_1 \cdot \tfrac{|R_1|}{|C_1'|}.
    \end{align*}
Again the conditions of \cref{lem:BPSprescribed_forest_embedding_bipartite} are satisfied, now with
$G_{\ref{lem:BPSprescribed_forest_embedding_bipartite}} = C_1'$, $R_{\ref{lem:BPSprescribed_forest_embedding_bipartite}} = R_1$, $(k_1)_{\ref{lem:BPSprescribed_forest_embedding_bipartite}} = \ell_1$ and $(k_2)_{\ref{lem:BPSprescribed_forest_embedding_bipartite}} = \ell_2$. Applying the conclusion of this lemma with $F_{\ref{lem:BPSprescribed_forest_embedding_bipartite}} = T-B$, $(A_1)_{\ref{lem:BPSprescribed_forest_embedding_bipartite}} = (V(T-B)\cap \even_T(r))\cup \{r\}$,
    and $(A_2)_{\ref{lem:BPSprescribed_forest_embedding_bipartite}} = V(T-B)\cap \odd_T(r)$, we have that $T-B$ can be embedded into $C_1'$ such that $(V(T-B)\cap \even_T(r))\cup \{r\}$ goes into $\bigcup V(X)$ and $V(T-B)\cap \odd_T(r)$ goes into $\bigcup V(Y)$. Moreover, 
    the lemma also allows us to assume that $r$ is mapped into the set $\bigcup V(X)\cap N_G(x)$.
    
   Now, since $B$ is a single component of $T-r$, there is a unique neighbour $u$ of $r$ in $B$. By \cref{fact:greedy_embedding}, we embed $B$ such that $u$ is mapped to $x$ and $B-u$ is mapped into $C_2'$, 
   noting that this is compatible with the placement of $r$ and thus completes an embedding of $T$ in $G$.

    So, suppose otherwise that $|B|\leq \alpha k$, in particular this means that all components of $T-r$ have at most $\alpha k$ vertices. We claim that there exists a partition of $T-r$ into two vertex-disjoint subforests $F_1 $ and $F_2$ such that $|F_2|\leq |F_1|\leq (1+\alpha)k/2$. Indeed, let us group the components of $T-r$ into two disjoint subforests $F_1$ and $F_2$, such that the size of $F_1$ is maximised subject to $|F_1|\leq (1+\alpha)k/2$. If additionally we have $|F_1|<(1-\alpha)k/2$, then every component of $F_2$ must have size greater than $\alpha k$, otherwise we would add the component to $F_1$, contradicting the maximality property. Since this is not possible, then we must have $|F_2| = k - |F_1| \leq (1+\alpha)k/2$, as desired.

     Now map $r$ to $x$. For each $i\in\{1,2\}$, the graph $T[V(F_i)\cup\{r\}]$ is a tree rooted at $r$. We have $\delta(C_i)\geq(1+\alpha)k/2\geq |F_i|$, while \ref{item:BPS3} gives $d_G(x,C_i')\geq\sqrt\gamma|C_i'|\geq\Delta$ for sufficiently large $k$. Hence \cref{fact:greedy_embedding} embeds $T[V(F_i)\cup\{r\}]$ into $G[V(C_i)\cup\{x\}]$ with $r$ mapped to $x$. Since $C_1$ and $C_2$ are vertex-disjoint, the two embeddings meet only at $x$ and together form a copy of $T$, a contradiction. This completes Case 2, and the proof of the lemma.
\end{proof}

\section[Exact 2/3]{Exact \boldmath{$2/3$}: Proof of \cref{thm:exact2/3}}\label{sec:2/3}

In this section, we prove \cref{thm:exact2/3}. With the help of \cref{lemma:no-proper-bad}, starting with a graph $G$ satisfying $\delta(G)\geq 2k/3$ and $\Delta(G)\geq k$, we find a subgraph in $G$ that looks somewhat like the extremal example in \cref{fig:extremal2k/3}. Similarly to previous proofs, we take a highest degree vertex $x$ in $G$ and find two rich subgraphs in $G$ (as per \cref{defn:good_subgraph}), and call them $C_1$ and $C_2$, in which $x$ has almost all of its neighbours. In particular, the subgraph induced on $\{x\}\cup V(C_1)\cup V(C_2)$ is composed of two vertex-disjoint subgraphs $C_1$ and $C_2$ that both have at most $100k$ vertices and have minimum degree at least $(2/3-\eps)k$ for some small $\eps>0$, as well as an additional vertex $x$ that has at least $(1-\eps)k$ neighbours in $V(C_1)\cup V(C_2)$. 

Using this approach allows us to focus our attention only on host graphs of this form, and we tackle the problem in this setting by applying a sequence of stability arguments. There are still quite a few cases to check, for example how well spread the neighbours of $x$ are between $C_1$ and $C_2$, and whether either of $C_1$ or $C_2$ is close to being bipartite or not. Throughout this section, we introduce arguments for each of these cases as their own independent lemmas, providing motivation before each one, on how they will relate to the vertex $x$ and subgraphs $C_1$ and $C_2$ as described in this paragraph. We conclude the section with the proof of \cref{thm:exact2/3}, applying each of these intermediary lemmas in turn according to the scenario we are faced with. 

\subsection[A rich subgraph has a little more than k vertices]{A rich subgraph has a little more than $k$ vertices}

Firstly, we focus on graphs with just over $k$ vertices, and provide a short proof that having high minimum degree is sufficient for such a graph to contain every bounded degree forest that is almost spanning. This will be used to show that both $C_1$ and $C_2$ (and more generally the $(2/3-\eps)$-periphery of a rich subgraph) cannot have size in the interval $[k+1,1.15k]$. We use as a black box a well known theorem due to Koml\'os, S\'arközy and Szemer\'edi from 2001.

\begin{theorem}[Koml\'os, S\'arközy and Szemer\'edi {\cite{KOMLÓS_SÁRKÓZY_SZEMERÉDI_2001}}]\label{thm:KSSspanning1/2}
    For all $\eps>0$ there exists $c>0$ such that for all sufficiently large $n\in \N$ the following holds. If $G$ is an $n$-vertex graph with  $\delta(G)\geq (1/2+\eps)n$, then $G$ contains a copy of every $n$-vertex tree $T$ with $\Delta(T)\leq cn/\log n$.
\end{theorem}
\begin{corollary}\label{lem:KSSapplication}
    For all $\Delta\in \N$ and sufficiently large $k\in \N$ the following holds. If $G$ is a graph with $|G|\in [k+1,1.15k]$ and $\delta(G)\geq 0.6k$, then $G$ contains a copy of every forest $F$ with at most $k+1$ vertices and $\Delta(F)\leq \Delta$.
\end{corollary}

\begin{proof}  Add vertices and edges to $F$ in order to obtain a tree $T$ such that $|T| = |G|$ and $\Delta(T)\leq \max\{\Delta,2\}$. We have $\delta(G)\geq 0.6k> (1/2+0.02)\vert G\vert$ so that by \cref{thm:KSSspanning1/2} $G$ contains a copy of every spanning tree with maximum degree at most $\max\{\Delta,2\}$. In particular, $G$ contains a copy of $T$ and hence of $F$.
\end{proof}

\subsection[Neighbours in both sides of a bipartite rich subgraph]{Neighbours in both sides of a bipartite rich subgraph}

Let us now turn our attention to the case where one of the rich subgraphs, say $C_1$, is bipartite and the maximum degree vertex $x$ has several neighbours in both parts of $C_1$. In order to embed a bounded degree tree $T$ into this structure, we find a central vertex $r\in V(T)$ and start embedding the components of $T-r$ back and forth between the bipartition classes of $C_1$. We want to group these components to ensure there is always enough space for the bipartition classes of $T$ to fit in their respective sides of $C_1$. To do this, we first find a way to split up our tree accordingly.
Recall that for a tree $T$ and a vertex $v\in V(T)$, $\textnormal{Even}_T(v)$ is the set of vertices in $T$ with even distance from $v$ (excluding $v$ itself), and $\textnormal{Odd}_T(v)$ is the set of vertices in $T$ with odd distance from $v$.

\begin{lemma}\label{lem:even_odd_2/3}
    Every tree $T$ contains a vertex $r$ such that the components $B_1,\dots,B_m$ of $T - r$ can be partitioned into two classes $\mathcal{C}_1$ and $\mathcal{C}_2$ such that for both $j\in \{1,2\}$, we have
\begin{equation}\label{eq:Odd_even_class_sizes}
        \sum_{B_i \in \mathcal{C}_j} |V(B_i)\cap \mathrm{Even}_T(r)| + \sum_{B_i \notin \mathcal{C}_j} |V(B_i)\cap \mathrm{Odd}_T(r)|\leq \left(\frac{2}{3}- \frac{1}{3\Delta(T)}\right)|T|+\frac{1}{2}.
    \end{equation}
\end{lemma}

\begin{proof}
    For adjacent vertices $u$ and $v$, let $B_{u,v}$ be the component of $T-u$ that contains $v$. We write $a_u$ for $|\mathrm{Even}_T(u)|$, $b_u$ for $|\mathrm{Odd}_T(u)|$, $a_{u,v}$ for $|V(B_{u,v})\cap \mathrm{Even}_T(u)|$ and $b_{u,v}$ for $|V(B_{u,v})\cap \mathrm{Odd}_T(u)|$. Starting with an arbitrary vertex $v_1$, construct a sequence of vertices $(v_i)_{i\in \N}$ such that $v_{i+1}$ maximizes $|a_{v_i,w}-b_{v_i,w}|$ over all neighbours $w$ of $v_i$. Because $T$ is finite and cycle-free, there exists $i\in \N$ such that $v_{i+2}=v_i$, and we fix $u:=v_i$ and $v:=v_{i+1}$.

    We may assume without loss of generality that $|a_{u,v}-b_{u,v}|\leq |a_{v,u}-b_{v,u}|$. Note that
    \begin{align*}
        a_{v,u}&= |V(B_{v,u})\cap \mathrm{Even}_T(v)|=|\mathrm{Odd}_T(u)\setminus V(B_{u,v})|=b_u-b_{u,v},\\
        b_{v,u}&= |V(B_{v,u})\cap \mathrm{Odd}_T(v)|=1+|\mathrm{Even}_T(u)\setminus V(B_{u,v})|=1+a_u-a_{u,v}.
    \end{align*}
    Therefore, 
    \begin{equation*}
        |a_{u,v}-b_{u,v}|\leq |a_{v,u}-b_{v,u}|=|b_u-b_{u,v}-(a_u-a_{u,v})-1|.
    \end{equation*} 
    Letting $\ell=1+\sum_{w\in N_T(u)} |a_{u,w}-b_{u,w}|$, and noting that $a_u = a_{u,v} + \sum_{w\in N(u)\setminus \{v\}}a_{u,w}$ and $b_u = b_{u,v} + \sum_{w\in N(u)\setminus \{v\}}b_{u,w}$ the above yields
    \begin{equation*}
        \ell \geq 1+|a_{u,v}-b_{u,v}| + \Biggl|\sum_{w\in N_T(u)\setminus \{v\}}a_{u,w}-b_{u,w}\Biggr| \geq 1+|a_{u,v}-b_{u,v}| + |a_u-a_{u,v}-(b_u-b_{u,v})|\geq 2 |a_{u,v}-b_{u,v}|.
    \end{equation*}
    However, since $v$ maximizes $|a_{u,w}-b_{u,w}|$ over all neighbours of $u$, this means that $|a_{u,w}-b_{u,w}|\leq \ell/2$ for all $w\in N_T(u)$. If $\ell=1$, then all the quantities $|a_{u,w}-b_{u,w}|$ are zero, and we may take any partition $N_T(u)=\mathcal{D}_1\cup\mathcal{D}_2$. Otherwise $\ell\geq2$, and we may apply \Cref{lem:sum_partition} to $(|a_{u,w}-b_{u,w}|)_{w\in N_T(u)}$ to obtain a partition $N_T(u)=\mathcal{D}_1\cup \mathcal{D}_2$ so that for $j\in \{1,2\}$,
    \begin{equation*}
        \sum_{w\in \mathcal{D}_j} |a_{u,w}-b_{u,w}| \leq \frac{2\ell}{3}.
    \end{equation*}
    Let 
    \begin{equation*}
        \mathcal{C}_j=\{w\in N_T(u):w\in \mathcal{D}_j \text{ and } a_{u,w}\geq b_{u,w}\} \cup \{w\in N_T(u):w\notin \mathcal{D}_j \text{ and } a_{u,w}< b_{u,w}\},
    \end{equation*}
    i.e., we want the bigger partition class of the component $B_{u,w}$ to contribute to \cref{eq:Odd_even_class_sizes} for $j$, when $w\in \mathcal{D}_j$. This gives
    \begin{equation}\label{eq:sum_ab}
        \sum_{w \in \mathcal{C}_j} a_{u,w} + \sum_{w \in N_T(u)\setminus \mathcal{C}_j} b_{u,w} = \sum_{w\in N_T(u)} \min (a_{u,w},b_{u,w}) + \sum_{w\in \mathcal{D}_j} |a_{u,w}-b_{u,w}|\leq \frac{|T|-\ell}{2}+\frac{2\ell}{3},
    \end{equation}
    which is at most $|T|/2+\ell/6$. However, by~\cref{fact:tree_part_sizes}, $\ell$ is bounded by
\begin{equation*}
    \ell \leq 1+ \sum_{w\in N_T(u)} \left(|B_{u,w}|-\frac{2(|B_{u,w}| -1 )}{\Delta}\right)\leq 1+ (|T| -1)- \frac{2(|T|-1)}{\Delta} + \frac{2d(u)}{\Delta}\leq |T| -\frac{2|T|}{\Delta}+3.
\end{equation*}
This implies that \cref{eq:sum_ab} is at most
\begin{equation*}
    \frac{|T|}{2} + \frac{\ell}{6} \leq \frac{|T|}{2} + \frac{|T|}{6} - \frac{|T|}{3\Delta} + \frac{1}{2} = \left(\frac{2}{3} - \frac{1}{3\Delta}\right)|T|+\frac{1}{2}.
\end{equation*}
Noting that the left hand side of \cref{eq:sum_ab} is the same as the left hand side of \cref{eq:Odd_even_class_sizes} taking $r=u$, this proves the lemma.
\end{proof}
As hinted at earlier, we now use \cref{lem:even_odd_2/3} to embed trees into bipartite graphs with minimum degree close to $2k/3$ and with an additional vertex having at least $\Delta$ neighbours in both parts.
\begin{lemma}\label{lem:bipartite2/3-eps}
    For all $\Delta,k\in \N$ with $k>6\Delta$ the following holds. Let $G$ be a graph and $x\in V(G)$ be such that $\delta(G-x)\geq (2/3-1/6\Delta)k$ and $G-x$ is bipartite with parts $Y_1,Y_2$. If $x$ has at least $\Delta$ neighbours both in $Y_1$ and in $Y_2$, then $G$ contains a copy of every $k$-edge tree $T$ with $\Delta(T)\leq \Delta$. 
\end{lemma}

\begin{proof}
    Let $T$ be as in the statement and apply \Cref{lem:even_odd_2/3} to obtain a root $r$ in $T$ and a partition of the components $B_1,\dots,B_m$ of $T-r$ into classes $\mathcal{C}_1$ and $\mathcal{C}_2$, satisfying \cref{eq:Odd_even_class_sizes}. We will greedily construct an embedding $\phi:T\hookrightarrow G$ such that
    \begin{enumerate}[label = \upshape{(\roman*)}]
        \item $\phi(r)=x$ and $\phi (N_T(r))\subseteq N_G(x)$,
        \item $\phi\left(\bigcup_{B_i \in \mathcal{C}_j} V(B_i)\cap \text{Odd}_T(r)\right) \subseteq Y_j$ for both $j\in \{1,2\}$, and
        \item $ \phi\left(\bigcup_{B_i \in \mathcal{C}_j} V(B_i)\cap \text{Even}_T(r)\right)\subseteq Y_{3-j}$ for both $j\in \{1,2\}$.
    \end{enumerate}

We order vertices of $T$ such that their distance to $r$ is ascending, and embed them using this order, starting with $\phi(r) = x$. 
Note that for both $j\in \{1,2\}$, the set $\bigcup_{B_i\in \mathcal{C}_j}V(B_i)\cap N_T(r)$ contains at most $\Delta$ vertices. Thus it is possible to greedily embed all vertices from $\bigcup_{B_i\in \mathcal{C}_j}V(B_i)\cap N_T(r)$ into $N_G(x)\cap V(Y_j)$, maintaining properties (i) and (ii). Suppose now we wish to embed a vertex $u\in V(T)$ at distance $d\geq 2$ from $r$, having already embedded all closer vertices, including the unique neighbour $u'$ of $u$ that has distance $d-1$ from $r$. If $d$ is odd and $u$ is contained in one of the components of $\mathcal{C}_j$, then $\phi(u')\in Y_{3-j}$ by property (iii), and so all of its neighbours are in $Y_j$. Assuming all previously embedded vertices do not contradict properties (ii) or (iii), by \cref{eq:Odd_even_class_sizes}, the number of vertices already embedded into $Y_{j}$ is at most 
\begin{equation*}
    \sum_{B_i \in \mathcal{C}_j} |V(B_i)\cap \text{Odd}_T(r)| + \sum_{B_i 
    \in \mathcal{C}_{3-j}} |V(B_i)\cap \text{Even}_T(r)|-1\leq \left(\frac{2}{3}- \frac{1}{3\Delta}\right)|T|-\frac{1}{2}.
\end{equation*}
Since all neighbours of $\phi(u')$ in $G-x$ belong in $Y_j$, and $\delta(G-x)\geq (2/3-1/6\Delta)k$, we have
\begin{equation*}
    |N_{G-x}(\phi(u'))\cap Y_j|  -  \left(\frac{2}{3}- \frac{1}{3\Delta}\right)|T|+ \frac{1}{2}> \frac{k}{6\Delta} - 1>0.
\end{equation*}
Therefore we can greedily find an unused neighbour of $\phi(u')$ in $Y_{j}$ to embed $u$ into, so that (ii) holds. By a symmetric argument, if $d$ is even, we can greedily embed $u$ to a vertex in $Y_{3-j}$ in order to satisfy (iii).
\end{proof}

\subsection[Many neighbours in one large rich subgraph]{Many neighbours in one large rich subgraph}

Our next lemma will be used to rule out the case where $x$ sends almost all of its neighbours to only one rich subgraph, $C_1$, if $|C_1|$ is large. It is proved without too much difficulty using the machinery of Besomi, Pavez-Sign\'e and Stein given in \cref{subsec:BPS}, and a few of the simple facts about the preservation of cut-density to reduced graphs from \cref{sec:regularity}.
\begin{lemma}\label{lem:one_good_bipartite_reduced}
    For all $\Delta\in \N$ there exists $\eps>0$ such that for all $\rho>0$ 
    and sufficiently large $k\in \N$ \textup{(}in terms of $\Delta$, $\eps$ and $\rho$\textup{)} the following holds. Let $T$ be a $k$-edge tree with $\Delta(T)\leq \Delta$. Let $G$ be a $\rho$-cut-dense graph such that $\delta(G)\geq (2/3-\eps)k$ and $\Delta(G)\geq (1-\eps)k$. If $|G|\in[1.1k,100k]$ then $G$ contains a copy of $T$.
\end{lemma}

\begin{proof}
Let $\eps>0$ be sufficiently small in terms of $\Delta$, let $\rho>0$, and let $\gamma>0$ be sufficiently small in terms of $\eps$, $\Delta$ and $\rho$. Let $M_0$ and $N_0$ be the outputs of \cref{lemma:regularity-degree} when applied with $\gamma_{\ref{lemma:regularity-degree}}=\gamma$ and $(m_0)_{\ref{lemma:regularity-degree}}=\gamma\inv$. Let $k$ be sufficiently large in terms of all other parameters. Suppose that $G$ and $T$ satisfy the conditions of the lemma but that $T$ does not embed in $G$. Let $k_1$ denote the size of the larger bipartition class of $T$.

Let $x\in V(G)$ be such that $d_G(x)= \Delta(G)\geq (1-\eps)k$. By \cref{fact:cutdense_subgraph}, $G-x$ is $\rho/2$-cut-dense. We apply \cref{lemma:regularity-degree} with $\gamma_{\ref{lemma:regularity-degree}}=\gamma$ and $\eta_{\ref{lemma:regularity-degree}}=5\sqrt{\gamma}$ to $G-x$ to obtain a subgraph $H \subseteq G-x$ and a $(\gamma,5\sqrt{\gamma})$-partition of $H$, with corresponding $(\gamma,5\sqrt{\gamma})$-reduced graph $R$, satisfying $|R|\leq M_0$. Note that $R$ satisfies \ref{item:BPS1} and \ref{item:BPS2} from \cref{thm:BPS_combined}.
In particular, $|H|>(1-\gamma)|G|$ and $d_{H}(v)>d_G(v) - (\gamma +5\sqrt{\gamma})|G|\geq (2/3 - 2\eps)k$ for all $v\in V(H)$. Furthermore, $x$ has at least $(1-2\eps)k$ neighbours in $H$. By \cref{fact:cutdense_subgraph}, $H$ is $(\rho/2-2\gamma-10\sqrt{\gamma})$-cut-dense. By \cref{fact:cutdense_reduced_connected} it follows that $R$ is connected.  

Since $|G|\geq 1.1k$, $R$  must be bipartite as otherwise we arrive at a contradiction to \ref{item:BPS1}.
Thus, we may assume $R$ is bipartite with parts $X$ and $Y$. \ref{item:BPS2} implies that both $X$ and $Y$ have at most $(1+\sqrt[4]{\gamma})k_1 \cdot |R|/|H|$ vertices. Define $A= \bigcup V(X)$ and $B=\bigcup V(Y)$. Since $x$ has at least $(1-2\eps)k$ neighbours in $A\cup B$, it follows easily that $x$ has at least $\Delta$ neighbours in both $A$ and in $B$.

Every $U\in V(R)$ is an independent set in $H$ by definition of a regular partition. So, if there exists an edge $uv$ of $H$ contained in either $A$ or $B$, then there exist distinct $U,U'\in X$ or $U,U'\in Y$ respectively, such that $u\in U$ and $v\in U'$. Again by definition of a regular partition, we must have $d(U,U') >5\sqrt{\gamma}$, so that $UU'\in E(R)$, contradicting that $X\cup Y$ is a bipartition of $R$. It follows that $H$ is also bipartite with parts $A$ and $B$.

    Finally, assuming $2\eps<1/6\Delta$, we have $\delta(H)\geq (2/3 - 1/6\Delta)k$, and therefore applying \cref{lem:bipartite2/3-eps} with $G_{\ref{lem:bipartite2/3-eps}}=H\cup \{x\}$, $(Y_1)_{\ref{lem:bipartite2/3-eps}}=A$ and $(Y_2)_{\ref{lem:bipartite2/3-eps}}=B$, we find a copy of $T$ in $G$. This concludes the proof of the lemma. 
\end{proof}

\subsection[Escaping to another rich subgraph]{Many neighbours in two rich subgraphs, and there is a path out of one into another rich subgraph}

Our next scenario surrounds the following case. Suppose the high degree vertex $x$ has $\Delta$ neighbours in both $C_1$ and $C_2$ (playing the roles of $A$ and $B_1$ in the next lemma), that satisfy some additional connectivity conditions to be discussed later. If there exist adjacent vertices $a$ and $b$ in $G-x$ such that $a$ has at least $2\Delta$ neighbours in $C_1$ and $b$ has at least $2\Delta$ neighbours in a different rich subgraph ($B_2$ in the lemma), then every bounded degree $k$-edge tree can be embedded. One can think of the edge $ab$ as a means to escape out of $C_1$ and reach another rich subgraph, which will also have high minimum degree. So, we can split off a branch from $T$ and place it outside $C_1$, preventing $C_1$ from reaching capacity (a bit like how there is a space issue in the extremal example from \cref{fig:extremal2k/3}). When proving the next lemma, we start by finding a central vertex $r\in V(T)$ as in \cref{lem:deletedvertex_compsk/2} and reduce to the case where the two largest branches of $T-r$ have size very close to $k/3$, as this is the situation where we incur most problems, again showing similarities to the extremal example. The embedding process is depicted in \cref{fig:path_embedding}.

\begin{lemma}\label{lemma:embedding_via_path}
For all $\Delta\in \N$ there exists $\eps>0$
such that for all sufficiently large $k\in \N$  the following holds. Let $T$ be a $k$-edge tree with $\Delta(T)\leq \Delta$.
 Let $G$ be a graph containing a vertex $x$ and sets $A,B_1,B_2 \subset V(G)
 \setminus \{x\}$ such that $x$ has at least $3\Delta$ neighbours in both $A$ and $B_1$,  $\delta(G[A])\geq (2/3-\eps)k$, and for both $i\in \{1,2\}$ we have $\delta(G[B_i])\geq (2/3-\eps)k$ and $A\cap B_i = \emptyset$. Additionally suppose that for all $\ell \in \N$ satisfying $\ell\leq 2k/5$, and for all $A'\subseteq A$ satisfying $|A'|\geq |A|- 2\Delta$, if $y,z\in A'$ are distinct vertices, then $G[A']$ contains a $yz$-path of length in the interval $[\ell+1,\ell+24000]$. 
 If $G - x$ contains vertices $a,b$ such that $ab\in E(G)$, $a\in L_{2\Delta}(A)$ and $b\in L_{2\Delta}(B_2)$, then $G$ contains a copy of $T$.
 \end{lemma}

 \begin{proof}
    Let $k\inv \ll \eps\ll \Delta\inv$.
    Suppose to the contrary that the conditions of the lemma are satisfied but that $T$ does not embed in $G$.
    By \cref{lem:deletedvertex_compsk/2} there exists $r\in V(T)$ such that all components $S_1,\dots,S_t$ of $T-r$ have size at most $\ceil{k/2}$, and consider it to be the root of $T$. Note that $\sum_{i\in [t]}|S_i| = k$ and $t\leq \Delta$, and without loss of generality we may assume $|S_1|\geq |S_2|\geq \dots \geq |S_t|$. 
      
\begin{claim}\label{claim:splittable}
    There does not exist $I \subset [t]$ such that $(1/3+\eps)k \leq \sum_{i\in I}|S_i| \leq  (2/3 - \eps)k$.  
\end{claim}

\begin{proofclaim}
    Suppose such an $I$ exists, let $F_1$ be the subforest of $T$ given by $\bigcup_{i\in I}S_i$, and let $F_2$ denote the subforest given by $\bigcup_{i\notin I}S_i = (T -r) - F_1$, noting that $|F_2| = k-|F_1|$. Since $|F_1|\in [(1/3+\eps)k, (2/3 - \eps)k]$, we also have $|F_2| \in [(1/3+\eps)k, (2/3 - \eps)k]$.
    Applying \cref{fact:greedy_embedding} with $G[A \cup \{x\}]$, $F_1\cup \{r\}$ and $(2/3-\eps)k$ playing the roles of $G$, $T$ and $d$ respectively, we find a copy of $F_1\cup \{r\}$ in $G[A \cup \{x\}]$ such that $r$ is embedded at $x$ and the forest $F_1$ is embedded into $G[A]$. Using an identical argument, $F_2\cup \{r\}$ can be embedded into $G[B_1\cup \{x\}] $ such that $r$ is embedded at $x$.  Since $A$ and $B_1$ are vertex-disjoint, this yields a valid embedding of $T$ in $G$, a contradiction.
\end{proofclaim}

\begin{claim}\label{claim:sizesS_1andS_2}
    We have $t\geq 3$ and $(1/3-2\eps)k\leq |S_2|\leq |S_1|\leq (1/3+\eps)k$.
\end{claim}

\begin{proofclaim}
It follows from \cref{claim:splittable} that $|S_i|\leq (1/3+\eps)k$ for all $i\in [t]$, whence it is clear that $t\geq 3$. Let $J\subset [t]$ be such that $\sum_{j\in J}|S_j|$ is maximised amongst all subsets of $[t]$ satisfying $\sum_{j\in J}|S_j|\leq (2/3 - \eps)k$, noting that \cref{claim:splittable} implies $\sum_{j\in J}|S_j|<(1/3+\eps)k$. For every $i\in [t]\setminus J$, we must have $|S_i|>(1/3-2\eps)k$, as otherwise we could replace $J$ by $J\cup \{i\}$ and contradict the maximality property. Since $J \neq [t]$, such an $i$ exists, and in particular the largest component $S_1$ of $T-r$ has size greater than $(1/3-2\eps)k$.

Next, note that if for any $j\in [t] \setminus \{1\}$ we have $3\eps k \leq |S_j| \leq (1/3-2\eps)k$, then $(1/3+\eps)k\leq |S_1|+|S_j|\leq  (2/3-\eps)k$ so taking $I = \{1,j\}$ gives a contradiction to the \cref{claim:splittable}. 
On the other hand, since $t\leq \Delta$ and $|\bigcup_{j=2}^tS_j| = k - |S_1| \geq (2/3-\eps)k$, then by averaging we must have $|S_2|\geq (2/3-\eps)k/(\Delta -1) \geq 3\eps k$. It follows that $|S_2|>(1/3 - 2\eps)k$, proving the claim.
\end{proofclaim}
Note by \cref{claim:sizesS_1andS_2} that $\Delta\geq 3$.
We now construct a path $P_T = p_0p_1\dots p_m$ in $T$ where $p_0 =r$, $p_m$ is a leaf, and $p_1,\dots,p_m$ all belong to the component $S_1$ as follows. Recall that $T(p_i)$ is the subtree of $T$ induced on $p_i$ and its descendants. Choose $p_1$ to be the unique neighbour of $r$ in $S_1$, and for $i\geq 2$ such that $p_{i-1}$ is not a leaf in $T$, let $p_{i}$ be a child of $p_{i-1}$ that maximises $|T(p_i)|$ (choosing arbitrarily if there is a tie). Note that $T(p_{0}) \supset T(p_{1})\supset \dots \supset T(p_{m})$. Since every vertex in $T$ except the root has at most $\Delta -1$ children, then for every $i\in [m-1]$, we have
\begin{equation}\label{eq:size_decrease}
     |T(p_{i+1})|\geq \frac{|T(p_{i})|}{\Delta}.
\end{equation}
Let $\ell$ be the maximum element of $[m]$ such that $|T(p_{\ell})|>k/6$, noting that $\ell\geq 1$ since $|T(p_{1})| 
=|S_1|>k/6$. So, $|T(p_{\ell+1})|\leq k/6$. On the other hand, by \eqref{eq:size_decrease}, we have $|T(p_{\ell+1})|>k/6\Delta$.

Before we begin our embedding process, we reserve some disjoint sets of vertices for future use. First choose distinct $y,y'\in (A\cap N_G(x))\setminus \{a,b\}$. Then, choose $Z_A\subseteq (A\cap N_G(a))\setminus \{b,y,y'\}$ of size $\Delta$, as well as  $Z_B\subseteq (B_2\cap N_G(b))\setminus \{a\}$ of size $\Delta$, and finally $Y_B\subseteq (B_1\cap N_G(x))\setminus (Z_B\cup \{a,b,y,y'\})$ of size $\Delta$. These choices are possible because~\cref{claim:sizesS_1andS_2} gives $\Delta\geq 3$ and by assumptions on $a$, $b$ and $x$. The conditions $A\cap B_1= A\cap B_2 = \emptyset$ ensure the required disjointness of the reserve sets.

Let $a'\in Z_A$ be chosen arbitrarily. Let
\begin{equation*}
    A' = (A\setminus (Z_A \cup \{a,b,y\})) \cup \{a',y'\}.
\end{equation*}
 Since $\Delta\geq 3$, we have $|A'|\geq |A| - (\Delta+3)\geq |A|-2\Delta$. Since $\ell\leq |S_1|<2k/5$, by assumption $G[A']$ contains a $y'a'$-path $P_{A'}$ of length $L$ for some $L\in [\ell+1,\ell+24000]$. We may assume that $k>6\Delta^{24003}$, which together with \eqref{eq:size_decrease} and the fact that $|T(p_{\ell})|>k/6$, implies $P_T$ has length exceeding $\ell + 24003\geq L +3$. By \eqref{eq:size_decrease}, we have
\begin{equation}\label{eq:size_T_p*}
    \frac{k}{6\Delta^{24003}} \leq |T(p_{L+3})|\leq |T(p_{\ell+1})|\leq \frac{k}{6}.
\end{equation}

We now construct an embedding $\phi$ of $T$ in $G$ by embedding different parts of the tree separately, so that each is compatible with all previously embedded vertices. The result is depicted in \cref{fig:path_embedding}. We note that it is possible that $B_1$ and $B_2$ intersect.

\begin{figure}[b]
    \centering
\tikzset{every picture/.style={line width=0.75pt}} 
\begin{tikzpicture}[x=0.75pt,y=0.75pt,yscale=-1,xscale=1]

\draw  [color = orange!80!black]  (170.23,35.02) -- (235.53,51.43) ; 
 \draw  [color = pink!70!black]  (168.09,35.61) -- (219.47,60.28) ;
\draw  [color = purple!60!black]  (168.09,35.61) -- (209.5,72) ;
\draw  [fill={rgb, 255:red, 0; green, 0; blue, 0 }  ,fill opacity=1 ] (165.14,35.02) .. controls (165.14,32.25) and (167.42,30) .. (170.23,30) .. controls (173.04,30) and (175.31,32.25) .. (175.31,35.02) .. controls (175.31,37.79) and (173.04,40.04) .. (170.23,40.04) .. controls (167.42,40.04) and (165.14,37.79) .. (165.14,35.02) -- cycle ;
\draw    (173.96,156.07) -- (140.25,148.85) ;
\draw    (173.96,156.07) -- (148.82,137.78) ;
\draw [line width = 0.6mm]  (189.3,158.07) -- (169.51,154.14) ;
\draw [color=teal]   (225.7,160.18) -- (191.3,158.07) ;
\draw [color=blue!80!black]  (191.3,158.07) -- (225.62,148.76) ;
\draw [color=black]  (191.3,158.07) -- (216.03,184.38) ;
\draw  [color=green!50!black]  (220.5,172.58) -- (191.3,158.07) ;
\draw  [fill={rgb, 255:red, 255; green, 255; blue, 255 }  ,fill opacity=1 ] (208.41,177.1) .. controls (204.65,158.09) and (217.08,139.4) .. (236.16,135.37) .. controls (255.25,131.35) and (273.76,143.5) .. (277.52,162.52) .. controls (281.27,181.53) and (268.85,200.22) .. (249.76,204.25) .. controls (230.68,208.28) and (212.16,196.12) .. (208.41,177.1) -- cycle ;
\draw  [fill={rgb, 255:red, 255; green, 255; blue, 255 }  ,fill opacity=1 ] (206.49,85.19) .. controls (206.49,64.71) and (222.54,48.11) .. (242.35,48.11) .. controls (262.15,48.11) and (278.21,64.71) .. (278.21,85.19) .. controls (278.21,105.67) and (262.15,122.28) .. (242.35,122.28) .. controls (222.54,122.28) and (206.49,105.67) .. (206.49,85.19) -- cycle ;
\draw   (173.96,156.07) -- (138.57,158.71) -- (154.32,125.52) -- cycle ;
\draw  [fill={rgb, 255:red, 255; green, 255; blue, 255 }  ,fill opacity=1 ] (58.09,120.62) .. controls (58.09,91.58) and (80.86,68.03) .. (108.94,68.03) .. controls (137.02,68.03) and (159.79,91.58) .. (159.79,120.62) .. controls (159.79,149.66) and (137.02,173.2) .. (108.94,173.2) .. controls (80.86,173.2) and (58.09,149.66) .. (58.09,120.62) -- cycle ;
\draw  [color={rgb, 255:red, 144; green, 19; blue, 254 }  ,draw opacity=1 ][fill={rgb, 255:red, 144; green, 19; blue, 254 }  ,fill opacity=1 ] (167.81,154.73) .. controls (167.81,152.4) and (169.55,150.51) .. (171.69,150.51) .. controls (173.84,150.51) and (175.58,152.4) .. (175.58,154.73) .. controls (175.58,157.07) and (173.84,158.96) .. (171.69,158.96) .. controls (169.55,158.96) and (167.81,157.07) .. (167.81,154.73) -- cycle ;
\draw  [color={rgb, 255:red, 245; green, 166; blue, 35 }  ,draw opacity=1 ][fill={rgb, 255:red, 245; green, 166; blue, 35 }  ,fill opacity=1 ] (188.46,159.38) .. controls (188.01,157.09) and (189.35,154.88) .. (191.45,154.43) .. controls (193.56,153.99) and (195.63,155.48) .. (196.08,157.77) .. controls (196.53,160.05) and (195.19,162.27) .. (193.09,162.71) .. controls (190.99,163.16) and (188.91,161.66) .. (188.46,159.38) -- cycle ;
\draw [color={rgb, 255:red, 17; green, 108; blue, 219 }  ,draw opacity=1 ][line width=1.5]    (127.41,82.22) .. controls (116.7,83.53) and (120.45,99.3) .. (122.05,106.78) .. controls (123.66,114.26) and (101.71,138.88) .. (120.98,154.38) ;
\draw  [color={rgb, 255:red, 17; green, 108; blue, 219 }  ,draw opacity=1 ][fill={rgb, 255:red, 17; green, 108; blue, 219 }  ,fill opacity=1 ] (117.1,155.29) .. controls (117.1,152.95) and (118.84,151.06) .. (120.98,151.06) .. controls (123.13,151.06) and (124.87,152.95) .. (124.87,155.29) .. controls (124.87,157.62) and (123.13,159.51) .. (120.98,159.51) .. controls (118.84,159.51) and (117.1,157.62) .. (117.1,155.29) -- cycle ;
\draw   (212.92,83) .. controls (208.53,78.04) and (211.46,67.87) .. (219.47,60.28) .. controls (227.48,52.7) and (237.53,50.57) .. (241.92,55.52) .. controls (246.32,60.48) and (243.38,70.65) .. (235.37,78.24) .. controls (227.37,85.83) and (217.31,87.96) .. (212.92,83) -- cycle ;
\draw   (220.87,185.07) .. controls (213.83,182.93) and (210.65,172.28) .. (213.77,161.27) .. controls (216.89,150.26) and (225.13,143.07) .. (232.17,145.2) .. controls (239.21,147.34) and (242.39,157.99) .. (239.27,169) .. controls (236.14,180.01) and (227.91,187.2) .. (220.87,185.07) -- cycle ;
\draw   (116.95,163.89) .. controls (110.47,158.51) and (113.07,144.06) .. (122.76,131.6) .. controls (132.44,119.15) and (145.54,113.41) .. (152.01,118.79) .. controls (158.48,124.17) and (155.88,138.63) .. (146.2,151.08) .. controls (136.52,163.54) and (123.42,169.27) .. (116.95,163.89) -- cycle ;
\draw [color=orange!80!black]  (464.63,35.62) -- (516.01,60.29)   (466.77,35.03) -- (526.5,53) ;
\draw  [color=pink!70!black]  (464.63,35.62) -- (516.01,60.29) ;
\draw  [color=purple!60!black]  (464.63,35.62) -- (504.24,69.74) ;
\draw    (449.75,67.74) -- (466.77,35.03) ;
\draw  [fill={rgb, 255:red, 0; green, 0; blue, 0 }  ,fill opacity=1 ] (461.68,35.03) .. controls (461.68,32.26) and (463.96,30.01) .. (466.77,30.01) .. controls (469.58,30.01) and (471.85,32.26) .. (471.85,35.03) .. controls (471.85,37.8) and (469.58,40.04) .. (466.77,40.04) .. controls (463.96,40.04) and (461.68,37.8) .. (461.68,35.03) -- cycle ;
\draw    (168.35,36.13) -- (127.41,82.22) ;
\draw  [color={rgb, 255:red, 17; green, 108; blue, 219 }  ,draw opacity=1 ][fill={rgb, 255:red, 17; green, 108; blue, 219 }  ,fill opacity=1 ] (123.52,82.22) .. controls (123.52,79.89) and (125.26,78) .. (127.41,78) .. controls (129.55,78) and (131.29,79.89) .. (131.29,82.22) .. controls (131.29,84.56) and (129.55,86.45) .. (127.41,86.45) .. controls (125.26,86.45) and (123.52,84.56) .. (123.52,82.22) -- cycle ;
\draw    (95.68,79.1) -- (165.14,36.13) ;
\draw  [color=green!70!black ,draw opacity=1 ][fill=green!70!black  ,fill opacity=1 ] (91.8,79.1) .. controls (91.8,76.77) and (93.54,74.88) .. (95.68,74.88) .. controls (97.83,74.88) and (99.57,76.77) .. (99.57,79.1) .. controls (99.57,81.44) and (97.83,83.33) .. (95.68,83.33) .. controls (93.54,83.33) and (91.8,81.44) .. (91.8,79.1) -- cycle ;
\draw  [fill={rgb, 255:red, 255; green, 255; blue, 255 }  ,fill opacity=1 ] (355.25,94.07) .. controls (349.85,84.88) and (356.67,70.4) .. (370.49,61.72) .. controls (384.3,53.04) and (399.88,53.45) .. (405.28,62.64) .. controls (410.68,71.83) and (403.86,86.32) .. (390.04,95) .. controls (376.23,103.68) and (360.65,103.26) .. (355.25,94.07) -- cycle ;
\draw [color={rgb, 255:red, 17; green, 108; blue, 219 }  ,draw opacity=1 ][line width=1.5]    (448.43,69.81) -- (464.31,172.2) ;
\draw  [fill={rgb, 255:red, 255; green, 255; blue, 255 }  ,fill opacity=1 ] (514.16,86.74) .. controls (504.04,78.31) and (500.97,64.87) .. (507.3,56.73) .. controls (513.64,48.58) and (526.99,48.81) .. (537.12,57.24) .. controls (547.25,65.67) and (550.32,79.1) .. (543.98,87.25) .. controls (537.64,95.4) and (524.29,95.17) .. (514.16,86.74) -- cycle ;
\draw [color={rgb, 255:red, 128; green, 128; blue, 128 }  ,draw opacity=1 ]   (444.59,69.11) -- (422.29,80.99) ;
\draw [color={rgb, 255:red, 128; green, 128; blue, 128}  ,draw opacity=1 ][fill={rgb, 255:red, 155; green, 155; blue, 155 }  ,fill opacity=1 ]   (446.72,69.38) -- (432.71,82.72) -- (426.04,89.07) ;
\draw [color={rgb, 255:red, 128; green, 128; blue, 128 }  ,draw opacity=1 ]   (451.93,94.79) -- (434.3,107.05) ;
\draw [color={rgb, 255:red, 128; green, 128; blue, 128 }  ,draw opacity=1 ]   (451.8,95.89) -- (438.7,115.41) ;
\draw [color={rgb, 255:red, 128; green, 128; blue, 128 }  ,draw opacity=1 ]   (436.68,105.12) -- (412.26,110.94) ;
\draw [color={rgb, 255:red, 128; green, 128; blue, 128 }  ,draw opacity=1 ]   (436.68,105.12) -- (419.46,123.01) ;
\draw    (507.84,168.02) -- (468.17,172.69) ;
\draw [color={rgb, 255:red, 128; green, 128; blue, 128 }  ,draw opacity=1 ]  (464.31,172.2) -- (433.58,177.08) ;
\draw [color={rgb, 255:red, 128; green, 128; blue, 128 }  ,draw opacity=1 ]   (507.84,168.02) -- (494.6,159.36) ;
\draw [color={rgb, 255:red, 128; green, 128; blue, 128 }  ,draw opacity=1 ]  (464.31,172.2) -- (441.07,188.15) ;
\draw [color={rgb, 255:red, 128; green, 128; blue, 128 }  ,draw opacity=1 ]    (494.6,159.36) -- (475.07,151.31) ;
\draw [color={rgb, 255:red, 128; green, 128; blue, 128 }  ,draw opacity=1 ]  (483.81,141.26) -- (494.6,159.36) ;
\draw  [color={rgb, 255:red, 17; green, 108; blue, 219 }  ,draw opacity=1 ][fill={rgb, 255:red, 17; green, 108; blue, 219 }  ,fill opacity=1 ] (460.45,171.71) .. controls (460.73,169.4) and (462.68,167.74) .. (464.81,168.01) .. controls (466.94,168.28) and (468.44,170.38) .. (468.17,172.69) .. controls (467.89,175.01) and (465.94,176.67) .. (463.81,176.4) .. controls (461.68,176.13) and (460.18,174.03) .. (460.45,171.71) -- cycle ;
\draw  [color={rgb, 255:red, 17; green, 108; blue, 219 }  ,draw opacity=1 ][fill={rgb, 255:red, 17; green, 108; blue, 219 }  ,fill opacity=1 ] (444.59,69.11) .. controls (444.87,66.8) and (446.82,65.14) .. (448.95,65.41) .. controls (451.08,65.68) and (452.58,67.78) .. (452.31,70.09) .. controls (452.03,72.41) and (450.08,74.07) .. (447.95,73.8) .. controls (445.82,73.53) and (444.32,71.43) .. (444.59,69.11) -- cycle ;
\draw  [color=blue!60!black]   (539.91,173.35) -- (595.38,154.48) ;
\draw [color=teal]  (538.52,174.92) -- (587.2,171.17) ;
\draw [color=green!50!black] (538.52,174.92) -- (582.06,185.29) ;
\draw  [fill={rgb, 255:red, 255; green, 255; blue, 255 }  ,fill opacity=1 ] (565.72,172.52) .. controls (565.81,161.2) and (579.09,152.13) .. (595.38,152.26) .. controls (611.67,152.4) and (624.8,161.69) .. (624.72,173.01) .. controls (624.63,184.33) and (611.35,193.4) .. (595.06,193.27) .. controls (578.77,193.14) and (565.63,183.85) .. (565.72,172.52) -- cycle ;
\draw  [line width = 0.6mm] (538.52,174.92) -- (507.84,168.02) ;
\draw  [color={rgb, 255:red, 245; green, 166; blue, 35 }  ,draw opacity=1 ][fill={rgb, 255:red, 245; green, 166; blue, 35 }  ,fill opacity=1 ] (534.61,175.21) .. controls (534.61,172.88) and (536.35,170.99) .. (538.49,170.99) .. controls (540.64,170.99) and (542.38,172.88) .. (542.38,175.21) .. controls (542.38,177.55) and (540.64,179.44) .. (538.49,179.44) .. controls (536.35,179.44) and (534.61,177.55) .. (534.61,175.21) -- cycle ;
\draw  [color={rgb, 255:red, 144; green, 19; blue, 254 }  ,draw opacity=1 ][fill={rgb, 255:red, 144; green, 19; blue, 254 }  ,fill opacity=1 ] (503.96,168.02) .. controls (503.96,165.68) and (505.7,163.79) .. (507.84,163.79) .. controls (509.99,163.79) and (511.73,165.68) .. (511.73,168.02) .. controls (511.73,170.35) and (509.99,172.24) .. (507.84,172.24) .. controls (505.7,172.24) and (503.96,170.35) .. (503.96,168.02) -- cycle ;
\draw    (389.69,66.93) -- (461.68,36.13) ;
\draw  [color=green!70!black  ,draw opacity=1 ][fill=green!70!black,fill opacity=1 ] (389.69,66.93) .. controls (389.69,64.59) and (391.43,62.7) .. (393.58,62.7) .. controls (395.72,62.7) and (397.46,64.59) .. (397.46,66.93) .. controls (397.46,69.26) and (395.72,71.15) .. (393.58,71.15) .. controls (391.43,71.15) and (389.69,69.26) .. (389.69,66.93) -- cycle ;
\draw [draw opacity=0]   (80.84,12.13) -- (262.29,13.23) ;
\draw [draw opacity=0]   (372.02,13.23) -- (553.48,14.34) ;
\draw [draw opacity=0]   (96.9,211.39) -- (278.35,212.5) ;
\draw [draw opacity=0]   (399.86,211.39) -- (581.32,212.5) ;

\draw (40.62,106.83) node [anchor=north west][inner sep=0.75pt]  [font=\small] [align=left] {$\displaystyle A\ $};
\draw (164.73,13) node [anchor=north west][inner sep=0.75pt]   [align=left] {$\displaystyle x$};
\draw (279.23,76.08) node [anchor=north west][inner sep=0.75pt]  [font=\small] [align=left] {$\displaystyle B_{1} \ $};
\draw (279.23,165.42) node [anchor=north west][inner sep=0.75pt]  [font=\small] [align=left] {$\displaystyle B_{2} \ $};
\draw (159.64,160.51) node [anchor=north west][inner sep=0.75pt]  [font=\small,color={rgb, 255:red, 144; green, 19; blue, 254 }  ,opacity=1 ] [align=left] {$\displaystyle a$};
\draw (183.28,163.09) node [anchor=north west][inner sep=0.75pt]  [font=\small,color={rgb, 255:red, 144; green, 19; blue, 254 }  ,opacity=1 ] [align=left] {$\displaystyle \textcolor[rgb]{0.96,0.65,0.14}{b}$};
\draw (94.8,104.95) node [anchor=north west][inner sep=0.75pt]  [font=\small] [align=left] {$\displaystyle \textcolor[rgb]{0.07,0.42,0.86}{P_{A'}}$};
\draw (126.87,143.82) node [anchor=north west][inner sep=0.75pt]  [font=\small,color={rgb, 255:red, 144; green, 19; blue, 254 }  ,opacity=1 ] [align=left] {\textcolor[rgb]{0.07,0.42,0.86}{$\displaystyle a'$}};
\draw (111.05,70.77) node [anchor=north west][inner sep=0.75pt]  [font=\small,color={rgb, 255:red, 144; green, 19; blue, 254 }  ,opacity=1 ] [align=left] {\textcolor[rgb]{0.07,0.42,0.86}{$\displaystyle {\textstyle y'}$}};
\draw (217.08,158.52) node [anchor=north west][inner sep=0.75pt]  [font=\footnotesize] [align=left] {$\displaystyle Z_{B}$};
\draw (218.98,62.01) node [anchor=north west][inner sep=0.75pt]  [font=\footnotesize] [align=left] {$\displaystyle Y_{B}$};
\draw (133.05,124.56) node [anchor=north west][inner sep=0.75pt]  [font=\footnotesize] [align=left] {$\displaystyle Z_{A}$};
\draw (461.27,13) node [anchor=north west][inner sep=0.75pt]   [align=left] {$\displaystyle r$};
\draw (80.42,78.51) node [anchor=north west][inner sep=0.75pt]  [font=\small,color={rgb, 255:red, 126; green, 261; blue, 33 }  ,opacity=1 ] [align=left] {\textcolor{green!70!black}{$\displaystyle {\textstyle y}$}};
\draw (378.03,63.25) node [anchor=north west][inner sep=0.75pt]  [font=\footnotesize,color=green!70!black  ,opacity=1 ] [align=left] {\textcolor{green!70!black}{$\displaystyle u$}};
\draw (360.45,80.14) node [anchor=north west][inner sep=0.75pt]  [font=\footnotesize] [align=left] {$\displaystyle S_{2}$};
\draw (516,63.05) node [anchor=north west][inner sep=0.75pt]  [font=\footnotesize] [align=left] {$\displaystyle F$};
\draw (491.06,174.43) node [anchor=north west][inner sep=0.75pt]  [font=\footnotesize,color={rgb, 255:red, 144; green, 19; blue, 254 }  ,opacity=1 ] [align=left] {$\displaystyle p_{L+2}$};
\draw (524.62,181.32) node [anchor=north west][inner sep=0.75pt]  [font=\footnotesize,color={rgb, 255:red, 144; green, 19; blue, 254 }  ,opacity=1 ] [align=left] {$\displaystyle \textcolor[rgb]{0.96,0.65,0.14}{p_{L+3}}$};
\draw (569.98,166.67) node [anchor=north west][inner sep=0.75pt]  [font=\footnotesize] [align=left] {$\displaystyle T( p_{L+3}) \ \ $};
\draw (455.67,60.07) node [anchor=north west][inner sep=0.75pt]  [font=\footnotesize,color={rgb, 255:red, 144; green, 19; blue, 254 }  ,opacity=1 ] [align=left] {\textcolor[rgb]{0.07,0.42,0.86}{$\displaystyle p_{1}$}};
\draw (456.92,85.88) node [anchor=north west][inner sep=0.75pt]  [font=\footnotesize,color={rgb, 255:red, 144; green, 19; blue, 254 }  ,opacity=1 ] [align=left] {\textcolor[rgb]{0.07,0.42,0.86}{$\displaystyle p_{2}$}};
\draw (458.54,178.88) node [anchor=north west][inner sep=0.75pt]  [font=\footnotesize,color={rgb, 255:red, 144; green, 19; blue, 254 }  ,opacity=1 ] [align=left] {\textcolor[rgb]{0.07,0.42,0.86}{$\displaystyle p_{L+1}$}};
\draw (455.61,141.76) node [anchor=north west][inner sep=0.75pt]  [rotate=-82] [align=left] {$\displaystyle \dotsc $};
\draw (407.36,85.12) node [anchor=north west][inner sep=0.75pt]  [font=\footnotesize,color={rgb, 255:red, 65; green, 117; blue, 5 }  ,opacity=1 ] [align=left] {$\displaystyle \textcolor[rgb]{0.5,0.5,0.5}{R_{1}}$};
\draw (400.94,115.15) node [anchor=north west][inner sep=0.75pt]  [font=\footnotesize] [align=left] {$\displaystyle \textcolor[rgb]{0.5,0.5,0.5}{R_{2}}$};
\draw (475.37,127) node [anchor=north west][inner sep=0.75pt]  [font=\footnotesize,color={rgb, 255:red, 65; green, 117; blue, 5 }  ,opacity=1 ] [align=left] {$\displaystyle \textcolor[rgb]{0.5,0.5,0.5}{R_{L+2}}$};
\draw (408.5,180.81) node [anchor=north west][inner sep=0.75pt]  [font=\footnotesize  ,opacity=1 ] [align=left] {$\displaystyle \textcolor[rgb]{0.5,0.5,0.5}{R_{L+1}}$};
\draw (432.8,125.95) node [anchor=north west][inner sep=0.75pt]  [font=\small] [align=left] {$\displaystyle \textcolor[rgb]{0.07,0.42,0.86}{P_{T}}$};
\end{tikzpicture}
    \caption{\centering The tree $T$ on the right is embedded into the structure in $G$ on the left, so that the corresponding colours and placement align. Note that $B_1$ and $B_2$ need not actually be distinct.}
    \label{fig:path_embedding}
\end{figure}

We start by defining $\phi$ for the subpath $p_0\dots p_{L+1}$ of $P_T$ by $\phi(p_0)=\phi(r) = x$, $\phi(p_1) = y' \in N_G(x)$ and following the path $P_{A'}$ from $y'$ to $a'$ to embed all remaining vertices in order, ending with $\phi(p_{L+1}) =a'$. We continue by defining $\phi(p_{L+2}) = a$ and $\phi(p_{L+3}) = b$, valid since $a'a, ab\in E(G)$. For every $i\in [L+2]$, define $R_i = T(p_i) - T(p_{i+1})$. Note that $R_i$ is a subtree in $S_1$ and think of $R_i$ as being rooted at $p_i$. Observe that $R_i$ and $R_j$ are vertex-disjoint for all distinct $i,j$, and $T[\bigcup_{j\in [i-1]}V(R_j)] = S_1 - T(p_{i})$ for all $i\in [L+3]$.
By \cref{claim:sizesS_1andS_2} and \cref{eq:size_T_p*}, we have
\begin{equation}\label{eq:sum_size_Ri}
 \sum_{j=1}^{L+2}|R_j|= |S_1| - |T(p_{L+3})| \leq \left(\frac{1}{3} +\eps\right)k -\frac{k}{6\Delta^{24003}}
    \leq \left(\frac{1}{3} -3\eps\right)k,
\end{equation}
where we are assuming $\eps<(24\Delta^{24003})^{-1}$.

We will successively embed $R_1,\dots,R_{L+1}$ greedily into $G[A']$ to be compatible with the current image of~$\phi$. Suppose for some $i\in [L+1]$ that we have already extended $\phi$ in this way for all $R_1,\dots ,R_{i-1}$.
Let $A_i' = (A' \setminus (\phi(\bigcup_{j=1}^{i-1}R_j)\cup V(P_{A'})) \cup \{\phi(p_i)\}$.
By \eqref{eq:sum_size_Ri}, we have 
\begin{align*}
\delta(G[A_i']) \geq \delta(G[A']) -  \Biggl| \bigcup_{j=1}^{i-1}\phi(R_j) \cup V(P_{A'}) \Biggr| &\geq
\delta(G[A']) -\sum_{j=1}^{L+2}|R_j| 
\\
&\geq\left(2/3 -\eps\right)k - 2\Delta -\left(1/3 -3\eps\right)k \\ 
& \geq \left(1/3+\eps\right)k >|R_i|,
\end{align*}
since $|R_i| < |S_1|$ and using \cref{claim:sizesS_1andS_2}.
In particular, $\phi(p_i)$ has at least $\Delta$ neighbours in $A_i'$, and so by \cref{fact:greedy_embedding}, there exists a copy of $R_i$ in $A_i'$, with $p_i$ embedded at $\phi(p_i)$. Since this is compatible with the embedding $\phi$ defined for all previously embedded vertices, this extends $\phi$ for $R_i$ using this copy. This holds for all $i\in [L+1]$. 

For $R_{L+2}$, let $A_{L+2} = (A\setminus (\text{im}(\phi)\cup \{y\})) \cup \{a\}$, where $\text{im}(\phi)$ is the set of vertices currently in the image of~$\phi$. Now $Z_A\setminus \{a'\}$ is still unused, so $a = \phi(p_{L+2})$ still has at least $\Delta -1$ available neighbours in~$A_{L+2}$, whilst 
$\delta(G[A_{L+2}\setminus \{a\}])\geq \delta(G[A]) - \sum_{j\in [L+2]}|R_j| -3 > |R_{L+2}|-1$. Again using \cref{fact:greedy_embedding}, we find a copy of $R_{L+2}$ in $A_{L+2}$ with $p_{L+2}$ embedded at $a$, and this is a valid extension for $\phi$.
Altogether we have $\phi(S_1 - T(p_{L+3})) \subseteq G[A \cup \{a\}]$ and $\phi(p_{L+2}) = a$.

Next, we embed the remainder of $S_1$ into $B_2$. Let $B_2' = B_2\setminus (Y_B\cup \{a\})$ and note that $b= \phi(p_{L+3})$ has $\Delta$ available neighbours in $Z_B\subseteq B'_2$. Furthermore by~\eqref{eq:size_T_p*}, $|T(p_{L+3})|\leq k/6$ and so $\delta(G[B_2']) \geq \delta(B_2) - \Delta -1 > |T(p_{L+3})|$. Thus we can embed $T(p_{L+3})$ into $G[B_2'\cup \{b\}]$ using \cref{fact:greedy_embedding} with $p_{L+3}$ embedded at $b$. Again this extends the embedding $\phi$ in a suitable way, so that now all of $S_1$ is embedded into $G$.

Our next aim is to embed $S_2$ into the set of unused vertices in $A$. Let $A^* = A\setminus (\phi(S_1 - T(p_{L+3}))\cup \{b\})$. Observe that $y\in A^*$ is a neighbour of $x$, and
note that $|\phi(S_1- T(p_{L+3}))| = \sum_{i\in [L+2]}|R_i| \leq (1/3 - 3\eps)k$ by \eqref{eq:sum_size_Ri}. We have
\begin{equation*}
    \delta(G[A^*])\geq \delta(G[A]) - (1/3 - 3\eps)k -1 >|S_2|,
\end{equation*}
so~\cref{fact:greedy_embedding} embeds $S_2$ into $A^*$ with its unique vertex $u\in N_T(r)\cap V(S_2)$ mapped to $y$.

It remains to embed the forest $F$ given by $\bigcup_{j=3}^tS_j$ into the set of unused vertices in~$B_1$.  We know that $Y_B$ is contained in~$B_1$ and is unused so far, and so~$x$ has at least $\Delta-2\geq t-2=d_{T[F\cup\{r\}]}(r)$ available neighbours in~$B_1$. Note that $|F| = k - |S_1|-|S_2|\leq (1/3+4\eps)k$ by \cref{claim:sizesS_1andS_2}. We have $\delta(G[B_1 \setminus \text{im}(\phi)])\geq (2/3-\eps)k - k/6-1\geq |F|-1$, where now $\text{im}(\phi)$ is updated as the current set of vertices in the image of~$\phi$. Applying \cref{fact:greedy_embedding} one final time we embed $T[F \cup \{r\}]$ into this subgraph $G[(B_1 \setminus \text{im}(\phi))\cup\{x\}]$ with $r$ rooted at $x$. Extending $\phi$ to account for this completes our copy of $T$, and proves the lemma.
\end{proof}

In order to apply \cref{lemma:embedding_via_path}, we first need to show that, in every almost spanning subgraph of either $C_1$ or $C_2$, any two vertices can be connected by a path of length in a small fixed range. To prove this, we will need the following result, telling us that cut-density is preserved under some random sparsification. A \defn{$p$-random} subset of $V(G)$ is one in which every vertex is selected independently with probability $p$. 
 \begin{lemma}[Pokrovskiy {\cite{AP_hyperstability}}]\label{cut_dense_random_subgraph}
     Let $n\inv \ll p,q$. Let $G$ be an $n$-vertex, $q$-cut-dense graph. If $S$ is a $p$-random subset of $V(G)$ then with high probability $G[S]$ is $p^{20q^{-3}}q^3/400$-cut-dense.
 \end{lemma}
 We now prove \cref{lemma:path_every_length}, which allows us to find the fixed length paths needed for \cref{lemma:embedding_via_path}. For a vertex $w\in V(G)$ and a subset $A\subseteq V(G)$, we write \defn{$d_G(w,A)$} to mean $|N_G(w)\cap A|.$ 

\begin{lemma}\label{lemma:path_every_length}
    For all $\rho>0$ and sufficiently large $k \in \N$ the following holds. Let $G$ be a $\rho$-cut-dense graph on $n\leq 100k$ vertices. If $\delta(G)\geq k/2$ then for every natural number $\ell \leq 2k/5$ and every pair of distinct vertices $y,z\in V(G)$, there exists an $yz$-path of length in the interval $[\ell+1,\ell+24000]$ in $G$.
\end{lemma}

\begin{proof}
Let $\rho>0$ and $k\in \N$ be sufficiently large. Let $G$ be as given with vertices $y,z\in V(G)$, and fix $\ell \leq 2k/5$. Let $G_0 = G-\{y,z\}$. By~\cref{fact:cutdense_subgraph}, $G_0$ is $\rho/2$-cut-dense. Let $\rho' = \frac{(\rho/2)^3}{400(15)^{20(\rho/2)^{-3}}}$. We make the following claim.

\begin{claim}
    There exist disjoint sets $A_1$ and $A_2$ in $V(G_0)$ such that 
\begin{enumerate}[label=\upshape{(\roman*)}]
\item $d_G(w,A_i)\geq k/40$ for every $w\in V(G)$ and $i\in \{1,2\}$;
\item $d_G(w,V(G_0)\setminus (A_1\cup A_2))\geq 2k/5$ for every $w\in V(G)$; and
   \item $G[A_i]$ is $\rho'$-cut-dense for each $i\in \{1,2\}$.
\end{enumerate}
\end{claim}

\begin{proofclaim}
Construct two disjoint random subsets $A_1$ and $A_2$ as follows. Independently for each vertex in $V(G_0)$, place it in $A_1$ with probability $1/15$, in $A_2$ with probability $1/15$, and in neither with probability $13/15$. Fix $i\in \{1,2\}$. By \cref{cut_dense_random_subgraph}, with high probability $G[A_i]$ is $\rho'$-cut-dense. For each $w\in V(G)$, we have 
\begin{equation*}
    \mathbb{E}\left[d_G(w,A_i)\right] = \frac{1}{15}|N_G(w)\cap V(G_0)| \geq \frac{d_G(w)-2}{15} \geq \frac{k}{31}
\end{equation*}
for sufficiently large $k$, and writing $U = V(G_0)\setminus (A_1\cup A_2)$, we have 
\begin{equation*}
    \mathbb{E}\left[d_G(w,U)\right] = \frac{13}{15}|N_G(w)\cap V(G_0)| \geq \frac{13}{15}(d_G(w) - 2) \geq \frac{5k}{12}.
\end{equation*}
Standard Chernoff bounds show that each of the inequalities in (i) and (ii) fails for a specific $w$ with probability $e^{-\Omega(k)}.$ Since $|G|\leq 100k$, a union bound together with (iii) gives a choice of $A_1$ and $A_2$ satisfying all properties.
\end{proofclaim}
As in the proof of the claim, take $U= V(G_0)\setminus (A_1\cup A_2)$, let $G' = G[U]$ and  $H_i = G[A_i]$. It follows from (ii) that $\delta(G')\geq 2k/5\geq \ell$. Therefore we can greedily construct a path $Q_0$ in $G'$ of length $\max\{\ell - 4,0\}$. Let $u$ and $v$ denote its endvertices and note that $Q_0$ does not intersect with any of $A_1$, $A_2$ and $\{y,z\}$ by choice of $U$.

By (i) we have $\delta(H_i)\geq k/40$ for each $i\in \{1,2\}$. Moreover by (i), there exist distinct $y',u'\in A_1$ with $yy',uu'\in E(G)$  and similarly distinct $z',v'\in A_2$ with $zz',vv'\in E(G)$. 
For each $i\in \{1,2\}$, (iii) implies $H_i$ is connected. By \cref{thm:mindeg_diameter}, we have $$\diam(H_i)\leq 
\left\lfloor \frac{3|H_i|}{\delta(H_i)+1}\right \rfloor -1\leq \left\lfloor\frac{3n}{k/40+1}\right \rfloor-1 \leq 11998.$$

Thus there exists a $y'u'$-path $Q_1$ in $H_1$, and a $v'z'$-path $Q_2$ in $H_2$, both of length at most $11998$. We may assume that these paths were chosen to be the shortest of all possible. The sets $\{y,z\}$, $A_1$, $A_2$ and $U$ are pairwise disjoint, so $P = yy'Q_1u'uQ_0vv'Q_2z'z$ gives a $yz$-path in $G$. It is not hard to check that it has length in the interval $[\ell+1,\ell+24000]$, as desired.
\end{proof}

\subsection[Proof of the theorem]{Proof of \cref{thm:exact2/3}}

We now have all of the ingredients needed to prove \cref{thm:exact2/3}.
\begin{proof}[Proof of \cref{thm:exact2/3}]
Let $\Delta\in \N$ and let $\eps>0$ be sufficiently small in terms of $\Delta$. For convenience take $\eta = \eps/2$ and let $\rho$ be the output of \cref{lemma:no-proper-bad} when applied with $\eps_{\ref{lemma:no-proper-bad}} = \eps$ and $\eta_{\ref{lemma:no-proper-bad}}=\eta$. Finally let $k$ be sufficiently large in terms of all other parameters. 

Let $T$ be a $k$-edge tree with $\Delta(T)\leq \Delta$. Suppose that $G$ is a graph satisfying $\delta(G)\geq \floor{2k/3}$ and $\Delta(G)\geq k$ such that $G$ does not contain a copy of $T$. Since $\delta(G)\geq ((2/3-\eps/2) +\eps/4)k$, we may apply \cref{lemma:no-proper-bad} to $G$ with 
$\eps_{\ref{lemma:no-proper-bad}} = \eps/4$,
$\eta_{\ref{lemma:no-proper-bad}} = \eta$ and
$a_{\ref{lemma:no-proper-bad}} = 2/3-\eps/2$
to obtain a collection of vertex-disjoint $(2/3-3\eps/4, \rho,k)$-rich subgraphs $C_1,\dots,C_m$ such that every $C_i$ is $\eta k$-closed, and for every $v\in V(G)$ there exist $i,j\in [m]$ for which $|N_G(v)\setminus (C_i\cup C_j)|< \eta k$. Let $x\in V(G)$ satisfy $d_G(x) = \Delta(G)\geq k$, and without loss of generality assume $|N_G(x)\setminus (C_1\cup C_2)|<\eta k$. Note if the indices $i,j$ for $x$ are equal, then this weaker assumption still holds.

  \textbf{Case 1:} $x$ has at least $\eta k$ neighbours in both $C_1$ and $C_2$.

We make the following claim.
  \begin{claim}\label{claim:path-between-components}
      Let $i\in \{1,2\}$ and $j\in [m]\setminus \{i\}$. If there exists $a,b\in V(G)$ such that $ab\in E(G)$, $a\in L_{\eta k }(C_i)$ and $b\in L_{\eta k}(C_j)$ then $x\in \{a,b\}$.      
  \end{claim}
  \begin{proofclaim}
  Let $i,j, a,b$ be as in the claim and suppose $x\notin \{a,b\}$. Note that every $A'\subset V(C_i)$ with $|A'|\geq |C_i|- 2\Delta$ satisfies $\delta(G[A'])\geq \delta(C_i) - 2\Delta \geq k/2$. Applying \cref{lemma:path_every_length} with $G_{\ref{lemma:path_every_length}} = G[A']$, we deduce that for every pair of vertices $y,z \in A'$ and every $\ell \leq 2k/5$, there is a $yz$-path of length in $[\ell+1,\ell+24000]$ in $G[A']$. Then the original graph $G$ satisfies the assumptions of \cref{lemma:embedding_via_path} with $C_i$, $C_{3-i}$ and $C_j$ playing the roles of $A$, $B_1$ and $B_2$ respectively. So, $G$ contains a copy of $T$, a contradiction.  
  \end{proofclaim}

For $i\in \{1,2\}$ all neighbours of every $a\in L_{\eta k}(C_i)\setminus \{x\}$ lie within $L_{\eta k}(C_i)$. To see this, suppose that $a\in L_{\eta k}(C_i)\setminus \{x\}$ has a neighbour $b\notin L_{\eta k}(C_i)$. Clearly $b\neq x$ by assumption of Case 1. There exist $j_1,j_2\in [m]$ such that all but $\eta k$ neighbours of $b$ are in $V(C_{j_1})\cup V(C_{j_2})$, and in particular, there exists $j\in [m]$ such that $b$ has at least $\eta k$ neighbours in $V(C_j)$. Since $b\in L_{\eta k}(C_j)\setminus L_{\eta k}(C_i)$, we have $i\neq j$ and therefore we arrive at a contradiction with \Cref{claim:path-between-components}. Thus, we have established that the minimum degree of both $G[L_{\eta k}(C_1)\setminus \{x\}] $ and $G[L_{\eta k}(C_2) \setminus \{x\}]$ is at least $\lfloor 2k/3 \rfloor -1$. Suppose now that there exists a vertex $a\in L_{\eta k}(C_1)\cap L_{\eta k}(C_2)\setminus \{x\}$. Then $a \in L_{\eta k}(C_1)$ and $a$ has a neighbour $b\in V(C_2)\setminus \{x\}$. In particular, $b\in L_{\eta k }(C_2)$, contradicting \cref{claim:path-between-components}. 
It follows that $G[L_{\eta k}(C_1)\setminus\{x\}]$ and $G[L_{\eta k}(C_2)\setminus\{x\}]$ are vertex-disjoint graphs.

By \Cref{cor:partition_subforests_T-v}\ref{item:Fsize1} there exists $r\in V(T)$ and a partition of the components of $T-r$ into two vertex-disjoint subforests $F_1$ and $F_2$, both with at most $\lfloor 2k/3 \rfloor$ vertices. For both $i\in \{1,2\}$, applying \cref{fact:greedy_embedding}, we can greedily embed the tree $F_i\cup \{r\}$ rooted at $r$ into $G[L_{\eta k}(C_i)]$, such that $r$ is
mapped to $x$, together forming a complete copy of $T$ in $G$.

  \textbf{Case 2:} $x$ has less than $\eta k$ neighbours in one of $C_1$ or $C_2$.

  Without loss of generality assume that $|N_G(x)\cap C_2|< \eta k$. 
  Moreover $\Delta(C_1\cup \{x\})\geq |N_G(x)\cap C_1|\geq (1-\eps)k$.
    Note that $\delta(C_1\cup \{x\})\geq (2/3-3\eps/4)k$. Let $\rho' = \min\{\rho/2,(1-\eps)/100\}$. It is not hard to check that $G[C_1 \cup \{x\}]$ is $\rho'$-cut-dense. 
  If $|C_1|\geq 1.1k$, then we may apply \cref{lem:one_good_bipartite_reduced} with $G[C_1\cup \{x\}]$ playing the role of $G$ to find a copy of $T$. 
  
 This leaves the possibility that $|C_1|< 1.1k$. Let $k' = (1-\eps)k -1$. Let $S=L_{(2/3-4\eps)k}(C_1)$, and note that since $C_1$ is $\eta k$-closed, we have $\vert S\setminus C_1\vert < \eta k$, and thus $\vert S\vert < 1.12k < 1.15k'$.
 If $\vert S\vert > k$, we know that since $\delta(G[S])>(2/3-4\eps)k$, $G[S]$ contains $T$ by \Cref{lem:KSSapplication}.
 
 We may therefore assume that $\vert S\vert \leq k$, which, as $x\in S$ and $d_G(x)\geq k$, implies that $x$ has a neighbour $y\notin S$. The number of neighbours of $y$ in $S$ is less than $(2/3-4\eps)k+ \vert S\setminus C_1\vert < (2/3-3\eps)k$. However, all but at most $\eta k$ neighbours of~$y$ are in the union of two of the rich subgraphs, whence there has to exist $j\in [m]\setminus \{1\}$ such that $y\in L_{\eta k}(C_j)$. Furthermore, since all but $2\eta k =\eps k$ neighbours of~$x$ are in~$C_1$, the induced subgraph $G[S\cap N_G(x)]$ has size at least $(1-\eps)k = k'+1$ and minimum degree at least $(2/3-5\eps)k\geq 0.6k'$.

 Fix an edge $uv\in E(T)$ with the property that all components of $T-u$ have at most $\lceil k/2 \rceil$ vertices and such that the component $T_v$ of $T-u$ which contains~$v$ has size at least $k/\Delta$, possible by \cref{lem:deletedvertex_compsk/2} and since $d_T(u) \leq \Delta$. It follows that the number of vertices in the forest $F=T-(\{u\}\cup V(T_v))$ is between $\lfloor k/2 \rfloor$ and $(1-1/\Delta)k<(1-\eps)k = k'+1$. By \Cref{lem:KSSapplication} with $k_{\ref{lem:KSSapplication}} = k'$, there exists an embedding of~$F$ into $G[S\cap N_G(x)]$, which we can extend to an embedding of~$T$ into~$G$ by mapping~$u$ to~$x$, $v$ to~$y$, and $T_v$ greedily into $C_j$.
\end{proof}

\section[Second neighbourhood]{Second neighbourhood: Proof of \cref{thm:secondnbhd_result}}\label{sec:second_nbhd}
As briefly sketched in \cref{subsec:outline}, our approach for proving \cref{thm:secondnbhd_result} again begins with an application of \cref{lemma:no-proper-bad}. We start with a graph $G$ satisfying $\delta(G)\geq (1+\eps)k/2$ and containing a vertex $x\in V(G)$ whose first and second neighbourhoods both have size at least $(1+\eps)4k/3$. We find a collection of vertex-disjoint rich subgraphs $(C_i)_{i\in [m]}$ as in \cref{lemma:no-proper-bad}, and without loss of generality, $x$ has almost all of its (first) neighbours in $V(C_1)\cup V(C_2)$. 
We consider property \ref{item:main3} from \cref{lemma:no-proper-bad} for each of the second neighbours $z\in N_G^2(x)$, and think of the two rich subgraphs $C_i$ and $C_j$ in which $z$ has almost all of its neighbours as being `associated' with $z$. If there are only few members of $(C_i)_{i\in [m]}$ that are associated with some second neighbour of $x$, then we can consider the subgraph $\hat{G}$ of $G$ induced on all of the associated rich subgraphs together with $\{x\}\cup N_G(x)\cup V(C_1) \cup V(C_2)$. Then $\hat{G}$ will have bounded size, minimum degree above $k/2$, and will contain enough first and second neighbours of $x$ to conclude via \cref{thm:secondneighbourhood_dense}.  

Now suppose on the other hand that there are many elements of $(C_i)_{i\in [m]}$ that are associated with some second neighbour of $x$, then we find a large matching in $G$, where for every edge, one endpoint is in $V(C_1)\cup V(C_2)$, and the other endpoint is a second neighbour of $x$ lying outside of these subgraphs, that has a unique subgraph in $(C_i)_{i\in [m]}$ associated with it. We use these matching edges to `escape' from $C_1\cup C_2$ to many different rich subgraphs, that are all vertex-disjoint. 

In order to embed a tree $T$ in this setup, we will split up $T$ into a subtree $S$, a matching $M$, and some additional components of controlled size. We will then forcibly embed the matching $M$ onto the matching edges in $G$ mentioned above, the subtree $S$ into $C_1\cup C_2$, and the remaining components will be suitably small to fit in the associated rich subgraphs. 

\begin{lemma}\label{lem:treesplit_MSF}
 For all $\Delta\in \N$ and sufficiently large $k \in \N$ the following holds. Let $T$ be a $k$-edge tree with $\Delta(T)\leq \Delta$. 
There exists a matching $M$, a tree $S$, and a forest $F$ in $T$ whose edge sets form a partition of $E(T)$, and satisfy the following properties:
    \begin{enumerate}[label=\upshape{(P\arabic*)}, leftmargin =\widthof{P10000}]
        \item $S$ and $F$ are vertex-disjoint, \label{item:MSF1}
        \item every edge of $M$ has exactly one endpoint in $S$ and one endpoint in $F$, \label{item:MSF2}
        \item the tree $S$ and every component of $F$ each have at most $\ceil{k/2}$ vertices, \label{item:MSF3}
        \item all vertices of $V(M)\cap V(S)$ belong to the same bipartition class in $T$,  \label{item:MSF5}
        \item the smallest subtree of $S$ containing $V(M)\cap V(S)$ contains at most $\Delta^{4\Delta+1}$ vertices. \label{item:MSF6}
    \end{enumerate}
\end{lemma}

\begin{proof}
        Let $r \in V(T)$ be chosen according to \cref{lem:deletedvertex_compsk/2} so that all components of $T - r$ have size at most $\ceil{k/2}$. For $j\in \N$, let \defn{$D_j(r)$} denote the set of vertices at distance exactly $j$ from $r$. 
        
        Recall that for a vertex $w$, $T(w)$ is the subtree of descendants of $w$ including $w$ itself. For every $j\in \N$ and every $v\in D_j(r)$ that is not a leaf, fix a child $\firstborn{v}$ of $v$ such that $|T(\firstborn{v})|$ is maximal among the children of $v$. We call $\firstborn{v}$ the \defn{escape vertex} of $v$.  
        
        Let $B=\bigcup_{j=1}^{2\Delta} D_{2j}(r)$, let $B_0$ be the set of non-leaves in $B$, and let $X$ be the set of vertices in $B_0$ that are \emph{not} a descendant of the escape vertex of a vertex in $B_0$, i.e., 
        \begin{equation*}
            X\coloneqq B_0 \setminus \bigcup_{v\in B_0} V(T(\firstborn{v})).
        \end{equation*}
        We define further
        \begin{align*}
            F&\coloneqq \bigcup_{v\in X} T(\firstborn{v}),&
            M&\coloneqq T\left[\bigcup_{v\in X} \{v,\firstborn{v}\} \right],&
            S&\coloneqq T\left[ V(T)\setminus \bigcup_{v\in X} V(T(\firstborn{v}))\right].
        \end{align*}

        See \cref{fig:MSF_decomp} for an example of the decomposition.
        It is clear that $S$ and $F$ are disjoint, so that \ref{item:MSF1} is satisfied. Furthermore, it follows from the definition of $X$ that $X\subset V(S)$. Thus, every edge in $M$ consists of a vertex $v\in V(S)$ and its unique escape vertex $\firstborn{v}\in V(F)$, and $M$ is a matching satisfying \ref{item:MSF2}. Since $V(M)\cap V(S)=X\subset B$, and the vertices in $B$ are those whose distance to $r$ is even and at most $4\Delta$, and $\Delta(T)\leq \Delta$, we know that \ref{item:MSF5} and \ref{item:MSF6} hold.

        As for \ref{item:MSF3}, it is clear that every component of $F$ is contained in a component of $T-r$ and therefore has size at most $\ceil{k/2}$ by choice of $r$. 
        It remains to show that $|S|\leq \ceil{k/2}$. Let $R$ be the subtree of $T$ induced by the vertices of distance at most $4\Delta$ from $r$, which we know has size at most $\Delta^{4\Delta +1}$. Thus, if we suppose towards a contradiction that $|S|> \ceil{k/2}$, then $|S\setminus R|>3k/8$ as $k$ is sufficiently large. However, this implies that
        \begin{equation*}
            3k/8<|S\setminus R|\leq \sum_{v\in X\cap D_{4\Delta}(r)} (|T(v)|-1)< \frac{1}{2\Delta}\sum_{j=1}^{2\Delta}\sum_{v\in X\cap D_{2j}(r)} (|T(v)|-1)=\frac{1}{2\Delta}\sum_{v\in X} (|T(v)|-1).
        \end{equation*}
        On the other hand, by definition of the escape vertices, we obtain
        \begin{equation*}
            |F|=\sum_{v\in X} |T(\firstborn{v})|\geq \sum_{v\in X} \frac{|T(v)|-1}{\Delta},
        \end{equation*}
        which, by the above, is more than $3k/4$. 
        However, since $S$ and $F$ are vertex-disjoint, they cannot both be larger than $\ceil{k/2}$, and the contradiction completes the proof.
\end{proof}

\begin{figure}
    \centering
    \begin{tikzpicture}[
    u/.style={fill=black,minimum size =4pt,ellipse,inner sep=0pt},
    node distance=1.5cm,scale=0.5
    ]

 \filldraw[blue!30!white] (-10,6) circle (8pt); 
\filldraw[blue!30!white] (-4.5,6) circle (8pt); 
\filldraw[blue!30!white] (0,6) circle (8pt); 
\filldraw[blue!30!white] (4.5,6) circle (8pt); 
\filldraw[blue!30!white] (9.5,6) circle (8pt); 
 \filldraw[blue!30!white] (-9.5,4) circle (8pt); 
\filldraw[blue!30!white] (-6.5,4) circle (8pt); 
\filldraw[blue!30!white] (-1,4) circle (8pt); 
\filldraw[blue!30!white] (5.5,4) circle (8pt); 
\filldraw[blue!30!white] (7.5,4) circle (8pt); 
\filldraw[blue!30!white] (9.5,4) circle (8pt); 

\filldraw[orange!20!white,rotate around={-90:(-12.5,4)}] (-12.5, 4) ellipse (1.5cm and 1.9cm); 
\filldraw[orange!20!white,rotate around={90:(12.5,3.8)}] (12.5, 3.8) ellipse (1.5cm and 1.5cm); 
\filldraw[orange!20!white,rotate around={-90:(-4,4)}] (-4, 4) ellipse (1.3cm and 1.4cm); 
\filldraw[orange!20!white,rotate around={-90:(2,4)}] (2, 4) ellipse (1.3cm and 1.4cm); 

\node[u] (r) at (0,8) {}; 
\node[u] (a1) at (-6,7) {}; 
\node[u] (a2) at (0,7) {}; 
\node[u] (a3) at (6,7) {}; 
\node[u] (b1) at (-10,6) {}; 
\node[u] (b2) at (-4.5,6) {}; 
\node[u] (b3) at (0,6) {}; 
\node[u] (b4) at (4.5,6) {}; 
\node[u] (b5) at (9.5,6) {}; 
\node[u] (c1) at (-12.5,5) {}; 
\node[u] (c2) at (-8,5) {}; 
\node[u] (c3) at (-4,5) {}; 
\node[u] (c4) at (-1,5) {}; 
\node[u] (c5) at (2,5) {}; 
\node[u] (c6) at (4,5) {}; 
\node[u] (c7) at (6.5,5) {}; 
\node[u] (c8) at (9.5,5) {}; 
\node[u] (c9) at (12.5,5) {}; 
\node[u] (d1) at (-14,4) {}; 
\node[u] (d2) at (-12.5,4) {}; 
\node[u] (d3) at (-11,4) {}; 
\node[u] (d4) at (-9.5,4) {}; 
\node[u] (d5) at (-6.5,4) {}; 
\node[u] (d6) at (-5,4) {}; 
\node[u] (d7) at (-3,4) {}; 
\node[u] (d8) at (-1,4) {}; 
\node[u] (d9) at (1,4) {}; 
\node[u] (d10) at (3,4) {}; 
\node[u] (d11) at (5.5,4) {}; 
\node[u] (d12) at (7.5,4) {}; 
\node[u] (d13) at (9.5,4) {}; 
\node[u] (d14) at (11.5,4) {}; 
\node[u] (d15) at (13.5,4) {}; 
\node[u] (e1) at (-10.5,3) {}; 
\node[u] (e2) at (-8.5,3) {}; 
\node[u] (e3) at (-7.5,3) {}; 
\node[u] (e4) at (-5.5,3) {}; 
\node[u] (e5) at (-1,3) {}; 
\node[u] (e6) at (4.5,3) {}; 
\node[u] (e7) at (6.5,3) {}; 
\node[u] (e8) at (8.5,3) {}; 
\node[u] (e9) at (10.5,3) {}; 
\node[rotate=90] at (-12.5,3.2) {\textbf{\ldots}};
\node[rotate=90] at (-4,3.4) {\textbf{\ldots}};
\node[rotate=90] at (2,3.4) {\textbf{\ldots}};
\node[rotate=-90] at (12.5,3.4) {\textbf{\ldots}};

\draw[line width=1pt] (r) -- (a1);
\draw[line width=1pt] (r) -- (a2);
\draw[line width=1pt] (r) -- (a3);
\draw[line width=1pt] (a1) -- (b1);
\draw[line width=1pt] (a1) -- (b2);
\draw[line width=1pt] (a2) -- (b3);
\draw[line width=1pt] (a3) -- (b4);
\draw[line width=1pt] (a3) -- (b5);
\draw[line width=2pt,red] (b1) -- (c1);
\draw[line width=1pt] (b1) -- (c2);
\draw[line width=2pt,red] (b2) -- (c3);
\draw[line width=1pt] (b3) -- (c4);
\draw[line width=2pt,red] (b3) -- (c5);
\draw[line width=2pt,red] (b4) -- (c6);
\draw[line width=1pt] (b5) -- (c7);
\draw[line width=1pt] (b5) -- (c8);
\draw[line width=2pt,red] (b5) -- (c9);
\draw[line width=1pt] (c1) -- (d1);
\draw[line width=1pt] (c1) -- (d2);
\draw[line width=1pt] (c1) -- (d3);
\draw[line width=1pt] (c2) -- (d4);
\draw[line width=1pt] (c2) -- (d5);
\draw[line width=1pt] (c3) -- (d6);
\draw[line width=1pt] (c3) -- (d7);
\draw[line width=1pt] (c4) -- (d8);
\draw[line width=1pt] (c5) -- (d9);
\draw[line width=1pt] (c5) -- (d10);
\draw[line width=1pt] (c7) -- (d11);
\draw[line width=1pt] (c7) -- (d12);
\draw[line width=1pt] (c8) -- (d13);
\draw[line width=1pt] (c9) -- (d14);
\draw[line width=1pt] (c9) -- (d15);
\draw[line width=2pt, red] (d4) -- (e1);
\draw[line width=1pt] (d4) -- (e2);
\draw[line width=2pt, red] (d5) -- (e3);
\draw[line width=1pt] (d5) -- (e4);
\draw[line width=2pt, red] (d8) -- (e5);
\draw[line width=1pt] (d11) -- (e6);
\draw[line width=2pt, red] (d11) -- (e7);
\draw[line width=2pt, red] (d13) -- (e8);
\draw[line width=1pt] (d13) -- (e9);

\node[above=1mm] at (r) {$r$};  
    \end{tikzpicture}
    \caption{The top few levels of the tree decomposition detailed in the proof of \cref{lem:treesplit_MSF}. 
    Vertices shaded in blue belong to $X$. Edges coloured red belong to the matching $M$, these all lie between elements of $X$ and their escape vertices. Each component of $F$ is shaded in orange, these are the branches of $T$ that ``descend'' from an edge in $M$. The remaining edges of $T$ that do not belong in $F$ or $M$ define the subtree $S$.}
    \label{fig:MSF_decomp}
\end{figure}

Our next lemma will provide a mechanism for finding a copy of the subtree $S$ in $C_1\cup C_2$ (as mentioned at the start of this section) so that the vertices in $V(S)\cap V(M)$ are embedded into a predetermined set.

\begin{lemma}\label{lem:embedfor_matchingpoints}
For all $\Delta\in \N$ there exists $D\in\N$ such that for all $\eps>0$ and all sufficiently large $k\in\N$ \textup{(}in terms of $\Delta$, $\eps$ and $D$\textup{)} the following holds.
    Let $G$ be a graph with $\delta(G)\geq (1+\eps)k/2$ and $|G|<250k$. Let $P \subseteq V(G)$ be such that $|P|\geq D$. Let $T$ be a tree on at most $\ceil{k/2}$ vertices satisfying $\Delta(T)\leq \Delta$. Let $Q\subseteq V(T)$ be non-empty and such that
    all vertices in $Q$ are in the same bipartition class in $T$, and the smallest subtree of $T$ containing $Q$ has size at most $\Delta^{4\Delta+1}$. 
    Then there exists an embedding $\psi:T\hookrightarrow G$ such that $ \psi(v)\in P$ for all $v\in Q$.
\end{lemma}

\begin{proof}
    Let $\Delta\in \N$ and let $D\in \N$ be sufficiently large in terms of $\Delta$. Let $\eps>0$ and let $k\in \N$ be sufficiently large in terms of all other parameters. Choose an arbitrary subset $P_0\subseteq P$ of size exactly $D$ and replace $P$
by $P_0$. It is enough to find an embedding that maps $Q$ into this smaller set. Let $R$ be a random subset of $V(G)\setminus P$ of size $|P|$. 

For every $p\in P$, we have $d_G(p,G\setminus P)\geq d_G(p)-|P|\geq k/2$ for sufficiently large $k$, and so $\E[d_G(p,R)]\geq \frac{|R|}{|G|}d_G(p,G\setminus P) \geq \frac{|R|}{250k} \cdot \frac{k}{2 }\geq |R|/500$.
    By linearity, $\E[e_G(P,R)]\geq |P||R|/500 = |P|^2/500$. Thus there exists a set $R \subseteq V(G)\setminus P$ with $|R|=|P|$ such that the bipartite graph $H$ defined on $P\cup R$ with all edges between $P$ and $R$ satisfies $e(H)\geq |P|^2/500 \geq D^2/500$. By the commonly known result of Köv\'ari, S\'os and Tur\'an, $H$ contains a copy $K$ of the complete bipartite graph $K_{t,t}$ with $t = \Delta^{4\Delta+1}$, provided $D$ is sufficiently large in terms of $\Delta$.

    Let $T_Q\subseteq T$ be the smallest subtree of $T$ containing $Q$ so that $|T_Q|\leq \Delta^{4\Delta+1}$. There exists an embedding $\psi^*:T_Q\hookrightarrow K$ such that the bipartition class of $T_Q$ containing all vertices in $Q$ is embedded into $P$. We can extend $\psi^*$ to an embedding $\psi$ of $T$ in $G$ by mapping $V(T) \setminus V(T_Q)$ greedily into $V(G)\setminus V(K)$. This is possible since $|T|-|T_Q|\leq \ceil{k/2}$ and every $v\in V(G)$ has at least $(1+\eps)k/2 - 2t \geq \ceil{k/2}$ neighbours in $V(G)\setminus V(K)$. 
\end{proof}

Before proving \cref{thm:secondnbhd_result} we make one final simple observation about finding a large matching in a bipartite graph.

\begin{fact}\label{fact:bipartite_matching_size}
    For a bipartite graph $B$ with parts $X$ and $Y$, such that every vertex in $Y$ has degree at least 1, and every vertex in $X$ has degree at most $d$, there exists a matching of size $|Y|/d$.
\end{fact}

\begin{proof}
    Let $M$ be a largest matching in $B$, let $X' = X\cap V(M)$ and $Y' = Y\cap V(M)$. There is no edge between $X\setminus X'$ and $Y\setminus Y'$, so each $y\in Y$ has at least one neighbour in $X'$. Therefore, there are at least $|Y|$ edges from $Y$ to $X'$. Since each $x\in X'$ has at most $d$ neighbours in $Y$, then we have $|M| = |X'|\geq \frac{|Y|}{d}$, as desired.
\end{proof}

We are now ready to combine everything to prove \cref{thm:secondnbhd_result}.

\begin{proof}[Proof of \cref{thm:secondnbhd_result}]
    Let $\Delta\in \N$ and $\eps>0$. Let $D$ be the maximum of $\Delta^{4\Delta+1}$, and the output $D_{\ref{lem:embedfor_matchingpoints}}$ of \cref{lem:embedfor_matchingpoints} when applied with $\Delta$, noting that this is independent of $\eps$. We may therefore assume without loss of generality that $\eps<(100D\Delta)\inv$ as decreasing $\eps$ only strengthens \cref{thm:secondnbhd_result}. Let $\eta = \eps/8$ and let $\rho$ be the output of \cref{lemma:no-proper-bad} when applied with $\eps_{\ref{lemma:no-proper-bad}} = \eps/2$ and $\eta_{\ref{lemma:no-proper-bad}}=\eta$. Let $\gamma$ be sufficiently small in terms of $\eps$ and $\rho$.
    Let $M_0$ and $N_0$ be the outputs of \cref{lemma:regularity-degree} when applied with $\gamma_{\ref{lemma:regularity-degree}}=\gamma$ and $(m_0)_{\ref{lemma:regularity-degree}}=\gamma\inv$. Finally, let $k$ be sufficiently large in terms of all other parameters. 

    Let $T$ be a $k$-edge tree with $\Delta(T)\leq \Delta$. Suppose that $G$ and $x$ satisfy the assumptions of \cref{thm:secondnbhd_result}, and that $T$ does not embed in $G$. Apply \Cref{lemma:no-proper-bad} to $G$ with $\eps_{\ref{lemma:no-proper-bad}}=\eps/2$, $\eta_{\ref{lemma:no-proper-bad}}=\eta$ and $a_{\ref{lemma:no-proper-bad}} = 1/2$ to obtain vertex-disjoint $(1/2+\eps/8,\rho,k)$-rich subgraphs $C_1,\dots, C_m$ that are each $\eta k$-closed, such that 
    \begin{equation}\label{item:SN2}
        \text{for all } y\in V(G) \text{ there exist } i,j\in [m] \text { such that } |N_G(y)\setminus (C_i\cup C_j)|< \eta k,
    \end{equation}
    and, without loss of generality we have $\vert N_G(x) \setminus (C_1 \cup C_2)\vert < \eta k$.

Our first goal is to show that the number of neighbours of $x$ cannot be too unbalanced between~$C_1$ and~$C_2$.
\begin{claim}\label{claim3:balanced_neighbours}
   In $G$, $x$ has at least $\eps k/6$ neighbours in both $C_1$ and $C_2$.
\end{claim}

\begin{proofclaim}
    Without loss of generality suppose to the contrary that $|N_G(x)\cap C_2|< \eps k/6$. Then $|N_G(x)\cap C_1|>d_G(x)- \eta k -\eps k/6\geq (1+\eps/8)4k/3$.

    Since $\gamma$ was chosen sufficiently small in terms of $\eps$, fix $\theta>0$ such that $\gamma\leq \theta^4/10^{10}$, and $\gamma \ll \theta\ll \eps$.
    Apply \cref{lemma:regularity-degree} to $C_1-x$ with $\gamma_{\ref{lemma:regularity-degree}}=\gamma$ and $\eta_{\ref{lemma:regularity-degree}}=5\sqrt{\gamma}$ to obtain a subgraph $H_1 \subseteq C_1-x$, and an $(\gamma,5\sqrt{\gamma})$-regular partition of $V(H_1)$, with corresponding reduced graph $R$. 
    
    Since $\delta(C_1)\geq (1/2+\eps/8)k$, then $d_{H_1}(v)\geq d_{C_1-x}(v)- (\gamma+5\sqrt{\gamma}) |C_1-x| \geq (1+\theta)k/2$ for every $v\in V(H_1)$. Moreover, at most $\gamma|C_1 -x|$ are lost when passing to $H_1$, and so $x$ clearly still has at least $(1+\theta)4k/3$ neighbours in $H_1$.
    By \cref{fact:cutdense_subgraph}, $H_1$ is $(\rho/2 - 2\gamma-10\sqrt{\gamma})$-cut-dense, and thus it follows from \cref{fact:cutdense_reduced_connected} that $R$ is connected.
    Note that the conclusions of \cref{thm:BPS_combined} must hold applied to $G_{\ref{thm:BPS_combined}} = G[V(H_1)\cup \{x\}]$ with $\eps_{\ref{thm:BPS_combined}} = \theta$. However, clearly \ref{item:BPS3} cannot hold, as $R$ is connected and thus only has one component. Thus we immediately reach a contradiction.
\end{proofclaim}

Now, we wish to gain some control over where the second neighbours of $x$ lie. We say that a vertex $z\in V(G)$ is \defn{external} if there exists $j\in [m]\setminus \{1,2\}$ such that $z\in L_{\eta k}(C_j)$, and \defn{internal} otherwise. 
\begin{claim}\label{claim:internal_vertices}
    Every internal vertex has at least $(1+\eps/2)k/2$ neighbours in $V(C_1) \cup V(C_2)$.
\end{claim}

\begin{proofclaim}
    Suppose $z\in V(G)$ is internal, that is, there is no $j\in [3,m]$ for which $z\in L_{\eta k}(C_j)$. We can assume $|N_G(z)\setminus (C_1\cup C_2)|>\eta k$ as otherwise we are done. Then by \eqref{item:SN2} there exists $j\notin \{1,2\}$ and $i\in [m]$ for which $|N_G(z)\setminus (C_i\cup C_j)|<\eta k$. Since $z$ is internal, this implies $i\in \{1,2\}$. In particular,
        $|N_G(z)\setminus C_i| = |N_G(z)\setminus (C_i\cup C_j)| + |N_G(z)\cap C_j| <2\eta k,$
    from which the claim is clear by noting $|N_G(z)\cap C_i|\geq \delta(G) - 2\eta k\geq (1+\eps/2)k/2$.
\end{proofclaim}
From now on, let \defn{$\ext(x)$} denote the set of vertices in $N^2_G(x)$ that are external, and let \defn{$\inter(x)$} denote the set of vertices in $N_G^2(x)$ that are internal. 
\begin{claim}\label{claim:poor_nbhrs}
    Every vertex in $N_G(x)$ has fewer than $\Delta$ neighbours in $\ext(x)$, and furthermore, every vertex in $N_G(x)$ is internal. 
\end{claim}

\begin{proofclaim}
    Suppose to the contrary that there is a vertex $y\in N_G(x)$ and a set $Z \subseteq N_G(y)\cap \ext(x)$ of size $\Delta$. 
    Apply \cref{cor:partition_subforests_T-v}\ref{item:Fsize2} to $T$ to obtain $r\in V(T)$ such that $T - r$ has a partition into three subforests $F_1$, $F_2$ and $F_3$ each containing at most $\ceil{k/2}$ vertices and such that $F_3$ is a tree. For both $i\in\{1,2\}$, let $G_i = G[(C_i \cup \{x\})\setminus (Z\cup \{y\})]$. Then $\delta(G_i -x)\geq \ceil{k/2}$ and there exists an embedding $\phi_i$ of $F_i \cup \{r\}$ into $G_i$ with $r$ mapped to $x$, using \cref{fact:greedy_embedding} and since by \cref{claim3:balanced_neighbours} we have $d_{G_i}(x)\geq \eps k/6 - 1 - \Delta \geq \Delta$. 
    
    If $F_3 =\emptyset$ then $T$ embeds in $G$, so assume otherwise. Since $F_3$ is a tree, there is a unique vertex $u$ in the set $V(F_3)\cap N_T(r)$, and we consider $F_3$ to be rooted at $u$. Let $v_1,\dots,v_d$ denote the children of $u$ in $F_3$, where $d<\Delta$, and let $W_j$ denote the children of $v_j$ for each $j\in [d]$, noting again that $|W_j|<\Delta$. For each $z\in Z$ choose an arbitrary $i(z)\in [3,m]$ for which $z\in L_{\eta k}(C_{i(z)})$, which exists because $z$ is external. Construct an embedding $\phi_3$ of $F_3$ as follows. Embed $u$ at $y$ and $v_1,\dots,v_d$ into $Z$ arbitrarily. For each $j\in [d]$ in succession, letting $z_j= \phi_3(v_j)$, we can greedily embed $W_j$ into $C_{i(z_j)}\setminus Z$ since the set 
    $N_G(z_j) \cap V(C_{i(z_j)}) \setminus (Z\cup \phi_3(\bigcup_{j'<j}W_{j'}))$ has size at least $\eta k - \Delta^2 \geq \Delta$. Now since $\delta(C_{i(z)})\geq (1/2+\eps/8)k \geq |F_3|$ for each $z\in Z$, the remainder of the tree $F_3$ fits greedily in $\bigcup_{z\in Z}C_{i(z)}$.
    
     Since $\phi_1(F_1) \subset V(G_1 -x)$, $ \phi_2(F_2) \subset V(G_2 - x)$ and $\phi_3(F_3)\subset \bigcup_{i=3}^m V(C_i) \cup Z\cup \{y\}$ are pairwise disjoint, and the image of $x$ under all three embeddings is $r$, combining these yields an embedding of $T$ in $G$, a contradiction. This proves the first part of the claim.

     It is now easy to see that every $y\in N_G(x)$ is internal. Indeed, if there exists $j\in [3,m]$ such that $y\in L_{\eta k}(C_j)$, then all vertices in $N_G(y)\cap V(C_j)$ are in $\ext(x)$ due to the minimum degree condition in $C_j$, and this contradicts the first part of this claim.
\end{proofclaim}
Now, let us define an \defn{$x$-peripheral-matching} to be a matching $M$ in $G$ where each edge has one endpoint in $N_{G}(x)$ and one endpoint in $\ext(x)$, and such that there exists an injection $f:V(M)\cap \ext(x) \rightarrow [3,m]$ satisfying $w\in L_{\eta k}(C_{f(w)})$ for every $w\in V(M)\cap \ext(x)$.

\begin{claim}\label{claim:matching_outside}
    There exists an $x$-peripheral-matching $M$ with at least $D$ edges.
\end{claim}

\begin{proofclaim}
    For every $w\in \ext(x)$, note that by \eqref{item:SN2}, there are at most two distinct $j\in [m]$ for which $w\in L_{\eta k}(C_{j})$. Let $J = \{j\in [3,m]: L_{\eta k}(C_j) \cap \ext(x) \neq \emptyset\}$.
    
    First suppose that $|J|\leq 2\Delta D$. Recall that $x$ has at most $\eta k$ neighbours outside of $V(C_1)\cup V(C_2)$. Let $A \subseteq N^2_G(x)$ be chosen arbitrarily such that $(1+\eps)4k/3 \leq |A| < 2k$. Let $$\hat{G} =G\biggl[\{x\}\cup N_G(x) \cup A \cup \bigcup_{j\in J \cup \{1,2\}}V(C_j)\biggr],$$ noting that $|\hat{G}|\leq 1+ \eta k+2k + (|J|+2)100k \leq 400\Delta Dk <4\eps^{-1}k$. Observe that $|N^2_{\hat{G}}(x)|\geq |A| \geq (1+\eps)4k/3$ and furthermore, $|N_{\hat{G}}(x)|\geq d_G(x)- \eta k\geq  (1+\eps/2)4k/3$.

As for the minimum degree of $\hat{G}$, note that for vertices $v\in V(C_j)$ for some $j\in J\cup \{1,2\}$, we have $d_{\hat{G}}(v)\geq \delta(C_j) \geq (1+\eps/4)k/2$. Using \cref{claim:internal_vertices,claim:poor_nbhrs}, all vertices in $N_G(x)\cup \inter(x)$ have at least $(1+\eps/4)k/2$ neighbours in $V(\hat{G})$. Finally, our choice of $J$ and \cref{item:SN2} ensure that every $w\in \ext(x)$ also has at least $(1+\eps/4)k/2$ neighbours in $V(\hat{G})$. It follows that $\delta(\hat{G})\geq (1+\eps/4)k/2$, and we can apply \cref{thm:secondneighbourhood_dense} to $\hat{G}$ to obtain a copy of $T$ in $\hat{G}$, a contradiction.

Thus we may assume that $|J|>2\Delta D$. 
Consider an auxiliary bipartite graph $B_1$ with vertex classes $\ext(x)$ and $J$, such that there is an edge between $w\in \ext(x)$ and $j\in J$ if and only if $w\in L_{\eta k}(C_j)$. By \eqref{item:SN2} there are at most two rich subgraphs in the collection in which a vertex can have at least $\eta k$ neighbours, so $d_{B_1}(w)\leq 2$ for every $w\in \ext(x)$. On the other hand, $d_{B_1}(j)\geq 1$ for all $j\in J$ by choice of $J$. Applying \cref{fact:bipartite_matching_size} with
$X_{\ref{fact:bipartite_matching_size}} = \ext(x)$, $Y_{\ref{fact:bipartite_matching_size}} = J$ and $d_{\ref{fact:bipartite_matching_size}} = 2$,
we obtain a matching $M_1$ in $B_1$ of size at least $|J|/2$. Define $W = \ext(x)\cap V(M_1)$ so that $|W|\geq |J|/2$ and there exists an injection $f:W \rightarrow [3,m]$ such that $w\in L_{\eta k}(C_{f(w)})$ for each $w\in W$. 

Now define another graph $B_2$ to be the subgraph of $G$ with vertex classes $N_G(x)$ and $W$ and consisting only of edges that go between these two classes. Note that $B_2$ is bipartite because $N_G(x)\cap \ext(x) = \emptyset$ by \cref{claim:poor_nbhrs}. This claim also tells us that $d_{B_2}(y)\leq \Delta$ for every $y\in N_G(x)$. Conversely $d_{B_2}(w)\geq 1$ for every $w\in W$ since $w$ is a second neighbour of $x$. Applying \cref{fact:bipartite_matching_size} with $X_{\ref{fact:bipartite_matching_size}} = N_G(x)$, $Y_{\ref{fact:bipartite_matching_size}} = W$ and $d_{\ref{fact:bipartite_matching_size}} = \Delta$, we obtain a matching $M$ in $B_2$ of size at least $|W|/\Delta \geq |J|/2\Delta \geq D$. By restricting $f$ to the domain $V(M)\cap \ext(x)$, we see that $M$ is an $x$-peripheral-matching with the desired properties. 
\end{proofclaim}

Let $M$ be as in \cref{claim:matching_outside} and restrict it if necessary so that $|M| =D$. Let $Z_M = V(M)\cap \ext(x)$ and $P= V(M)\cap N_G(x)$. By \cref{claim:poor_nbhrs}, these sets are disjoint, and hence they both have size $D$. Let $G^* = G[\{x\}\cup N_G(x)\cup V(C_1) \cup V(C_2)]$. By \cref{claim:internal_vertices,claim:poor_nbhrs}, $\delta(G^*)\geq (1+\eps/4)k/2$, and further note that $|G^*|\leq 1+ \eta k + 100k+100k < 250k$. Let $G^{**} = G^* -Z_M$.
Since $D$ is independent of $k$, we have $\delta(G^{**})\geq (1+\eps/8)k/2$. Note that $P\subseteq V(G^{**})$ because $P\subseteq N_G(x)\setminus Z_M$.

Apply \cref{lem:treesplit_MSF} to partition $T$ into a matching $N$, a subtree $S$ and a subforest $F$ satisfying the properties \ref{item:MSF1} through \ref{item:MSF6}. Let $Q = V(N)\cap V(S)$. By \ref{item:MSF6}, the smallest subtree of $S$ containing $Q$ has size at most $\Delta^{4\Delta+1}$, and by \ref{item:MSF5}, all vertices in $Q$ are in the same bipartition class in $T$. We have $|N|=|Q|\leq \Delta^{4\Delta+1}\leq D$ by assumption. The following process for embedding $T$ using this structure is depicted in \cref{fig:SNmatching}.

\begin{figure}[tb]
    \centering
\tikzset{every picture/.style={line width=0.75pt}} 
\begin{tikzpicture}[x=0.75pt,y=0.75pt,yscale=-1,xscale=1]

\draw  [fill={rgb, 255:red, 255; green, 255; blue, 255 }  ,fill opacity=1 ] (137.5,95.12) .. controls (137.5,66.61) and (173.18,43.5) .. (217.2,43.5) .. controls (261.21,43.5) and (296.89,66.61) .. (296.89,95.12) .. controls (296.89,123.63) and (261.21,146.75) .. (217.2,146.75) .. controls (173.18,146.75) and (137.5,123.63) .. (137.5,95.12) -- cycle ;
\draw  [fill={rgb, 255:red, 0; green, 0; blue, 0 }  ,fill opacity=1 ] (221.5,91.5) .. controls (221.5,89.7) and (222.92,88.24) .. (224.66,88.24) .. controls (226.41,88.24) and (227.83,89.7) .. (227.83,91.5) .. controls (227.83,93.3) and (226.41,94.76) .. (224.66,94.76) .. controls (222.92,94.76) and (221.5,93.3) .. (221.5,91.5) -- cycle ;
\draw    (222.5,90.5) -- (253.75,98.54) ;
\draw    (222.5,90.5) -- (246.5,106.5) ;
\draw    (223.66,90.5) -- (240.66,118.5) ;
\draw    (224.66,91.5) -- (263.5,88.5) ;
\draw [color={rgb, 255:red, 17; green, 108; blue, 219 }  ,draw opacity=1 ][line width=2.25]    (282.73,90) -- (338.91,85.01) ;
\draw [color={rgb, 255:red, 17; green, 108; blue, 219 }  ,draw opacity=1 ][line width=2.25]    (266.29,116.33) -- (325.21,131.31) ;
\draw [color={rgb, 255:red, 17; green, 108; blue, 219 }  ,draw opacity=1 ][line width=2.25]    (251.22,130.4) -- (308.77,160.36) ;
\draw   (320.64,130.4) -- (362.8,131.6) -- (355.1,154.57) -- cycle ;
\draw   (305.25,158.76) -- (340.1,167.34) -- (324.6,188.81) -- cycle ;
\draw   (333.43,85.01) -- (373.28,71.29) -- (374.26,95.48) -- cycle ;
\draw    (364.03,80.92) -- (338.91,85.01) ;
\draw    (369.96,87.28) -- (338.91,85.01) ;
\draw    (357.18,137.21) -- (325.21,131.31) ;
\draw    (355.35,144.48) -- (325.66,132.67) ;
\draw    (336.17,176.25) -- (308.77,160.36) ;
\draw    (328.5,181.5) -- (308.77,160.36) ;
\draw  [fill={rgb, 255:red, 255; green, 255; blue, 255 }  ,fill opacity=1 ] (323.84,189.87) .. controls (323.84,177.33) and (334.06,167.17) .. (346.67,167.17) .. controls (359.28,167.17) and (369.51,177.33) .. (369.51,189.87) .. controls (369.51,202.4) and (359.28,212.57) .. (346.67,212.57) .. controls (334.06,212.57) and (323.84,202.4) .. (323.84,189.87) -- cycle ;
\draw  [fill={rgb, 255:red, 255; green, 255; blue, 255 }  ,fill opacity=1 ] (348.5,146.29) .. controls (348.5,133.76) and (358.72,123.59) .. (371.33,123.59) .. controls (383.95,123.59) and (394.17,133.76) .. (394.17,146.29) .. controls (394.17,158.83) and (383.95,168.99) .. (371.33,168.99) .. controls (358.72,168.99) and (348.5,158.83) .. (348.5,146.29) -- cycle ;
\draw  [fill={rgb, 255:red, 255; green, 255; blue, 255 }  ,fill opacity=1 ] (364.03,80.92) .. controls (364.03,68.39) and (374.25,58.23) .. (386.86,58.23) .. controls (399.47,58.23) and (409.7,68.39) .. (409.7,80.92) .. controls (409.7,93.46) and (399.47,103.62) .. (386.86,103.62) .. controls (374.25,103.62) and (364.03,93.46) .. (364.03,80.92) -- cycle ;
\draw    (419.2,60.04) .. controls (454.83,72.75) and (405.5,130.86) .. (422.86,153.55) ;
\draw    (361.66,221.19) .. controls (398.19,234.81) and (400.02,151.74) .. (422.86,153.55) ;
\draw [color={rgb, 255:red, 128; green, 128; blue, 128 }  ,draw opacity=1 ]   (110.07,150.82) -- (161.93,110.74) ;
\draw [shift={(163.5,109.5)}, rotate = 141.59] [color={rgb, 255:red, 128; green, 128; blue, 128 }  ,draw opacity=1 ][line width=0.75]    (10.93,-3.29) .. controls (6.95,-1.4) and (3.31,-0.3) .. (0,0) .. controls (3.31,0.3) and (6.95,1.4) .. (10.93,3.29)   ;
\draw [color={rgb, 255:red, 128; green, 128; blue, 128 }  ,draw opacity=1 ]   (211.03,173.07) -- (250.02,121.11) ;
\draw [shift={(251.22,119.51)}, rotate = 126.88] [color={rgb, 255:red, 128; green, 128; blue, 128 }  ,draw opacity=1 ][line width=0.75]    (10.93,-3.29) .. controls (6.95,-1.4) and (3.31,-0.3) .. (0,0) .. controls (3.31,0.3) and (6.95,1.4) .. (10.93,3.29)   ;
\draw [color={rgb, 255:red, 155; green, 155; blue, 155 }  ,draw opacity=1 ]   (262.18,180.34) -- (271.54,154.98) ;
\draw [shift={(272.23,153.1)}, rotate = 110.25] [color={rgb, 255:red, 155; green, 155; blue, 155 }  ,draw opacity=1 ][line width=0.75]    (10.93,-3.29) .. controls (6.95,-1.4) and (3.31,-0.3) .. (0,0) .. controls (3.31,0.3) and (6.95,1.4) .. (10.93,3.29)   ;
\draw [color={rgb, 255:red, 155; green, 155; blue, 155 }  ,draw opacity=1 ][fill={rgb, 255:red, 128; green, 128; blue, 128 }  ,fill opacity=1 ]   (298.72,217.56) -- (321.71,202.32) ;
\draw [shift={(323.38,201.22)}, rotate = 146.47] [color={rgb, 255:red, 155; green, 155; blue, 155 }  ,draw opacity=1 ][line width=0.75]    (10.93,-3.29) .. controls (6.95,-1.4) and (3.31,-0.3) .. (0,0) .. controls (3.31,0.3) and (6.95,1.4) .. (10.93,3.29)   ;
\draw [fill=black  ,fill opacity=1 ] 
(332.91,85.01) .. controls (332.91,83.35) and (334.25,82.01) .. (335.91,82.01) .. controls (337.56,82.01) and (338.91,83.35) .. (338.91,85.01) .. controls (338.91,86.67) and (337.56,88.01) .. (335.91,88.01) .. controls (334.25,88.01) and (332.91,86.67) .. (332.91,85.01) -- cycle ;
\draw  [color=black  ,draw opacity=1]
(332.91,85.01) .. controls (332.91,83.35) and (334.25,82.01) .. (335.91,82.01) .. controls (337.56,82.01) and (338.91,83.35) .. (338.91,85.01) .. controls (338.91,86.67) and (337.56,88.01) .. (335.91,88.01) .. controls (334.25,88.01) and (332.91,86.67) .. (332.91,85.01) -- cycle ;
\draw  [fill={rgb, 255:red, 0; green, 0; blue, 0 }  ,fill opacity=1 ] (319.21,131.31) .. controls (319.21,129.65) and (320.55,128.31) .. (322.21,128.31) .. controls (323.86,128.31) and (325.21,129.65) .. (325.21,131.31) .. controls (325.21,132.97) and (323.86,134.31) .. (322.21,134.31) .. controls (320.55,134.31) and (319.21,132.97) .. (319.21,131.31) -- cycle ;
\draw  [fill={rgb, 255:red, 0; green, 0; blue, 0 }  ,fill opacity=1 ] (302.77,159.36) .. controls (302.77,157.71) and (304.11,156.36) .. (305.77,156.36) .. controls (307.42,156.36) and (308.77,157.71) .. (308.77,159.36) .. controls (308.77,161.02) and (307.42,162.36) .. (305.77,162.36) .. controls (304.11,162.36) and (302.77,161.02) .. (302.77,159.36) -- cycle ;
\draw  [color=black  ,draw opacity=1 ][line width=1]  (253.75,98.54) .. controls (240.85,113.67) and (235.08,129.89) .. (240.88,134.77) .. controls (246.67,139.65) and (261.82,131.34) .. (274.72,116.21) .. controls (287.62,101.07) and (293.38,84.85) .. (287.59,79.98) .. controls (281.8,75.1) and (266.65,83.41) .. (253.75,98.54) -- cycle ;
\draw [color={rgb, 255:red, 155; green, 155; blue, 155 }  ,draw opacity=1 ]   (362.5,47) -- (347.96,61.64) ;
\draw [shift={(346.55,63.06)}, rotate = 314.82] [color={rgb, 255:red, 155; green, 155; blue, 155 }  ,draw opacity=1 ][line width=0.75]    (10.93,-3.29) .. controls (6.95,-1.4) and (3.31,-0.3) .. (0,0) .. controls (3.31,0.3) and (6.95,1.4) .. (10.93,3.29)   ;

\draw (158.6,82.73) node [anchor=north west][inner sep=0.75pt]   [align=left] {$\displaystyle G^{**} \ \ $};
\draw (208.78,80.35) node [anchor=north west][inner sep=0.75pt]   [align=left] {$\displaystyle x$};
\draw (256.97,95.9) node [anchor=north west][inner sep=0.75pt]  [opacity=1 ] [align=left] {$\displaystyle P$};
\draw (303.14,63.71) node [anchor=north west][inner sep=0.75pt]   [align=left] {$\displaystyle \textcolor[rgb]{0.07,0.42,0.86}{M}$};
\draw (373.85,71.23) node [anchor=north west][inner sep=0.75pt]   [align=left] {$\displaystyle C_{j_{F'}}$};
\draw (428.6,143.65) node [anchor=north west][inner sep=0.75pt]   [align=left] {$\displaystyle \geqslant D$ rich subgraphs};
\draw (91.44,151.73) node [anchor=north west][inner sep=0.75pt]  [color={rgb, 255:red, 128; green, 128; blue, 128 }  ,opacity=1 ] [align=left] {$\displaystyle S$};
\draw (192.57,170.88) node [anchor=north west][inner sep=0.75pt]  [color={rgb, 255:red, 128; green, 128; blue, 128 }  ,opacity=1 ] [align=left] {$\displaystyle \mathnormal{Q}$};
\draw (250,183.78) node [anchor=north west][inner sep=0.75pt]  [color={rgb, 255:red, 128; green, 128; blue, 128 }  ,opacity=1 ] [align=left] {$\displaystyle N$};
\draw (188.24,209.01) node [anchor=north west][inner sep=0.75pt]  [color={rgb, 255:red, 128; green, 128; blue, 128 }  ,opacity=1 ] [align=left] {{\fontfamily{ptm}\selectfont \textcolor[rgb]{0.61,0.61,0.61}{Components of }}\textcolor[rgb]{0.61,0.61,0.61}{$\displaystyle \mathnormal{F}$}};
\draw (368.83,122.31) node [anchor=north west][inner sep=0.75pt]  [rotate=-282.73] [align=left] {$\displaystyle \dotsc $};
\draw (302.15,117.95) node [anchor=north west][inner sep=0.75pt]  [color={rgb, 255:red, 17; green, 108; blue, 219 }  ,opacity=1 ,rotate=-282.73] [align=left] {$\displaystyle \dotsc $};
\draw (332,67) node [anchor=north west][inner sep=0.75pt] 
[align=left] {$\displaystyle z_{F'}$};
\draw (358.94,32.83) node [anchor=north west][inner sep=0.75pt]  [color={rgb, 255:red, 128; green, 128; blue, 128 }  ,opacity=1 ] [align=left] {$\displaystyle u_{F'}$};

\end{tikzpicture}
    \caption{\centering Embedding process using the $x$-peripheral-matching $M$.}
    \label{fig:SNmatching}
\end{figure}

By \ref{item:MSF3}, $|S|\leq \ceil{k/2}$. Thus we can apply \cref{lem:embedfor_matchingpoints} with $\Delta$, $ \eps/8$, $G^{**}$, $S$, $P$, $Q$ playing the roles of $\Delta$, $\eps$, $G$, $T$, $P$, $Q$ respectively to find an embedding $\psi:S\hookrightarrow G^{**}$ such that $\psi(v)\in P$ for all $v\in Q$. That is, there is a copy of $S$ in $G^{**}$ such that the endpoints of edges of $N$ that intersect $S$ are embedded into $P = V(M)\cap V(G^{**})$.

Each vertex of $P$ lies on a unique edge of the host matching $M$. For every edge of the tree matching $N$, extend $\psi$ across the unique edge of $M$ incident with the image of its endpoint in $S$. Because $\psi$ is injective and $M$ is a matching, these extensions are pairwise disjoint. Thus, for each component $F'\subset F$, its unique vertex $u_{F'}\in V(N)\cap V(F)$ is now embedded at a vertex $z_{F'}\in Z_M$. Let $j_{F'} = f(z_{F'})$ so that $z_{F'}\in L_{\eta k}(C_{j_{F'}})$. Since $f$ is injective, the indices $j_{F'}$ are distinct.

It remains to embed the rest of each component $F'$. The image of $S$ lies in $G^{**}\subseteq G^*$ and hence, outside $C_1\cup C_2$, can only use $x$ and vertices of $N_G(x)$. By~\cref{claim:poor_nbhrs}, no vertex of $N_G(x)$ belongs to any $C_j$ with $j\geq 3$. For each $F'$, therefore use the host
\begin{equation*}
    H_{F'} = G[(V(C_{j_{F'}}) \setminus (\{x\}\cup Z_M)) \cup \{z_{F'}\}].
\end{equation*}

The only vertices removed from $C_{j_{F'}}$ are $x$ and the $D = O_{\Delta}(1)$ other matching endpoints. So, we have $\delta(H_{F'} -z_{F'})\geq (1/2+\eps/16)k\geq |F'|$ using~\ref{item:MSF3}, whilst $z_{F'}$ has at least $\eta k/2$ neighbours in $H_{F'}$ using the definition of $f$. By~\cref{fact:greedy_embedding} we can greedily embed each component $F'\subseteq F$ into its own $H_{F'}$ with $u_{F'}$ rooted at $z_{F'}$. These extensions avoid the image of $S$ and all matching endpoints already used by other components. This completes the embedding of $T$ in $G$ and therefore concludes the proof of the theorem.
\end{proof}

\section{Conclusion}\label{sec:conclusion}

\textbf{Extensions to higher degree trees.}
Currently, the methods used in our overarching proof rely fairly heavily on $k$ being sufficiently large with respect to the constant $\Delta$. Naturally, one would like to push the ideas further to account for trees with higher maximum degrees. Since \cref{thm:CoverV2} is applied as a black box early in our proof, the only hope of doing this would be to first find an extension of this theorem to a larger family of trees. As discussed in more detail in \cite{AP_hyperstability}, there are examples showing that an analogue of \cref{thm:CoverV2} cannot be true if we consider trees of linear maximum degree, so perhaps it is more reasonable to see what can be done for trees with sublinear maximum degree. 

\textbf{Tightness of degree combinations.}
As mentioned in the introduction, the combination of minimum and maximum degree in \cref{conj:alpha} is asymptotically best possible for all $\alpha = 1/\ell$ when $\ell \geq 5$ is odd, seen in \cite{BPSmaxmin2020}. This is witnessed (for $k$ divisible by $\ell$) by the tree $T_{k,\ell}$ consisting of one root vertex with $\ell$ children that each have exactly $k/\ell-1$ distinct children of their own. The host graph $G$ avoiding $T_{k,\ell}$ is composed of the disjoint union of two copies of a complete bipartite graph with parts of size close to $(1+\alpha)k/2$ and $(1-\alpha)k$, and adding a vertex $x$ that has neighbours in the bigger bipartition class in both copies. 
When $\alpha\in (0,1/3)$ is not a unit fraction, the authors of \cite{BPSmaxmin2020} show that the conjecture is at least not far from being optimal by exhibiting a host graph with minimum degree and maximum degree close to $(1+\alpha)k/2$ and $2(1-\alpha-O(\alpha^2))k$ respectively, that does not contain all $k$-edge trees. 
It would be interesting to know the optimal combination of bounds in this regime, and furthermore, the effect of restricting the problem to a family of bounded degree trees, where the example using $T_{k,\ell}$ would not be permissible.  
 
A positive result relating to this, due to Besomi, Pavez-Sign\'e and Stein \cite{BPSmaxmin2020} (discussed further in \cref{subsec:BPS}), proves that every dense graph $G$ with minimum degree asymptotically above $k/2$ and maximum degree asymptotically above $4k/3$ contains a copy of every bounded degree tree, provided $G$ avoids a specific structure. Furthermore, we proved \cref{thm:mainalpha-exact} via a reduction to \cref{thm:alpha_weakened}, which assumes a slightly lower bound on the maximum degree at the cost of assuming a slightly larger minimum degree. 

\textbf{Further exact results.}
We remark that there seems to be more difficulty in proving exactly (non-asymptotically) the conjectures that only require the minimum degree of the host graph to be at least $k/2$, and not higher, for example \cref{conj:boundedk/2,conj:secondneighbourhood}. A notable exception is the result of Hyde and Reed \cite{HYDE2023}, proving that there exists a function $g(k)$ such that every graph with minimum degree at least $k/2$ and maximum degree at least $g(k)$ contains a copy of every $k$-edge tree. One reason for the complexity is that there are several different structures a host graph may take with minimum degree just below $k/2$ and high maximum degree, that do not embed all $k$-edge trees. The most straightforward instance is given by the complete bipartite graph $K_{m,n}$ where $m= \floor{k/2}-1$ and $n$ is arbitrarily large, which does not contain any balanced $k$-edge tree. Alternatively, if $G$ is the union of $d$ vertex-disjoint cliques each with $\floor{k/2}-1$ vertices, together with an additional vertex $x$ that has exactly one neighbour in each clique, then $G$ can never embed a path on $k$-edges, no matter how large $d$ is, since it would require us to embed vertices in at least three of the cliques, but all vertices in the path have degree at most two. Similar variations of this also exist that forbid more trees. It seems that proving exact versions of these conjectures is challenging, whilst a very interesting direction.

\textbf{Directed graph analogues.}
Several similar questions have been posed in the setting where the host graph is a directed graph, and the fixed size trees are directed/oriented and often satisfy some additional constraints. These types of problems have also received a lot of traction in recent years, see e.g.~\cite{kathapurkar2022spanning,KLIMOSOVA2023113515,grzesik2025antidirected,Stein_2024,stein_trujillo2025,stein_zarate2024,skokan2024alternating,chen2025long}. 
Closely related to the results in this paper, one can consider the minimum semi-degree of a directed graph $D$, to be the minimum of $\delta^+(D)$ and $\delta^-(D)$. 
Recall that \cref{conj:2k/3} was proven for spanning trees, i.e., in the case where $|G| = k+1$, by Reed and Stein \cite{ReedStein_spanning1,ReedStein_spanning2}. In a similar spirit, Stein \cite[Question 6.12]{Stein_2024} asked whether it is true that every sufficiently large $n$-vertex digraph of minimum semi-degree exceeding $2n/3$ that has a vertex $x$ satisfying $d^+(x) = d^-(x) = n-1$ contains every oriented tree on $n$ vertices.

More generally, Stein \cite[Problem 6.13]{Stein_2024} also posed the following problem: Determine the smallest $f(k)$ such that for every $k\in \N$, every oriented graph (or digraph) of minimum semi-degree exceeding $f(k)$ that has a vertex $x$ with both out-degree and in-degree at least $k$ contains each oriented $k$-edge tree. We must have $f(k) \leq k$ as if the minimum semi-degree is at least $k$, a greedy argument is sufficient. On the other hand,
a directed example obtained by modifying the extremal example in \cref{fig:extremal2k/3} shows that we must have $f(k)\geq 2k/3$. It would be very interesting to know if this lower bound is attained. As evidence towards this, Kontogeorgiou, Santos and Stein \cite{kontogeorgiou2025antidirected} recently showed that every large dense digraph of minimum semi-degree
asymptotically above $ 2k/3$ that has a vertex of out-degree and in-degree asymptotically above~$k$ contains every balanced anti-directed $k$-edge tree with bounded maximum degree. It seems plausible that the bounds are identical to those in the graph case, analogous to \cref{conj:2k/3}.

\section*{Acknowledgements}
The authors are grateful to the anonymous referees for helpful comments on an earlier version of this manuscript.

\small
\setlength{\baselineskip}{12pt}

\end{document}